\documentclass[11pt,twoside]{article}

\usepackage[utf8]{inputenc} 
\newcommand{\titel}{Algebraic hierarchy of logics unifying fuzzy logic and quantum logic}
\newcommand{\fett}[1]{\mbox{\textbf{\textit{#1}}}}
\usepackage{amsmath}
\usepackage{amssymb}
\usepackage{amsthm}
\usepackage{latexsym}

\usepackage{mathrsfs}
\usepackage{color}

\usepackage{mathptmx,times}

\usepackage{ifpdf}
\ifpdf
  \usepackage[pdftex]{graphicx}
  \usepackage[pdftitle={\titel},
              pdfauthor={Andreas de Vries},%
              pdfpagemode=UseNone,%
              colorlinks=false,hyperindex=true,plainpages=false,%
  ]{hyperref}
  \DeclareGraphicsExtensions{.pdf,.jpg,.png}
\else
  \usepackage{graphicx}
  \usepackage[pdftitle={\titel},
              pdfauthor={Andreas de Vries},%
              pdfpagemode=UseNone,%
              colorlinks=false,hyperindex=true,plainpages=false,%
              dvips=true,ps2pdf=true%
  ]{hyperref}
  \DeclareGraphicsExtensions{.eps}
\fi

\newcommand{\URL}  [1]{\footnotesize\textrm{#1}}

\newcommand{\code} [1]{{\footnotesize\tt #1}}

\newtheorem{satz}{Theorem}
\newtheorem{lem}[satz]{Lemma}
\newtheorem{defi}[satz]{Definition}

\newtheorem{axiom}[satz]{Axiom}
\newtheorem{cor}[satz]{Corollary}
\newtheorem{rem}[satz]{Remark}
\newtheorem{beisp}[satz]{Example}
\newtheorem{beispe}[satz]{Examples}

\newtheorem{exercise}{Problem} 
\newenvironment{definition} {\begin{defi}\begin{em}}{\end{em}\hfill{$\lozenge$}\end{defi}}
\newenvironment{beispiel} {\begin{beisp}\begin{em}}{\end{em}\hfill{$\lozenge$}\end{beisp}}
\newenvironment{beispiele} {\begin{beispe}\begin{em}}{\end{em}\hfill{$\lozenge$}\end{beispe}}

\newenvironment{bem} {\begin{rem}\begin{em}}{\hfill{$\lozenge$}\end{em}\end{rem}}

\newenvironment{pgm*} {\begin{tt}
                        \begin{footnotesize}
                          \begin{tabular}{l}
                      }
                      {    \\
                        \end{tabular}
                        \end{footnotesize}
                       \end{tt}}

\begin{document}

\title{\titel\\[1ex] --- Lecture Notes ---}

\author{Andreas de Vries}%

\author{Andreas de Vries%
\thanks{%
	e-Mail:
	\href{mailto:de-vries@fh-swf.de}{\URL{de-vries@fh-swf.de}},
	URL: 
	\href{http://www3.fh-swf.de/fbtbw/devries/deVries.htm}
	{\URL{http://www3.fh-swf.de/fbtbw/devries/deVries.htm}}%
}
	\\ 
	{\footnotesize 
		South Westphalia University of Applied Sciences, 
		Haldener Stra{\ss}e 182, 
		58095 Hagen, Germany
	}
}

\date{\today}

\maketitle

\begin{abstract}
	\noindent
	In this paper, a short survey about the concepts underlying general logics is given.
	In particular, 	a novel rigorous definition of a fuzzy negation as
	an operation acting on a lattice to render it into a fuzzy logic is presented.
	According to this definition, 
	a fuzzy negation satisfies the weak double negation condition,
	requiring double negation to be expansive, the antitony condition,
	being equivalent to the disjunctive De Morgan law and
	thus warranting compatibility of negation with the lattice operations,
	and the Boolean boundary condition stating that the universal bounds of the lattice
	are the negation of each other.
	From this perspective, the most general logics
	are fuzzy logics, containing as special cases
	paraconsistent (quantum) logics,
	quantum logics, intuitionistic logics, 
	and Boolean logics,
	each of which given by its own algebraic restrictions.
	New examples of a non-contradictory logic violating the conjunctive De Morgan
	law, and of a typical non-orthomodular fuzzy logic along with its explicit lattice
	representation are given.
\end{abstract}


\setcounter{tocdepth}{2}
\tableofcontents

\section{Introduction}
\begin{quote}
\begin{em}
	Logic is as empirical as geometry. 
	We live in a world with a non-classical logic.
\end{em}

\hfill{Hilary Putnam}
\end{quote}



\noindent
Logic is the science which investigates the principles governing
correct or reliable inference. It deals with propositions
and their relations to each other.
Besides the classical Boolean logic there has been established
various generalizations such as modal, intuitionistic, quantum, or fuzzy logic,
as well as propositional structures underlying substructural logics which focus 
on relaxations of structural rules governing
validity and provability.
The purpose of the present article is to give a brief unifying survey of
the algebraic interrelations of logics
which does not seem to have been presented before in this comprehensive form.
It does not intend to deal, however, with the wide aspect of
semantic algebras deriving logics, an important issue in considerations
both of fuzzy logics and of quantum logics.
It concentrates on the propositional structure of general logics
but noteworthy not on logical calculi neither on provability or model theory.
Although this work has been greatly influenced by some classical references
such as \cite{Birkhoff-1973, Mittelstaedt-1978},
it provides a wider spectrum by including fuzzy logics,
and thus complements modern approaches like \cite{Restall-2000} 
in revealing the algebraic hierarchy of logics as
propositional structures of substructural logics,
especially pointing out the decisive role of negation.

The mathematical concept underlying any logic is the notion of the
lattice, a partially ordered set with two binary operations 
forming an algebraic structure.
Establishing a lattice with an additional fuzzy negation operator yields a
fuzzy logic, and further algebraic requirements such as non-contradiction, 
paraconsistency,
orthomodularity, or distributivity then specify it to the different
classical and non-classical logics as sketched in Figure \ref{fig-logics}.
This comprehensive view on the different concepts of logics is
enabled by defining a fuzzy negation as a lattice operation satisfying
weak double negation, antitony
and the Boolean boundary condition.
Remarkably, antitony is equivalent to the \emph{disjunctive} De Morgan law,
but does not imply the conjunctive De Morgan law.
This definition is well established in fuzzy logic contexts 
\cite{Deschrijver-Kerre-2005}
and generalizes commonly used notions of a negation
as an involutive operation 
\cite{Dalla-Chiara-Giuntini-2002,Dvurecenskij-Pulmannova-2000}.
It is wide enough to include all common fuzzy, quantum, intuitionistic,
and classical logics.

Consequently, the next section of is paper starts with an outline of lattice theory
and subsequently gives the definition and important properties of
a fuzzy logic and its common sublogics.
In the following sections, various logics are considered in some more detail,
with emphasis on examples of a typical non-orthomodular fuzzy logic
and some special logics of quantum registers in unentangled states.


\subsection{Notation}
Every investigation, including the present one about logic,
has to be communicated by means of language.
The language being used is usually called the 
\emph{metalanguage}\index{metalanguage}.
It has to be distinguished carefully from the language
of the studied logic, the \emph{object language}\index{object language}
%
There are many different notations existing in the literature, so
Table \ref{tab-logic-functors} lists the symbols as they are used in the
present text.
\begin{table*}[htp]
\centering
\begin{footnotesize}
\begin{tabular}{|l|ll@{\ }l|}
	\hline
	&
	$A \Rightarrow B$ & if $A$ then $B$; $A$ only if $B$
	&
	\\
	&
	$A \Leftrightarrow B$ & $A$ if and only if $B$
	&
	\\
	\textbf{Metalanguage} 
	&
	$x := y$ & $x$ is defined as $y$
	&
	\\
	&
	$x=y$ & $x$ equals $y$
	&
	\\
	&
	$x \leqq y$ & $x$ precedes $y$, $x$ is less than or equal to $y$
	&
	\\ \hline
	&
	$\neg x$ & not $x$ & (negation)
	\\
	&
	$x \wedge y$ & $x$ and $y$ & (meet, conjunction)
	\\
	\textbf{Object language} 
	&
	$x \vee y$ & $x$ or $y$ & (union, disjunction)
	\\
	&
	$x \to y$ & 
	$x$ implies $y$,
	$\neg x \vee y$ 
	& (material implication)
	\\
	&
	$x \leftrightarrow y$ 
	& $x$ is equivalent to $y$, $(x \rightarrow y) \wedge (y \rightarrow x)$ 
	& (equivalence operation)
	\\ \hline
\end{tabular}
\end{footnotesize}
\vspace*{-1ex}
\caption{\label{tab-logic-functors}\footnotesize
	Logical functors of metalanguage and an exemplary object language,
	here Boolean logic.
}
\end{table*}

\section{Lattices}
Lattice theory is concerned with the properties of
a binary relation $\leqq$, to be read ``precedes or equals,''
``is contained in,''
``is a part of,'' or ``is less than or equal to.''
This relation is assumed to have certain properties, the most
basic of which leads to the following concept of a ``partially
ordered set,'' or ``poset.''

\begin{definition}
	A \emph{poset}\index{poset} $(X, \leqq)$
	is a set $X$ in which a binary relation
	$x \leqq y$ is defined, which satisfies,
	for all $x$, $y$, $z\in X$,
	\begin{eqnarray}
		\label{P1}
		\mbox{\emph{(Reflexivity)}} & &
		x \leqq x 
		\\
		\label{P2}
		\mbox{\emph{(Antisymmetry)}} & &
		\mbox{If $x \leqq y$ and $y \leqq x$, then $x=y$. }
		\\
		\label{P3}
		\mbox{\emph{(Transitivity)}} & &
		\mbox{If $x \leqq y$ and $y \leqq z$,
		then $x \leqq z$.}
	\end{eqnarray}
%
%
\end{definition}

If $x \leqq y$ and $x \ne y$, we write $x < y$, and say that $x$ ``precedes,''
``properly contains,'' or ``is less than'' $y$.
The relation $x \leqq y$ is also written $y \geqq x$, and reads
``$y$ succeeds or
contains $x$.'' We often write simply $X$ instead of 
$(X, \leqq)$ and speak of the poset $X$.
Some familiar examples of partially ordered sets are the following.

\begin{beispiele}
\label{bsp-posets}
	(a) $(\mathbb{R}, \leqq)$ is a poset, where $\leqq$ denotes ``less than or equal.''
	
	(b) Let $\mathscr{P}(\Omega)$ be the potential set of a set $\Omega$, i.e.,
	the set of all subsets of $\Omega$ including $\Omega$ itself and the empty
	set $\emptyset$. Then $(\mathscr{P}(\Omega), \subseteq)$ is a poset.
	
	(c) $(\mathbb{N}, \mid)$ is a poset. Here $\mathbb{N}$ is the set of positive
	integers, and $x \mid y$ denotes ``$x$ divides $y$.''
	
	(d) Let $\mathfrak{F}([a,b])$ be the set of all real-valued functions $f(x)$
	on the interval $[a.b]$ (where $a<b$), and let $f \leqq g$ mean that
	$f(x) \leqq g(x)$ for \emph{every} $x$ $\in$ $[a,b]$. Then
	$(\mathfrak{F}([a,b]), \leqq)$ is a poset.
	
	(e) Let $n>1$. Then $(\mathbb{R}^n, \preceq)$ is a poset, where $x \preceq y$ denotes 
	``componentwise less than or equal,'' i.e., $x_j \leqq y_j$ for all 
	$j=1$, \ldots, $n$ where $x=(x_1$, \ldots, $x_n)$ and
	$y=(y_1$, \ldots, $y_n)$. Note that there exist many points $x$, $y\in\mathbb{R}^n$
	such that neither $x \preceq y$ nor $y \preceq x$ holds, for instance for
	$x=(1,2)$ and $y=(0,3)$.
\end{beispiele}

\begin{figure}[htp]
\centering
\begin{footnotesize}
\begin{tabular}{*{6}{c@{\quad }}}
	\includegraphics[height=12ex]{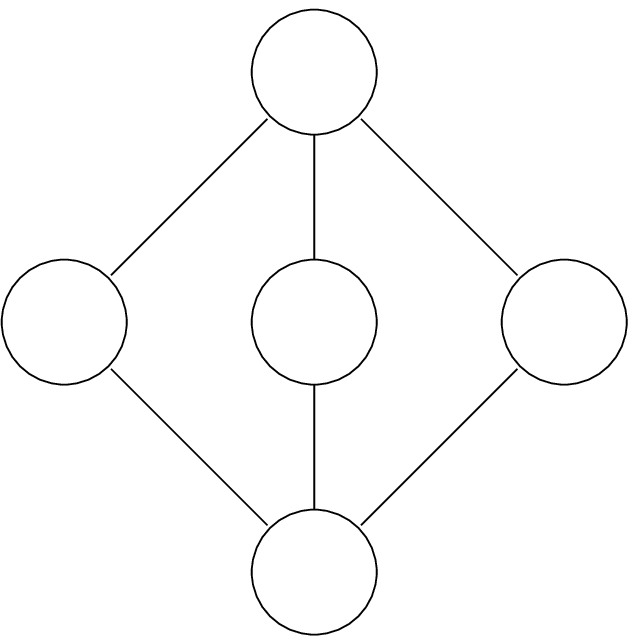}  &
	\includegraphics[height=12ex]{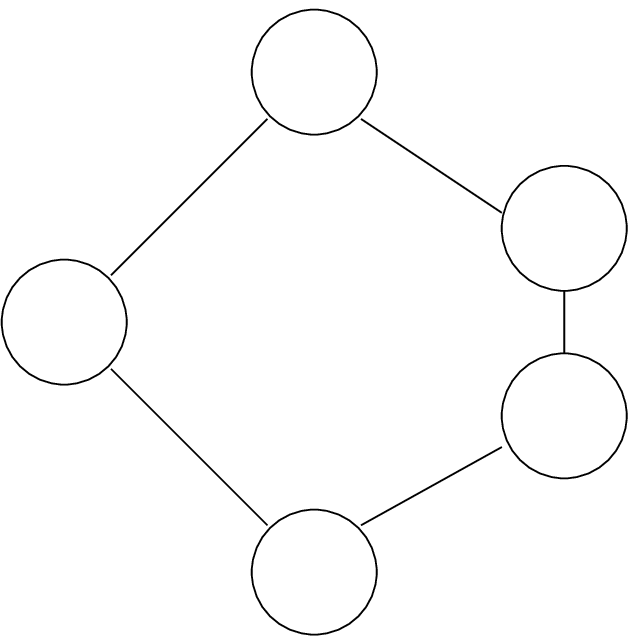}  &
	\includegraphics[height=12ex]{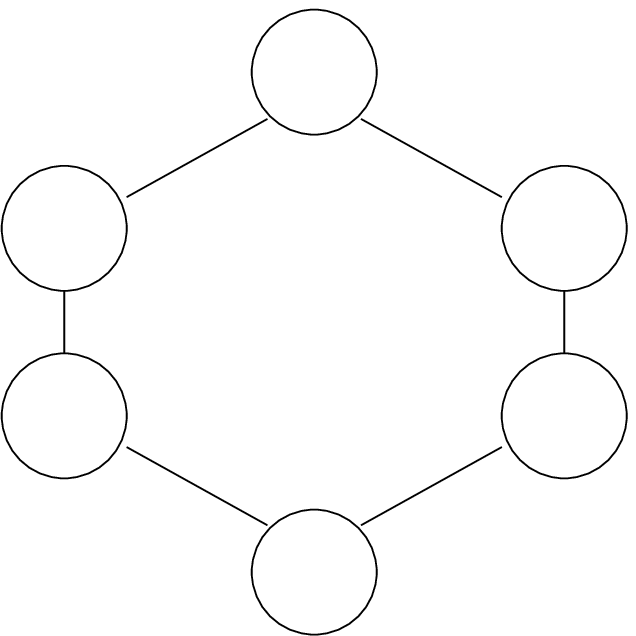}  &
	\includegraphics[height=12ex]{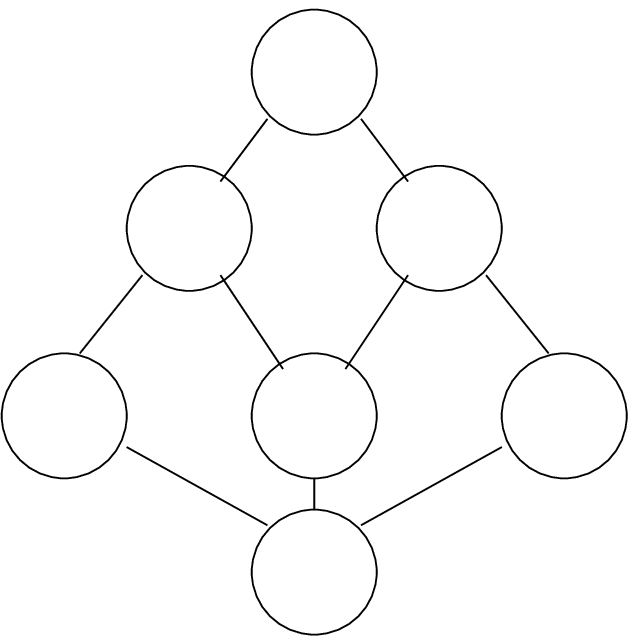}  &
	\includegraphics[height=16ex]{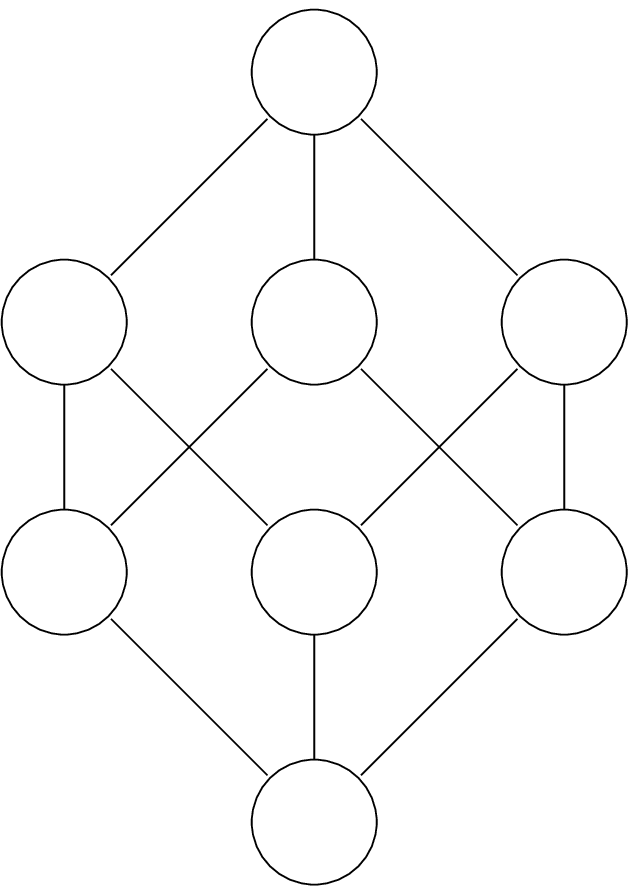} &
	\includegraphics[height=12ex]{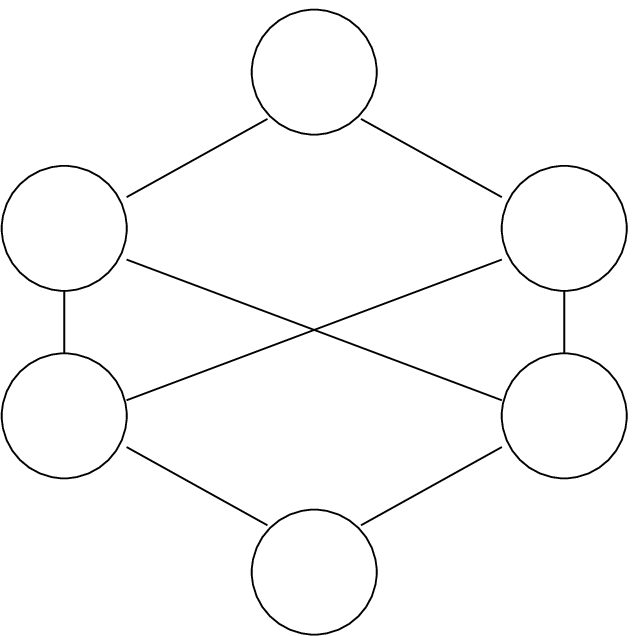}
	\\
	$M_5$ & $N_5$ & $O_6$ & $L_7$ & $\mathbf{2}^3$ & $P_6$
\end{tabular}
\end{footnotesize}
\caption{\label{fig-lattices}\footnotesize
	Hasse diagrams\index{Hasse diagram} of various posets.
	All depicted posets unless $P_6$ are lattices.%
	\index{M$_5$}\index{N$_5$}\index{O$_6$}\index{L$_7$}\index{Z$_2^3$}\index{P$_6$}%
	\index{$M_5$}\index{$N_5$}\index{$L_7$}\index{$P_6$}%
	\index{$\mathbf{2}^3$}%
	}
\end{figure}

A common way to depict a poset with finitely many elements is by means of a 
Hasse diagram, cf.\ Fig.\ \ref{fig-lattices}.
Here each element of the poset is represented by a circle,
placing $y$ higher than $x$ whenever $x$ $<$ $y$.
If $x<y$ and there exists no $z\in X$ such that $x<z<y$, then a
straight line is drawn between $x$ and $y$. (In this case, one often says
``$y$ covers $x$'').

A poset $X$ can contain at most one element $O \in X$ which satisfies 
$O \leqq x$ for all $x$ $\in$ $X$. For if $O$ and $\tilde O$ are two such elements,
then $O \leqq \tilde O$ and also $\tilde O \leqq O$, i.e., $O=\tilde O$ by (\ref{P2}).
If such an element $O$ exists, it is called the 
\emph{least element}\index{least element of a poset} of $X$.
With the analogous reasoning, the 
\emph{greatest element}\index{greatest element of a poset}
of $X$, if it exists, is denoted by $I$ and satisfies
$x \leqq I$ for all $x \in X$.
If both $O$ and $I$ $\in$ $X$ exist, they are called 
\emph{universal bounds}\index{universal bounds of a poset}
of $X$, since then $O \leqq x \leqq I$ for all $x$ $\in$ $X$.
Such an element $O$ is also called the \emph{zero}\index{zero of a poset}
of the poset $L$, and $I$ the \emph{unit}\index{unit of a poset} of the
poset.
In connection with logics, being special posets as we shall see below,
$O$ is also called the ``absurdity.'' 

In a poset $X$ with a least element $O\in X$,
elements $x$ $\in$ $X$ satisfying
$O < x$ such that there is no $y\in X$ with
$O < y < x$ are called \emph{atoms}\index{atom} or
\emph{points} of $X$.
In logics we thus can state that 
``atoms immediately follow from absurdity.''

An \emph{upper bound}\index{upper bound} of a subset 
$Y \subseteq X$ of a poset $X$ is an element $a\in X$
with $y \leqq a$ for every $y\in Y$.
The \emph{least upper bound}\index{sup} $\sup Y$ is an upper bound
contained in every other upper bound. By (\ref{P2}), $\sup Y$ is
unique if it exists.
Note that $Y=P_6 \setminus \{I\}$ in Figure \ref{fig-lattices}
does not contain an upper bound.
The notions of \emph{lower bound}\index{lower bound} and 
\emph{greatest lower bound}\index{inf}
$\inf Y$ are defined analogously. Again by (\ref{P2}), $\inf Y$ is unique if it exists.

\begin{definition}
	A \emph{lattice}\index{lattice} is a poset $L$ such that any elements
	$x$, $y$ $\in$ $L$ 
	have a unique greatest lower bound, denoted $x \wedge y$, 
	and a unique least upper bound $ x \vee y$, i.e.,
	\begin{equation}
		\label{def-wedge-vee}
		x \wedge y = \inf \{x, y\},
		\qquad
		x \vee y = \sup \{x, y\}.
	\end{equation}
	The operation $\wedge$\index{$\wedge$} 
	is also called \emph{meet}\index{meet operation},
	and the operation $\vee$\index{$\vee$} is called \emph{join}\index{join operation}.
	A lattice $L$ is \emph{complete}\index{complete}, 
	when for any set $D \subseteq L$ the bounds 
	$\sup D$ and $\inf D$ exist;
	it is \emph{$\sigma$-complete}\index{sigma-complete}, 
	when for any countable set $D \subseteq L$ the bounds 
	$\sup D$ and $\inf D$ exist.
	A lattice is \emph{atomic} if every element is a join of atoms.
\end{definition}

The lattice condition will later ensure that the logical operations of 
conjunction\index{conjunction} ($\wedge$) and disjunction\index{disjunction} ($\vee$)
are well-defined for any pairs of propositions.
For instance, 
the two atoms of the poset $P_6$ in
Figure \ref{fig-lattices} do not have unique least upper bounds,
hence $P_6$\index{P$_6$}\index{$P_6$} is not a lattice.
Moreover, it follows that finite sets of pairwise disjoint
(``orthogonal'')
propositions also have a well-defined disjunction.
The completeness condition ensures that the latter is true also for 
countable sets of pairwise disjoint propositions. This is arguably not an 
essential requirement for a logic (must logic be infinitary?), 
but it allows for probability measures to be defined on an infinite 
lattice or poset, since it is customary to require that the probability of a 
countable set of disjoint events be well-defined and equal to the countable 
sum of the probabilities of the disjoint events 
($\sigma$-additivity of probabilities).

\begin{beispiel}
\label{bsp-lattices}
	The poset $(\mathscr{P}(\Omega), \subseteq)$ in Example \ref{bsp-posets}(b)
	is a lattice, where for any family
	$\mathscr{A}$ $=$ \{$A_1$, $A_2$, \ldots \} of subsets 
	$A_1$, $A_2$, \dots $\subseteq$ $\Omega$ we have
	\begin{equation}
		\inf \mathscr{A} = \bigcap_{j} A_j,
		\quad
		\sup \mathscr{A} = \bigcup_{j} A_j.
	\end{equation}
	Note that $\mathscr{A}$ $\subseteq$ $\mathscr{P}(\Omega)$.
	Especially, we have $A \wedge B$ $=$ $A \cap B$ and $A \vee B$ $=$ $A\cup B$
	for two subsets $A$, $B$ $\subseteq$ $\Omega$.
	Especially, $\mathscr{P}(\Omega)$ is a complete lattice.
\end{beispiel}

\begin{beispiel}
\label{bsp-vector-space}
	For a field $\mathbb{K}$, 
	let
	$L(\mathbb{K}^n)$ $=$ \{$V$ $\subseteq$ $\mathbb{K}^n$: $V$ is a vector space\}
	be the set of all subspaces of the vector space $\mathbb{K}^n$.
	Moreover, for any set $A\subseteq \mathbb{K}^n$ let span\,$A$ denote the intersection
	of all subspaces of $\mathbb{K}^n$ which contain $A$, i.e.,
	\begin{equation}
		\textrm{span}\,A 
		= 
		\cap \{V \subseteq L(\mathbb{K}^n): A \subseteq V\}.
	\end{equation}
	span\,$A$ is also called the
	``linear hull'' of $A$.
	Then with $\leqq$ defined as the usual set inclusion and
	with the following definitions
	for two subspaces $V$, $W$ $\subseteq$ $\mathbb{K}^n$,
	\begin{equation}
		\label{ortholattice-wedge-vee}
		V \wedge W = V \cap W,
		\qquad
		V \vee W = \textrm{span}\, (V \cup W),
	\end{equation}
	the set $L(\mathbb{K}^n)$ is a lattice.
	It has universal bounds $O=\{0\}$ and $I=\mathbb{K}^n$.
\end{beispiel}

\begin{lem}
	Let $L$ be a lattice $L$, and $x$, $y \in L$ with $x \leqq y$. Then
	for all $z\in L$,
	\begin{equation}
		\label{L-B}
		x \wedge z \leqq y \wedge z.
	\end{equation}
\end{lem}
\begin{proof}
	We have $x \wedge z \leqq z$ and
	$x \wedge z \leqq x \leqq y$ by (\ref{def-wedge-vee}),
	hence (\ref{L-B}) by (\ref{def-wedge-vee}) again.
\end{proof}

The binary operations $\wedge$ and $\vee$ in lattices have important algebraic
properties, some of them analogous to those of ordinary multiplication
and addition.

\begin{satz}
	In a lattice $L$, the operations of meet and join satisfy the following laws,
	whenever the expressions referred to exist.	
	\begin{eqnarray}
		\label{L1}
		\mbox{(Idempotent laws)}\index{idempotent laws} & &
		x\wedge x = x, \quad x \vee x = x.
		\\
		\label{L2}
		\mbox{(Commutative laws)}\index{commutative laws} & &
		x \wedge y = y \wedge x, \quad x \vee y = y \vee x.
		\\
		\label{L3} 
		\mbox{(Associative laws)}\index{associative laws} & &
		x \wedge ( y \wedge z ) = ( x \wedge y ) \wedge z, \quad
		x \vee ( y \vee z ) = ( x \vee y ) \vee z.
		\\
		\label{L4}
		\mbox{(Laws of absorption)}\index{absorption} & &
		x \wedge ( x \vee y ) = x \vee ( x \wedge y ) = x.
	\end{eqnarray}
	(The laws of absorption are often also called ``laws of contraction.'')
	Moreover, 
	\begin{equation}
	\label{consistency}
	\mbox{(Consistency)} \qquad
		x \leqq y
		\iff
		x \wedge y = x
		\iff x \vee y = y.
		\qquad \ \
	\end{equation}
\end{satz}
\begin{proof}
	The idempotence and the commutativity laws are evident from (\ref{def-wedge-vee}).
	The associativity laws (\ref{L3}) follow since $x \wedge (y \wedge z)$ and
	($x \wedge y) \wedge z$ are both equal to $\sup\{x,y,z\}$ whenever all
	expressions referred to exist.
	The equivalence between $x \leqq y$, $x \wedge y = x$, and $x \vee y = y$
	is easily verified. Thus $x \geqq y$ is equivalent to
	$x \wedge y = y$ and $x \vee y = x$, and this implies (\ref{L4}).
\end{proof}

It can be proved that the identities (\ref{L1}) -- (\ref{L4}) completely
charactarize lattices \cite[Theorem I.8]{Birkhoff-1973}.
In fact Dedekind,\index{Dedekind}
who first considered the concept of
a lattice (``Dualgruppe'') at the end of the 19th
century, used (\ref{L1}) -- (\ref{L4}) to define lattices.

\begin{satz}[Principle of Duality]
	\label{duality-principle}\index{duality principle}
	Given any valid formula over a lattice, the dual formula obtained by interchanging
	$\leqq$ with $\geqq$, and simultaneously $\wedge$ with $\vee$, is also valid.
\end{satz}
\begin{proof}
	Since for any elements $x$, $y$ of the lattice we have
	$x \leqq y$ if and only if $y \geqq x$, the poset structure
	with respect to $\geqq$ is isomorphic 
	to the poset structure
	$\leqq$, but with $\wedge$ and $\vee$ interchanged.
\end{proof}

The dual of a lattice is simply its Hasse diagram (Fig.~\ref{fig-lattices})
turned upside down, illustrating the principle of duality.
In fact, the two poset structures $\leqq$ and $\geqq$
of a lattice are tied up to each other by the laws of associativity, absorption, 
and consistency so strongly that they are inescapably dual.
The following theorem concerns relations of 
modality and distributivity which are
valid in every lattice.

\begin{satz}
	Let $L$ be a lattice.
	For all $x$, $y$, $z$ $\in$ $L$ we then have
	the ``modular inequality''
	\begin{equation}
		\label{modular-inequality}
		x \vee (y \wedge z) \leqq (x \vee y) \wedge z
		\qquad
		\mbox{if } x \leqq z,
	\end{equation}
	and the ``distributive inequalities''
	\begin{eqnarray}
		\label{5}
		x \wedge (y \vee z)
		&\hspace*{-.5em} \geqq \hspace*{-.5em}&
		(x \wedge y) \vee (x \wedge z)
		\\
		\label{5'}
		x \vee (y \wedge z)
		&\hspace*{-.5em} \leqq \hspace*{-.5em}&
		(x \vee y) \wedge (x \vee z) 
	\end{eqnarray}
\end{satz}
\begin{proof}
	If $x \leqq z$, we have with $x \leqq x \vee y$ that
	$x \leqq (x \vee y) \wedge z$.
	Also $y \wedge z \leqq y \leqq x \vee y$ and $y \vee z \leqq z$.
	Therefore, $y \wedge z \leqq (x \vee y) \wedge z$, i.e. 
	$x \vee (y \wedge z) \leqq (x \vee y) \wedge z$, which is
	(\ref{modular-inequality}).
	 
	Clearly $x \wedge y \leqq x$, and $x \wedge y \leqq y \leqq y \vee z$;
	hence $x \wedge y \leqq x \wedge (y \vee z).$ 
	Also $x \wedge z \leqq x$, $x \wedge z \leqq z \leqq y \vee z$;
	hence $x \wedge z \leqq x \wedge (y \vee z)$. 
	That is, $x \wedge (y \vee z)$ is an upper bound of $x \wedge y$ and
	$x \wedge z$, from which (\ref{5}) follows.
	The inequality (\ref{5'}) follows from (\ref{5}) by the principle of
	duality.
\end{proof}

\subsection{Distributive lattices} 
\begin{figure}[htp]
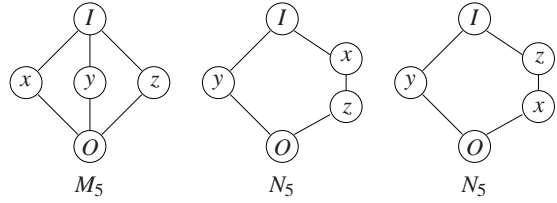

\centering
\begin{footnotesize}
\unitlength1ex
\begin{tabular}{*{3}{c}}
	\begin{picture}(15,15)
		\put(0,0){\includegraphics[height=15ex]{./M5}}
		\put( 7.5,14.3){\makebox(0,0)[t]{$I$}}
		\put( 1.0, 7.6){\makebox(0,0)[l]{$x$}}
		\put( 7.5, 7.5){\makebox(0,0){$y$}}
		\put(14.0, 7.6){\makebox(0,0)[r]{$z$}}
		\put( 7.5, 0.7){\makebox(0,0)[b]{$O$}}
	\end{picture}
	&
	\begin{picture}(15,15)
		\put(0,0){\includegraphics[height=15ex]{./N5}}
		\put(7.5,14.3){\makebox(0,0)[t]{$I$}}
		\put( 1.0, 7.5){\makebox(0,0)[l]{$y$}}
		\put(14.0,10.0){\makebox(0,0)[r]{$x$}}
		\put(14.0, 5.0){\makebox(0,0)[r]{$z$}}
		\put(7.5, 0.7){\makebox(0,0)[b]{$O$}}
	\end{picture}
	&
	\begin{picture}(15,15)
		\put(0,0){\includegraphics[height=15ex]{./N5}}
		\put(7.5,14.3){\makebox(0,0)[t]{$I$}}
		\put( 1.0, 7.5){\makebox(0,0)[l]{$y$}}
		\put(14.0, 9.8){\makebox(0,0)[r]{$z$}}
		\put(14.0, 5.2){\makebox(0,0)[r]{$x$}}
		\put(7.5, 0.7){\makebox(0,0)[b]{$O$}}
	\end{picture}
	\\
	$M_5$ & $N_5$ & $N_5$
\end{tabular}
\end{footnotesize}
\caption{\label{fig-M5-N5}\footnotesize
	Non-distributive lattices.
	}
\end{figure}

In many lattices, and thus in many logics, the analogy between the lattice
operations $\wedge$, $\vee$ and the arithmetic operations $\cdot$, $+$ includes
the distributive law $x(y+z)$ $=$ $xy+xz$. In such lattices, the distributive
inequalities (\ref{5}) and (\ref{5'}) can be sharpened to identities.
These identities do not hold in all lattices; for instance, they fail in the
lattices $M_5$ and $N_5$ in Figure \ref{fig-M5-N5}.
We now study distributivity, which in a lattice is symmetric with respect to
$\wedge$ and $\vee$ due to the duality principle (Theorem \ref{duality-principle}), 
which is not the case
in ordinary algebra, where $a+(bc)$ $\ne$ $(a+b)(a+c)$ due to the 
priority precedence
of multiplication ($\cdot$) and addition ($+$).

\begin{definition}
	\label{def-distributive}
	A lattice $L$ is called \emph{distributive}\index{distributive}
	if the following identity holds.
	\begin{equation}
		\label{L6'}
		x \wedge ( y \vee z) = ( x \wedge y) \vee (x \wedge z)
		\quad
		\mbox{for all } x,y,z \in L.
	\end{equation}
\end{definition}

\begin{satz}
	\label{satz-equivalence-L6'-L6''}
	In any lattice $L$, the identity (\ref{L6'}) is equivalent
	to
	\begin{equation}
		\label{L6''}
		x \vee ( y \wedge z) = ( x \vee y) \wedge (x \vee z)
		\quad
		\mbox{for all } x,y,z \in L.
	\end{equation}
\end{satz}
\begin{proof}
	We prove (\ref{L6'}) $\Rightarrow$ (\ref{L6''}). The converse
	(\ref{L6''}) $\Rightarrow$ (\ref{L6'})	follows analogously.
	\[
	\begin{array}{r@{\ =\ }ll}
		(x \vee y) \wedge (x \vee z)
		& [(x \vee y) \wedge x] \vee [(x \vee y) \wedge z]
		& \mbox{by (\ref{L6'})}
		\\
		& x \vee [z \wedge ( x \vee y)]
		& \mbox{by (\ref{L4}), (\ref{L2})}
		\\
		& x \vee [ (z\wedge x) \vee (z \wedge y)]
		& \mbox{by (\ref{L6'})}
		\\
		& [x \vee [z \wedge x)] \vee (z \wedge y)
		& \mbox{by (\ref{L3})}
		\\
		& x \vee (z \wedge y)
		& \mbox{by (\ref{L4})}
	\end{array}
	\]
\end{proof}

However, in nondistributive lattices the truth of (\ref{L6'}) for \emph{some} 
elements $x$, $y$, $z$
does not imply them obeying (\ref{L6''}), as the
two variants of $N_5$ in Figure \ref{fig-M5-N5} show.
An important property of distributive lattices is the following.

\begin{satz}
	\label{theo-distributive-uniqueness}
	A lattice $L$ is distributive if and only if the following property
	is satisfied for all $a$, $x$, $y\in L$:
	\begin{equation}
		\label{condition-distributivity}
		a \wedge x = a \wedge y \quad \mbox{and} \quad a \vee x = a \vee y
		\quad \mbox{imply} \quad x=y.
	\end{equation}
\end{satz}
\begin{proof}
	Suppose the lattice to be distributive.
	Then using repeatedly the Equations (\ref{L4}), (\ref{L2}), and (\ref{L6'}), we have
	\begin{eqnarray*} 
		x 
		& =
		& x \wedge (a \vee x) = x \wedge (a \vee y) = (x \wedge a) \vee (x \wedge y)
		\\ 
		& =
		& 
		(a \wedge y) \vee (x \wedge y) = (a \vee x) \wedge y = (a \vee y) \wedge y = y.
	\end{eqnarray*}
	The converse is proved in \cite[§II.7]{Birkhoff-1973}.
\end{proof}
Expressions involving the symbols $\wedge$, $\vee$ and elements of a lattice are
called \emph{lattice polynomials}\index{polynomial! lattice}\index{lattice! polynomial}

\begin{figure}[htp]
\centering
	\unitlength1ex
	\begin{picture}(11.5,12)
		\put(0,0){\includegraphics[width=12ex]{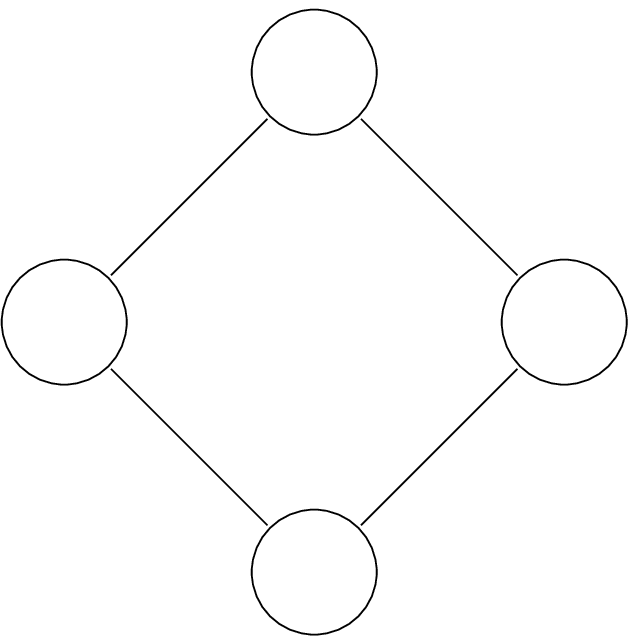}}
		\put(6,11.2){\makebox(0,0)[t]{$u$}}
		\put( 1.2,6.1){\makebox(0,0){$x$}}
		\put(10.9,6.1){\makebox(0,0){$y$}}
		\put(6,.8){\makebox(0,0)[b]{$v$}}
	\end{picture}
	\caption{\label{fig-F2} \footnotesize
	The lattice $F_2$ $\cong$ ${\mathbf{2}}^2$.
	For ${\mathbf{2}}^2$, we have $v=00$, $x=01$, $y=10$, and
	$u=11$.
	}
\end{figure}

\begin{lem}
	In any lattice $L$, the sublattice $S$ generated by two elements $x$ and $y$
	consists of $x$, $y$, $u$, and $v$, where $u = x \vee y$ and $v=x \wedge y$,
	as in Figure \ref{fig-F2}.
\end{lem}
\begin{proof}
	By (\ref{L4}), $x\wedge u$ $=$ $x$; 
	by (\ref{L3}), (\ref{L1}), 
	$x \vee u$ $=$ $x \vee (x \vee y)$ $=$ $(x \vee x) \vee y$ $=$ $x \vee y$ $=$ $u$.
	The other cases are analogous, using symmetry in $x$ and $y$ and duality. 
\end{proof}

A \emph{lattice morphism}\index{morphism} is a mapping 
$\mu$: $L$ $\to$ $K$ from a lattice $L$ to a lattice $K$ which preserves
the meet and join operations, i.e.,
$\mu(x \wedge y)$ $=$ $\mu(x) \wedge \mu(y)$,
$\mu(x \vee y)$ $=$ $\mu(x) \vee \mu(y)$
for all $x$, $y$ $\in$ $L$.

\begin{cor}
	\label{cor-sublattice-F2}
	Let $F_2 \cong {\mathbf{2}}^2$ be the lattice of Figure \ref{fig-F2}, and let
	$a$, $b$ $\in$ $L$ be two elements of an arbitrary lattice.
	Then the mapping $x \mapsto a$, $y \mapsto b$ can be extended to a 
	lattice morphism $\mu: F_2 \to L$.
\end{cor}

The preceding results are usually summarized in the statement that $F_2$ is the
``free lattice'' with generators $x$, $y$. It has just four elements and
is distributive; in fact it is a Boolean lattice.

Lattice polynomials in three or more variables can be extremely complicated.
However, in a distributive lattice any polynomial can be brought to a
normal form, similarly as a real or complex polynomial can be written
as a sum of products, 
$p(x_1$, \ldots, $x_n)$ $=$ $\sum_{i=1}^r \left(\prod_{j=1}^{s(i)} x_{ij} \right)$
or as a product of its divisors,
$p(x_1$, \ldots, $x_n)$ $=$ $\prod_\delta \left(\sum_{j\in T_\delta} x_j \right)$.

\begin{satz}
	\label{theo-normal-forms}
	In a distributive lattice $L$, every polynomial $p:L^n \to L$
	of $n$ variables is equivalent to a join of meets, and dually:
	\begin{equation}
		p(x_1, \ldots, x_n)
		= \bigvee_{\alpha \in A} \Big( \bigwedge_{i\in S_\alpha} x_i \Big)
		= \bigwedge_{\delta \in D} \Big( \bigvee_{j\in T_\delta} x_j \Big),
	\end{equation}
	where $S_\alpha$ and $T_\delta$ are nonempty sets of indices.
\end{satz}
\begin{proof}
	Each single $x_i$ can be so written, where $A$ (or $D$, respectively) is the
	family of sets consisting of the single element set $\{x_i\}$.
	On the other hand, we have by (\ref{L1})--(\ref{L3})
	\begin{equation}
		\bigvee_{\alpha \in A} \Big( \bigwedge_{i\in S_\alpha} x_i \Big)
		\vee \bigvee_{\beta \in B} \Big( \bigwedge_{i\in S_\beta} x_i \Big)
		= \bigvee_{\gamma \in A\cup B} \Big( \bigwedge_{i\in S_\gamma} x_i \Big).
	\end{equation}
	Using the distributive law, we have similarly
	\begin{equation}
	\label{polynomial-distributive-law}
		\bigvee_{\alpha \in A} \Big( \bigwedge_{i\in S_\alpha} x_i \Big)
		\vee \bigvee_{\beta \in B} \Big( \bigwedge_{i\in S_\beta} x_i \Big)
		= \bigvee_{\gamma \in A \times B} 
		\Big( \bigwedge_{i\in S_\alpha \cup S_\beta} x_i \Big).
	\end{equation}
	The assertion follows from (\ref{L2}) and (\ref{L6'}), combined with the relation
	$(\bigwedge_S x_i) \wedge (\bigwedge_T x_i)$ $=$ $\bigwedge_{S\cup T} x_i$,
	which follows from (\ref{L1})--(\ref{L3}).
\end{proof}

Equation (\ref{polynomial-distributive-law}) is the lattice generalization
of the distributive law of ordinary algebra,
\[
	\Big(\sum_{\alpha\in A} x_\alpha \Big) \Big( \sum_{\beta\in B} y_\beta\Big)
	=
	\sum_{(i,j)\in A\times B} x_i y_j.
\]

\section{Logics}
A general logic is now going to be introduced as
a lattice with universal bounds and a
special operation, the fuzzy negation. Here the fuzzy negation of a
lattice element $x$ is the square root of
a unique supremum $x''$ of $x$, cf.\ Eq.~(\ref{logic-1}).
Considering the elements of the lattice as propositions,
each proposition then implies its double negation, but not always vice versa.
Furthermore, the negation is antitone, a property which turns out to be equivalent
to the disjunctive De Morgan law, and the Boolean boundary condition
for the universal bounds holds true. If nothing else is assumed,
then the lattice is a fuzzy logic. If in addition the law of non-contradiction
holds, then it is a logic. In this way, the notion ``fuzzy logic'' includes
the propositional structures of 
substructural logics focussing on validity and premise combination
\cite{Restall-2000},
infinite-valued fuzzy logics of {\L}ukasiewiczian
type \cite[\S2.3.2]{Buldt-1999},
as well as quantum and distributive logics, in particular 
Boolean
algebras or nonclassical Heyting%
-Brouwerian
(``intuitionistic'') logics
\cite[\S\,XII.3]{Birkhoff-1973} in which
``tertium non datur'' or ``reductio ad absurdum'' are not valid.
A fuzzy logic is not necessarily
distributive, not even orthomodular or paraconsistent.
Quantum logics will turn out to be orthomodular,
but not necessarily distributive.

\begin{quote}
\begin{em}
	Creation comes when you learn to say no.
\end{em}

\hfill{Madonna, \emph{The Power of Goodbye}}
\end{quote}

\begin{definition}
	\label{def-complement}
	Let $L$ be a lattice with universal bounds $0$ and $1$, i.e.,
	$0 \leqq x \leqq 1$ for all $x \in L$.
	A mapping 
	$': L \to L$, $x \mapsto x'$, is called \emph{fuzzy negation}\index{negation},
	if the following relations hold for all $x$, $y \in L$:
	\index{weak double negation}%
	\index{negation! weak double -}\index{double negation! weak -}
	\begin{eqnarray}
		\label{logic-1}
		\mbox{\emph{(Weak double negation)}}
		\label{double-negation}
		 & &
		 \, \, \,
		 x \leqq (x')',
		\\ 
		\mbox{\emph{(Antitony)}}
		& &
		\,
		\label{antitony}
		y' \leqq x' \quad \mbox{if } x \leqq y,
		\qquad \qquad
		\\ 
		\mbox{\emph{(Boolean boundary condition)}}%
		\index{Boolean boundary condition}%
		& &
		\label{Boolean-limit}
		0' = 1, \quad 1' = 0.
	\end{eqnarray}
	The pair $(L,{}')$ then is called a 
	\emph{fuzzy logic}\index{fuzzy logic}\index{logic}\index{logic! fuzzy},
	and the elements $x\in L$ are called \emph{propositions}\index{proposition}.
	If for the fuzzy negation the ``law of non-contradiction''\index{law of non-contradiction}%
	\index{non-contradiction}
	\begin{equation}
		x \wedge x' = 0
	\end{equation}
	holds for all $x\in L$, then it is called 
	\emph{(non-contradictory) negation}%
	\index{non-contradictory -}\index{non-contradictory negation}
	and $(L,{}')$ is a \emph{logic}.
	As long as misunderstanding is excluded, we shortly write $L$ instead of $(L,{}')$.

	Algebraically, the element 
	$x'\in L$ in a (non-contradictory) logic $L$ 
	is called a \emph{pseu\-do-com\-ple\-ment}\index{pseudo-complement} of $x\in L$;
	if in addition $x \vee x' = 1$, then $x'$ is called
	a \emph{complement}\index{complement} of $x\in L$.
	In general, a lattice is called
	\emph{(pseudo-) complemented} if all its elements have (pseudo-) complements.
	A mapping $'$: $L$ $\to$ $L$, $x$ $\mapsto$ $x'$ 
	in a (pseudo-) complemented lattice $L$, assigning
	to each element $x$ a (pseudo-) complement, is called 
	\emph{(pseudo-) complementation}.%
	\index{pseudo-complementation}\index{complementation}
	If the (pseudo-) complementation is bijective, the lattice is called
	\emph{uniquely (pseudo-) complemented}\index{uniquely complemented}.
%
\end{definition}

Defined this way, a logic is a special fuzzy logic.
In any fuzzy logic $L$ 
the relation $x \leqq y$ will be interpreted as the statement
``$x$ implies $y$.''
The propositions $x \wedge y$ and $x \vee y$ will be interpreted as
``$x$ and $y$'' and ``$x$ or $y$,'' respectively. 
The universal bounds of a logic are usually denoted by 0 and 1, the
proposition 1 expresses truth, and the proposition 0 expresses falsehood
or absurdity.

\begin{satz}
	\label{theo-negation-antitony}
	Let be $L$ a lattice and $': L \to L$ a mapping satisfying (\ref{double-negation})
	for all $x \in L$. Then $'$ is antitone if and only if the 
	disjunctive De Morgan law\index{disjunctive De Morgan law}%
	\index{De Morgan's law! disjunctive -}
	holds, i.e.,
	\begin{equation}
		\label{disjunctive-De-Morgan}
		(x \vee y)' = x' \wedge y'
		\qquad \mbox{%
			for all $x$, $y \in L$.%
		}
	\end{equation}
\end{satz}
\begin{proof}
	Suppose first the antitony of $'$, and let
	$u := x \vee y$ and $v := x' \wedge y'$ for arbitrary $x$, $y \in L$.
	Then $u \geqq x$ and $u \geqq y$, as well as $v \leqq x'$ and $v \leqq y'$.
	By the antitony, this means that 
	$u' \leqq x'$ and $u' \leqq y'$, as well as
	$v' \geqq x'$ and $v' \geqq y'$, i.e.,
	\begin{equation}
		u' \leqq x' \wedge y' = v,
		\label{eq-antitony-De-Morgan-1}
	\end{equation}
	as well as
	$v' \geqq x' \vee y' = u$.
	By (\ref{double-negation}) and the antitony, the last inequality yields
	$ 
		v \leqq v'' \leqq u' 
	$ 
	which means together with (\ref{eq-antitony-De-Morgan-1}) that $u' = v$.
	
	Assume, on the other hand, the disjunctive De Morgan law (\ref{disjunctive-De-Morgan}).
	Since $y = x \vee y$ for $x \leqq y$, we have 
	$y' = (x \vee y)' = x' \wedge y'$
	with (\ref{disjunctive-De-Morgan}),
	hence $y' \leqq x'$.
\end{proof}

\begin{satz}
	\label{theo-conjunctive-De-Morgan}
	In a general fuzzy logic $L$ the conjunctive De Morgan inequality
	\begin{equation}
		\label{conjunctive-De-Morgan-inequality}
		(x \wedge y)' \geqq x' \vee y'.
	\end{equation}
	holds for all $x$, $y \in L$.
\end{satz}
\begin{proof}
	Since
	$(x' \vee y')' = x'' \wedge y'' \geqq x \wedge y$
	by (\ref{disjunctive-De-Morgan}) and (\ref{double-negation}) 
	for all $x$, $y \in L$, we have
	$x' \vee y' \leqq ((x' \vee y')')' \leqq (x \wedge y)'$
	by (\ref{double-negation}) and the antitony of the negation.
\end{proof}

Note that the conjunctive De Morgan law (see Eq.~(\ref{conjunctive-De-Morgan} below) 
does not necessarily hold in a fuzzy logic, even not in a logic.
%

\begin{bem}
	If we abandon the Boolean boundary condiditon (\ref{Boolean-limit})
	on a negation, then very little is known about the fuzzy negations of 0 and 1
	in a general fuzzy logic.
	By the antitony (\ref{antitony})
	and by the general lattice property $0 \leqq x \leqq 1$ for all $x \in L$,
	we only can derive
	\begin{equation}
		1' \leqq x' \leqq 0' \qquad \mbox{for all $x\in L$.}
	\end{equation}
	A fuzzy negation with $0' = 0$ thus must be constant, i.e., $x'=0$ for all $x\in L$.
	On the other hand, a constant fuzzy negation $x' = x_0 \in L$ for all $x\in L$
	implies $x_0 = 0$, since otherwise we had $0'' = x_0 > 0$, 
	contradicting (\ref{double-negation}).
	In a logic with the law of non-contradiction
	$x \wedge x' = 0$, however, we have
	\begin{equation}
		\label{eq-pseudocomplement-0-1}
		1' = 0,
	\end{equation}
	since $1 \wedge y = y$ for all $y\in L$.
\end{bem}


\begin{bem}
Sometimes the notions ``strong'' and ``weak negation'' are used,
especially in the context of logic programming\index{logic programming}
\cite{Boley-2003}
and artificial intelligence\index{artificial intelligence},
motivated by the following ideas.
Intuitively speaking, 
strong negation captures the presence of explicit negative information,
while weak negation captures the absence of positive information.
In computer science, weak negation 
captures the computational concept of negation-as-failure 
(or ``closed-world negation'').

A strong negation\index{negation! strong -} ($'$) can be interpreted as ``impossible.'' 
Negating this, in turn, gives a weak ``not impossible''
assertion, so that $x$ implies $(x')'$, i.e, $x \leqq (x')'$,
but not vice versa. 
With this negation, the rule of bivalence $x \vee x'=1$ 
($x$ is true or impossible) does not necessarily hold  (since
$x$ is possible as long as it is not recognized as true), 
but the corresponding $x \wedge x' = 0$ (not both true and impossible) does.
Weak negation\index{Negation! weak -}, in contrast, can be regarded 
as ``unconfirmed.'' 
Negating this gives 
$x \leqq (x')'$
(if $x$ is true it is always unconfirmed that it is unconfirmed),
but not vice versa
$(x')' \leqq x$ 
(if it is not confirmed that $x$ is unconfirmed, then $x$ is certainly true).
However, it has $x \vee x'=1$ ($x$ is true or unconfirmed) holding,
since if $x$ is not true it is certainly not confirmed, 
but not necessarily
$x \wedge x' = 0$ (it is never confirmed that $x$ is true and unconfirmed), 
since $x$ may be true but not confirmed.
Defining such kind of weak negation therefore implies that \emph{tertium non datur}
$x \vee x'=1$ does hold, but the law of contradiction is not necessarily true,
$x \wedge x' \geqq 0$.
Likewise, the ability to speak of the \emph{uncertain} apparently may force a 
weakening of the \emph{tertium non datur} in some form, 
losing double negation, and also of the law of contradiction.
Both possibilities are enabled by the above concept of a fuzzy logic.

However, there is some confusion with the term ``strong negation.''
Sometimes it simply means the classical negation ``false'' $=$ ``not true.''
Negating this gives another strong assertion,
$(x')' \leqq x$ (if it is false that $x$ is false, then $x$ is true),
and vice versa. 
This has $x \vee x' = 1$ ($x$ is true or false)
and $x \wedge x' = 0$ ($x$ is not both true and false) holding.
\end{bem}

What are the reasons that the conjunctive version
(Eq.~(\ref{conjunctive-De-Morgan}) below) is not implied by the antitony 
of the negation?
In fact, it is easily proved that for lattices with a total order
(i.e., $\forall x$, $y$ either $x \leqq y$ or
$y \leqq x$) antitony, disjunctive De Morgan law and conjunctive De Morgan law
are equivalent. However, in a partially ordered set this is not necessarily true.
One of the simplest counterexamples is $M_5$.

\begin{beispiel}
	\label{bsp-M5-non-conjunctive}
	Let $L=M_5$\index{M$_5$} denote the modular lattice as in 
	Figure \ref{fig-M5} and
	\begin{figure}[htp]
	\centering
	\begin{footnotesize}
	\unitlength1ex
		\begin{picture}(15,15)
			\put(0,0){\includegraphics[height=15ex]{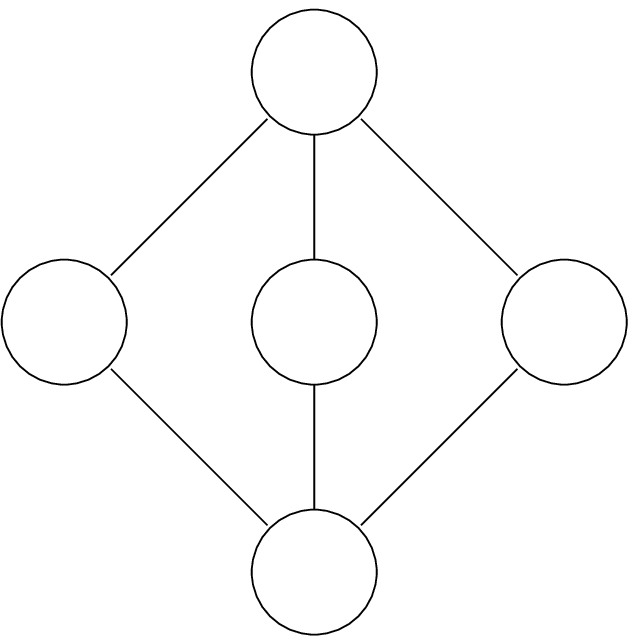}}
			\put( 7.5,14.3){\makebox(0,0)[t]{$1$}}
			\put( 1.0, 7.4){\makebox(0,0)[l]{$a$}}
			\put( 7.5, 7.5){\makebox(0,0){$b$}}
			\put(14.0, 7.4){\makebox(0,0)[r]{$c$}}
			\put( 7.5, 0.7){\makebox(0,0)[b]{$0$}}
		\end{picture}
	\end{footnotesize}
	\caption{\label{fig-M5}\footnotesize
		The non-distributive modular lattice $M_5$.
	}
	\end{figure}
	define the operation $': M_5 \to M_5$ by
	\begin{equation}
		\begin{array}{c|*{5}{c}}
			x      & 0 & a & b & c & 1
			\\ \hline
			x' & 1 & c & 0 & a & 0
		\end{array}
	\end{equation}
	Then $x'' = x$ for all $x \ne b$, but $b'' = 1 > b$,
	and antitony is easily verified.
	Since moreover $x \wedge x' = 0$ for all $x \in M_5$,
	${}'$ is a non-contradictory negation and $(M_5,{}')$ is a logic.
	However, 
	in contrast to the disjunctive De Morgan law (\ref{disjunctive-De-Morgan}),
	the conjunctive De Morgan law (\ref{conjunctive-De-Morgan})
	is not valid since, e.g.,
	$(a \wedge b)' = 1$ but $a' \vee b' = c$.
\end{beispiel}

\begin{satz}
	\textbf{\em (De Morgan's conjunctive law)}
	\label{theo-De-Morgan-fuzzy}
	If in a fuzzy logic $L$ we have
	$(x')'$ $=$ $x$ for all $x$ $\in$ $L$,
	then the conjunctive De Morgan\index{conjunctive De Morgan's law} law holds,
	\begin{equation}
		\label{conjunctive-De-Morgan}
		(x \wedge y)' = x' \vee y'
		.
	\end{equation}
\end{satz}
\begin{proof}
	Let $(x')' = x$.
	By (\ref{disjunctive-De-Morgan}) we obtain
	$(x' \vee y')' = x \wedge y$, hence
	$x' \vee y' = ((x' \vee y')')' = (x \wedge y)'$.
\end{proof}

\noindent
For a non-contradictory logic the above theorem has a stronger consequence.

\begin{satz}
	\textbf{\em (Tertium non datur)}
	\label{theo-tertium-non-datur}
	In a logic $L$, a negation with
	$(x')'$ $=$ $x$ for all $x$ $\in$ $L$ implies 
	the law {``tertium non datur''}\index{tertium non datur},
	or ``law of excluded middle,''\index{law of excluded middle}\index{excluded middle}
	$ 
		x \vee x' = 1.
	$ 
\end{satz}
\begin{proof}
	By Theorem \ref{theo-De-Morgan-fuzzy},
	De Morgan's laws (\ref{disjunctive-De-Morgan}) and (\ref{conjunctive-De-Morgan}) 
	hold, and hence
	the property $x \wedge x' = 0$ for all $x \in L$
	implies $(x' \vee x)' = x'' \wedge x' = 0$, i.e.,
	$x' \vee x = 1$ by (\ref{Boolean-limit}).
%
\end{proof}

\noindent
In general, a non-contradictory negation satisfying
the law \emph{tertium non datur} is called
\emph{ortho-negation}\index{ortho-negation} \cite{Restall-2000},
\emph{complemented negation}, or \emph{involutive negation},
\index{complemented logic}\index{involutive logic}%
and $(L,{}')$ a \emph{complemented logic}.

\begin{beispiel}
	\label{bsp-tertium}
	\cite[§2.2]{Svozil-1998}
	The \emph{tertium non datur} is a highly nontrivial assumption.
	An example for its nonconstructive feature is a proof of the following
	proposition: 
	\emph{%
		``There exist irrational numbers $x$, $y$ $\in$ $\mathbb{R}\setminus
		\mathbb{Q}$ with $x^y \in \mathbb{Q}$.''%
	}
	\emph{Proof:}
	Either $\sqrt{2}{}^{\sqrt{2}} \in \mathbb{Q}$, i.e., $x$ $=$ $y$ $=$ $\sqrt{2}$;
	or $\sqrt{2}{}^{\sqrt{2}} \notin \mathbb{Q}$, then 
	$(\sqrt{2}{}^{\sqrt{2}})^{\sqrt{2}} = 2 \in \mathbb{Q}$,
	i.e., $x$ $=$ $\sqrt{2}{}^{\sqrt{2}}$, $y$ $=$ $\sqrt{2}$.
	Q.E.D.
	The question whether or not $\sqrt{2}{}^{\sqrt{2}}$ is rational,
	however, remains unsolved in the proof.
\end{beispiel}

%
\begin{figure*}[htp]
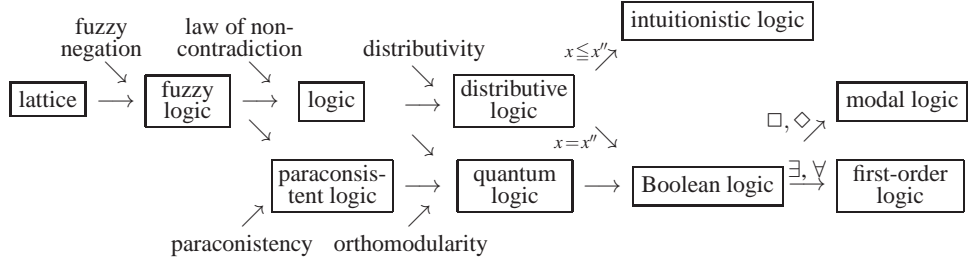

\centering
	\begin{footnotesize}%
	\begin{tabular}{*{11}{@{$\!$}l@{$\!$}}}
		\\[-1ex]
		 & &
		 \multicolumn{2}{r}{}
		 & \multicolumn{1}{c}{}
		 & 
		 \multicolumn{2}{l}{}
		 & 
		 & \hspace*{-3.5ex}\fbox{intuitionistic logic}
		 \\[-2.5ex]
		 \multicolumn{2}{r}{%
		 	\begin{tabular}{c} fuzzy\\[-.7ex] negation\end{tabular}%
		 \hspace*{-1.2em}}
		 &
		 \multicolumn{2}{r}{%
		 	\begin{tabular}{c} law of non-\\[-.7ex] contradiction\end{tabular}%
		 \hspace*{-1.5em}}
		 &
		\multicolumn{2}{r}{%
		 	\hspace*{-2.5em}%
		 	\begin{tabular}{c} \phantom{Satz}\\[-.7ex] distributivity\end{tabular}%
		 	\hspace*{-2.5em}%
		}
		 & &
		\\[-2.ex]
		\fbox{lattice} & 
		\begin{tabular}{c} 
			$\searrow$ \\ $\longrightarrow$ \\ \phantom{$\searrow$} 
		\end{tabular}
		&
		\begin{tabular}{|c|} \hline fuzzy\\[-.7ex] logic \\ \hline\end{tabular} 
		&
		\begin{tabular}{c} 
			$\searrow$ \\ $\longrightarrow$ \\ {$\searrow$} 
		\end{tabular}
		&
		\fbox{logic} & 
		\begin{tabular}{c}
			\\[-1ex]
			$\searrow$ \\ 
			$\longrightarrow$ \\[1ex] 
			$\searrow$ 
		\end{tabular}
		&
		\multicolumn{2}{l}{%
			\hspace*{.5ex}%
			\begin{tabular}{r} 
				${}^{x \, \leqq \, x''}\hspace*{-2.0ex}\nearrow$\\
				\hspace*{-1.7em}%
				\begin{tabular}{|@{\ }c@{\ }|}\hline 
					distributive \\[-.7ex] logic 
				\\ \hline \end{tabular}
				\phantom{$\longrightarrow$}\\[-.4ex] 
				${}_{x \, =\, x''}\hspace*{-1ex}\searrow$
				\\[.7ex]
			\end{tabular}
		}
		& & &
		 \fbox{modal logic}
		\\
		 & & & &
		\hspace*{-2.5ex}%
		\begin{tabular}{|@{\ }c@{\ }|} \hline 
		 	paraconsis- \\[-.7ex] tent logic
		\\ \hline\end{tabular}%
		 & \multicolumn{1}{l}{$\longrightarrow$} & 
		 \begin{tabular}{|c|}\hline quantum \\[-.7ex] logic \\ \hline \end{tabular}%
		 & 
		 \multicolumn{1}{r}{$\longrightarrow$}\hspace*{1.5ex}
		 & \hspace*{-2.5ex}\fbox{Boolean logic}
		&
		\hspace*{-4ex}
		\begin{tabular}{c}
			\\[-8ex]
			\hspace*{-2ex}${}^{\textstyle \Box, \Diamond}\hspace*{-1.5ex}\nearrow$
			\\[1ex]
			$\exists$, $\forall$\\[-1.5ex]
			$\longrightarrow$ 
		\end{tabular}
		& \begin{tabular}{|c|}\hline first-order \\[-.7ex] logic \\ \hline \end{tabular}
		\\[-1.ex]
		& & & 
		\multicolumn{1}{l}{$\nearrow$\hspace*{1ex}}
		& & 
		\multicolumn{1}{r}{$\nearrow$\hspace*{1ex}}
		& & & 
		\\[-0ex]
		& & & \hspace*{-5ex}paraconistency\hspace*{-8ex}%
		& & 
		\hspace*{-5ex}orthomodularity\hspace*{-8ex}%
		& & & & &
	\end{tabular}%
	\end{footnotesize}%
	\caption{\label{fig-logics} \footnotesize
		The algebraic hierarchy of logics.
		In particular, a Boolean logic is a special quantum logic, a quantum logic
		is a special fuzzy logic.
		By Theorem \ref{theo-quantum-logic}, a logic with the law \emph{tertium non datur}
		is a quantum logic.
		Establishing Boolean logic with the quantifiers $\exists$ and $\forall$
		yields first-order logic, and with the quantifiers $\Box$ and $\Diamond$
		modal logic.
	}
\end{figure*}
%
\begin{definition}
	\label{def-logic}
	A \emph{paraconsistent logic}\index{paraconsistent logic}\index{logic! paraconsistent -}
	is a fuzzy logic satisfying the 
	\emph{paraconsistency condition}\index{paraconsistency}
	\begin{equation}
		\label{paraconsistency}
		x = y
		\mbox{\qquad if $x \leqq y$ and $x' \wedge y = 0$.} 
	\end{equation}
	An \emph{intuitionistic logic}\index{intuitionistic logic}%
	\index{logic! intuitionistic -}
	is a distributive logic in which there exists propositions $x<(x')'$.
	A \emph{quantum logic}\index{quantum! logic}\index{logic! quantum -},
	or \emph{orthologic}\index{orthologic},
	is a logic satisfying the 
	\emph{orthomodular identity}\index{orthomodular}
	\begin{equation}
		\label{orthomodular-identity}
		x \vee (x' \wedge y) = y
		\mbox{\qquad if $x \leqq y$,} 
	\end{equation}
	A \emph{Boolean logic}\index{logic! Boolean -}
	is a complemented distributive logic.
\end{definition}

Therefore we obtain the algebraic structure of logics in Figure
\ref{fig-logics}.
Every distributive complemented lattice is othomodular,
since interchanging $x$ and $y$ and setting $z = x'$
for $x \leqq y$ in the distributive law (\ref{L6'}) with
$x \wedge y$ $=$ $x$,
$x \wedge x'$ $=$ 0
und $y = x \vee y$ yields
(\ref{orthomodular-identity}).

\begin{beispiel}
	\label{bsp-BN4-MO1}
	\textbf{(The logics BN$_4$ and MO$_1$)}
	A given a lattice may yield the propositional structure for
	more than one fuzzy logics, depending on the negation.
	\begin{figure}[htp]
	\centering
		\unitlength1ex
		\begin{picture}(11.5,12)
			\put(0,0){\includegraphics[width=12ex]{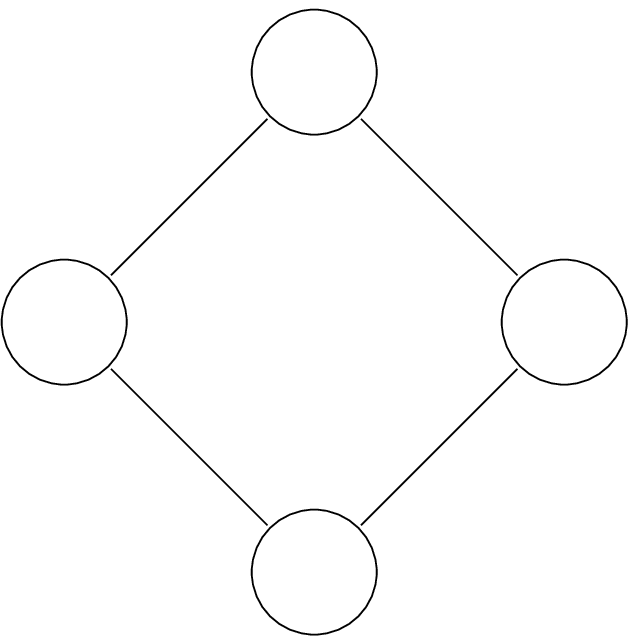}}
			\put(6,11.4){\makebox(0,0)[t]{$t$}}
			\put( 1.2,6.1){\makebox(0,0){$b$}}
			\put(10.9,6.1){\makebox(0,0){$n$}}
			\put(6,.2){\makebox(0,0)[b]{$f$}}
		\end{picture}
		\caption{\label{fig-F2-BN-4} \footnotesize
			The lattice $F_2$ $\cong$ ${\mathbf{2}}^2$,
			providing the propositional structure of the logics 
			BN$_4$ and MO$_1$.
		}
	\end{figure}
	A simple example is the Boolean lattice $\mathbf{2}^2$.
	Defining two negations $x^\sim$ and $x^\bot$
	by the following tables,
	\begin{equation}
		\begin{array}{c|cccc}
			x      & t & b & n & f 
			\\ \hline
			x^\sim & f & b & n & t
		\end{array}
		\qquad
		\begin{array}{c|cccc}
			x      & t & b & n & f 
			\\ \hline
			x^\bot & f & n & b & t
		\end{array}
	\end{equation}
	we obtain the logics\index{BN$_4$}\index{MO$_1$}
	\cite[§8.2]{Restall-2000}, \cite[§2.4]{Svozil-1998}
	\begin{equation}
		\mathrm{BN}_4 = (\{f,n,b,t\}, \ {}^\sim),
		\qquad
		\mathrm{MO}_1 = (\{f,n,b,t\}, \ {}^\bot),		
	\end{equation}
	Both negations may be illustrated geometrically, supposing
	$f$, $n$, $b$, $t$ as the four points on the unit circle $S^1$ in the 
	plane $\mathbb{R}^2$, viz.,
	$f=(0,-1)$, $n=(1,0)$, $b=(-1,0)$, $t=(0,1)$.
	The negation ${}^\sim$ then corresponds to the reflection in the 
	horizontal line \{($x$, 0)\} through the origin,
	whereas the negation ${}^\bot$ corresponds to a
	rotation around the origin by the angle $\pi$.
	The main differences between the logics BN$_4$ and MO$_1$ are that
	BN$_4$ is contradictory (e.g., $b \wedge b^\sim = b > f$)
	and that the \emph{tertium non datur} does not hold in
	BN$_4$ ($b \vee b^\sim = b < t$), whereas MO$_1$ is a classical Boolean
	logic.
	The fusion and material implication of BN$_4$ read:
	\begin{equation}
		\begin{array}{c|cccc}
			* & f & n & b & t
			\\ \hline
			 f    & f & f & f & f \\
			 n    & f & f & n & n \\
			 b    & f & n & b & t \\
			 t    & f & n & t & t
		\end{array}
		\qquad
		\begin{array}{c|cccc}
			\to   & f & n & b & t
			\\ \hline
			 f    & t & t & t & t \\
			 n    & n & t & n & t \\
			 b    & f & n & b & t \\
			 t    & f & n & f & t
		\end{array}
	\end{equation}
	The many-valued logic BN$_4$ considered by Dunn and Belnap was
	the result of research on relevance logic, but it also has significance for 
	computer science applications. The truth degrees may be interpreted as indicating,
	e.g., with respect to a database query for some particular state of affairs, 
	that there is
		no information concerning this state of affairs ($n$ $=$ $\emptyset$),
		information saying that the state of affairs fails ($f$ $=$ \{0\}),
		information saying that the state of affairs obtains ($t$ $=$ \{1\}),
		conflicting information saying that the state of affairs 
		obtains as well as fails ($b$ $=$ \{0, 1\}).
\end{beispiel}

\begin{beispiel}
\label{bsp-vector-space-2}
	The lattice $L(V_n)$ of all linear subspaces of an
	$n$-dimensional vector space $V_n$ (Example \ref{bsp-vector-space})
	is uniquely complemented since the orthogonal complement
	$V^\bot$ of any subspace $V$ satisfies
	$V \wedge V^\bot = \{0\}$ and $V \vee V^\bot = V_n$.
	Note that $O=\{0\}$ and $I=V_n$ are the universal bounds
	of $L(V_n)$.
	Moreover, it is orthomodular, and the complementation is
	involutive, $(x^\bot)^\bot = x$.
	Therefore,  $L(V_n)$ is an involutively complemented 
	logic in which De Morgan's laws hold.
\end{beispiel}

\begin{satz}
	\label{theo-quantum-logic}
	A logic is a quantum logic if and only if $x=x''$.
\end{satz}
\begin{proof}
	In a quantum logic we have by $y=x''$ in 
	(\ref{orthomodular-identity})
	that $x'' = x \vee (x' \wedge x'') = x \vee 0 = x$.
	Conversely, if $x=x''$, then from (\ref{disjunctive-De-Morgan}),
	(\ref{conjunctive-De-Morgan}) and (\ref{5}) we deduce
	\begin{equation}
		(x \vee (x' \wedge y))'
		= x' \wedge (x' \wedge y)'
		= x' \wedge (x \vee y')
		\geqq (x' \wedge x) \vee (x' \wedge y')
		= (x \vee y)',
	\end{equation}
	hence $x \vee (x' \wedge y) \geqq x \vee y$. But 
	$x \vee (x' \wedge y) \leqq (x \vee x') \wedge (x \vee y)$ by (\ref{5'}),
	and since with Theorem \ref{theo-tertium-non-datur} $x \vee x' = 1$, we have
	$x \vee (x' \wedge y) \leqq x \vee y$. We conclude
	$x \vee (x' \wedge y) = x \vee y = y$ for $x \leqq y$,
	i.e., (\ref{orthomodular-identity}) holds.
\end{proof}

With Theorem \ref{theo-tertium-non-datur}, in a quantum logic
therefore the law \emph{tertium non datur} holds, especially in Boolean logic.

\begin{beispiel}
	A complemented lattice, in which the orthomodular identity does not hold, 
	is $O_6$ in Figure \ref{fig-O6}.
\begin{figure}[htp]
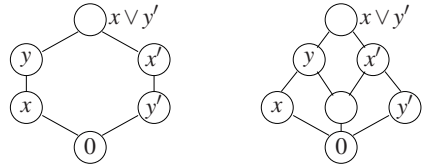

\centering
\begin{footnotesize}
	\unitlength1ex
	\begin{picture}(14,15)
		\put(0,0){\includegraphics[height=15ex]{./O6}}
		\put( 9.3,13.7){\makebox(0,0)[l]{$x \vee y'$}}
		\put( 1.5, 9.8){\makebox(0,0){$y$}}
		\put( 1.5, 5.3){\makebox(0,0){$x$}}
		\put(13.5,10.0){\makebox(0,0){$x'$}}
		\put(13.5, 5.3){\makebox(0,0){$y'$}}
		\put( 7.5, 1.6){\makebox(0,0){$0$}}
	\end{picture}
	\qquad \qquad
	\begin{picture}(14,15)
		\put(0,0){\includegraphics[height=15ex]{./L7}}
		\put( 9.3,13.7){\makebox(0,0)[l]{$x \vee y'$}}
		\put( 4.5, 9.8){\makebox(0,0){$y$}}
		\put( 1.5, 5.3){\makebox(0,0){$x$}}
		\put(10.5,10.0){\makebox(0,0){$x'$}}
		\put(13.5, 5.3){\makebox(0,0){$y'$}}
		\put( 7.5, 1.6){\makebox(0,0){$0$}}
	\end{picture}
	\end{footnotesize}
	\caption{\label{fig-O6}\footnotesize
		The non-orthomodular complemented lattice $O_6$\index{O$_6$}.
		If $x' \wedge y > 0$, then it forms the lattice $L_7$ which is
		orthomodular.
	}
\end{figure}
	Although $x<y$, we have
	$x \vee (x' \wedge y) = x \vee 0 = x \ne y$.
	If a lattice contains $O_6$ as a sublattice, it is non-orthomodular.
\end{beispiel}

\section{Fuzzy logics}
Propositions of every-day language, for instance ``Tom is big'' or 
``The room is cold,'' 
often are very vaguely defined and cannot always be answered by a definite 
yes or no.
Traditional mathematics avoid such vagueness
by precise definitions and clear notional scopes.
This concept undoubtedly is very successful, but it simply excludes vague 
notions such as ``big.''

The problem is not our principal inability to precise ``big'' or ``cold.''
It would be quite easy to determine that a man is ``big'' if he has a height
greater than 1.90\,m, and that a room is ``cold'' if its temperature
is smaller than 17${}^\circ$C.
The actual problem is that the boundary between ``big'' and ``normal''
is not exactly at 1.90\,m, or
the transition between ``cold'' and ``warm'' not abruptly at 17$^\circ$,
i.e., that there is a blurred transition  between these characteristics
in which they hold \emph{simultaneously}.
The traditional mathematical modelling is not able to represent this issue.
In 1965, 
Zadeh
therefore introduced the concept of fuzzy sets
\cite{Zadeh-1965}, leading to an entire mathematical branch,
in particular to fuzzy logic.

In a fuzzy logic $L$ a proposition $a$ usually is identified with a 
\emph{membership function}\index{membership function}
$\mu_A: X \to [0,1]$ for a subset $A \subset X$ of a given universe $X$.
The universe $X$ is an arbitrary set whose elements are
objects of the real world or its 
quantifiable (``measurable'') properties,
for instance a subset of $\mathbb{R}$ representing the height of a person or the 
temperature of a room, or a discrete subset of combinations of persons,
symptomes and diagnoses \cite[§7.1.13]{Hajek-1998}.
The value $\mu_A(x)$ 
specifies the degree of membership
which the element $x\in X$ has with respect to the 
subset $A \subset X$, 
i.e., the truth value of proposition $a$ for the element $x$.
This way, the membership function  $\mu_A$
generalizes the characteristic
function $\chi_A:X \to \{0,1\}$ of a set $A$,
\[
	\chi_A(x) = 
	\left\{ \begin{array}{ll} 
		0 & \mbox{if $x \in A$,} \\ 1 & \mbox{if $x \notin A$,}
	\end{array} \right.
\]
where ``$x \in A$'' is to be identified with the proposition
$\mu_A(x)$ $=$ ``the measured value $x$ is in $A$.''
In this paper, we directly identify the proposition $a(x)$ with its 
membership function, i.e.,
$a(x) = \mu_a(x).$
A fuzzy logic $L$ then is the set 
$L = \{a: X \to [0,1]\}$
of membership functions on a given set $X$.
The fuzzy logic is contradictory if and only if there exists at least
one proposition
$a \in L$ with $0<a(x)<1$ for some $x \in X$.
Usually, for a fuzzy logic $L$ the connectives $\wedge$ and
$\vee$ for all $g$, $f \in L$, 
are defined pointwise by
\begin{equation}
	f(x) \wedge g(x) 
	=
	\min( f(x), g(x)),
	\qquad 
	f(x) \vee g(x) 
	= 
	\max( f(x), g(x))
	.
	\label{eq-fuzzy-operations}
\end{equation}
Moreover, the constant functions
$0$ and $1 \in L$ are the universal bounds of $L$.

\subsection{t-norms and the derivation of negations}

\begin{definition}
	\cite{Hajek-1998,Klement-et-al-1999,Gottwald-Hajek-2005}
	A \emph{t-norm}\index{t-norm} (``triangular norm'') is a binary operation
	$*: [0,1]^2 \to [0,1]$ satisfying the following conditions
	for all $x$, $y$, $z \in [0,1]$:
	\begin{eqnarray}
		\mbox{\emph{(commutativity)}} \qquad &
		x * y = y * x
		\qquad
		\\
		\mbox{\emph{(associativity)}} \qquad &
		(x * y) * z = x * (y * z)
		\qquad
		\\
		\mbox{\emph{(monotony)}} \qquad &
		x * z \leqq y * z \ \mbox{ if } x \leqq y 
		\qquad
		\\
		\mbox{\emph{(boundary condition)}} \qquad &
		1 * x = x,
		\qquad
	\end{eqnarray}
\end{definition}

By the boundary condition we have especially $1 * 0 = 0$, i.e.,
by commutativity and monotony $0 * x = 0$.
In fuzzy set theory t-norms are used to model the intersection of
two fuzzy sets and therefore equivalently refer to the logical
term of conjunction. Consequently, t-norms may be used to derive
a logical structure.

\begin{lem}
	Let $*$ be a continuous t-norm. 
	Then for each pair $x$, $y \in [0,1]$, the element
	\begin{equation}
		x \to y \ := \ \sup \{ z \in [0,1] : x * z \leqq y \}
	\end{equation}
	is well-defined. It is called the \emph{residuum}\index{residuum} of the t-norm,
	or \emph{material implication}.
\end{lem}
\begin{proof}
	\cite[Lemma 2.1.4]{Hajek-1998}
\end{proof}

In this way, a given continuous t-norm uniquely implies a material implication
\cite[Lemma 2.1.4]{Hajek-1998}, 
and by the definition $x' := (x \to 0)$ also uniquely a negation.

In the context of logics, the t-norms are also called
\emph{fusion}\index{fusion} or \emph{multiplicative conjunction}\index{conjunction}.
In general, the fusion $f * g$ is the proposition which is false to the
degree that the sum of $f$ and $g$ is false,
whereas $f \to g$ is the proposition which is false to the extent that $f$ is 
truer than $g$.
It mirrors in the language of formulae the
behavior of the concatenation of premises $X$ in ``sequents'' 
$X \vdash A$ stating that $A$ can be derived from 
the structure $X$ (which may be a theory, e.g.),
by holding the ``introduction rule'' and the ``elimination rule''
\begin{equation}
	\frac
	{X \vdash A \quad Y \vdash B}
	{X;Y \vdash A * B}
	, \qquad \quad
	\frac
	{X \vdash A * B \quad Y(A;B) \vdash C}
	{Y(X) \vdash C}
\end{equation}
They mean that, if $X$ is a premise for $A$ and $Y$ a premise for $B$
then $X;Y$ is a premise for $A * B$; if on the other hand,
$X$ is enough for $A * B$ and $Y(A;B)$ is a premise for $C$ then
we can replace the reference to $A;B$ in $Y$ by a claim to $X$.
So, $A * B$ is the formula equivalent to the structure $A;B$.

We observe at this point that, in contrast to an implicative lattice
(cf.\ Def.~\ref{def-implicative-lattice}), the
residuum and therefore the negation is not determined by the lattice meet operation
$\wedge$ but by the additional fusion operation $*$.
The following continuous t-norms are important for fuzzy logics.

\begin{beispiel}
	\textbf{({\L}ukasiewicz logic)}%
	\index{Lukasiewicz logic}%
	\label{bsp-Lukasiewicz}
	The \emph{{\L}ukasiewicz t-norm} is defined as
	\begin{equation}
		f *_{\mbox{\footnotesize\L}} g = \max(f + g - 1,\ 0)
		.
		\label{eq-Lukasiewicz-fusion}
	\end{equation}
	It uniquely determines the \emph{{\L}u\-ka\-sie\-wicz implication}
	$\to_{\mbox{\footnotesize\L}}$ and the
	\emph{{\L}u\-ka\-sie\-wicz negation} $\neg_{\mbox{\footnotesize\L}}$
	as
	\begin{equation}
		f \to_{\mbox{\footnotesize\L}} g = \min(1 - f + g, 1),
		\qquad
		\neg_{\mbox{\footnotesize\L}} f = 1 - f
		.
		\label{eq-Luksasiewicz-implication-negation}
	\end{equation}
	Since for $0<f<1$ we have 
	$f \wedge \neg_{\mbox{\footnotesize\L}} f = \min(f, 1-f) > 0$
	as well as
	$f \vee \neg_{\mbox{\footnotesize\L}} f = \max(f, 1-f) < 1$,
	in the fuzzy logic $(L, \neg_{\mbox{\footnotesize\L}})$
	both the law of non-contradiction and tertium non datur do not hold.
	With $f = \neg_{\mbox{\footnotesize\L}} f$, 
	and Theorem \ref{theo-De-Morgan-fuzzy}, however,
	the conjunctive De Morgan law (\ref{conjunctive-De-Morgan}) is valid.
\end{beispiel}

\begin{beispiel}
	\textbf{(Gödel logic)}%
	\index{Gödel logic}%
	\label{bsp-Goedel}
	The \emph{Gödel t-norm} is defined as
	\begin{equation}
		f *_{\mbox{\footnotesize G}} g = \min(f,g) = f \wedge g
		.
		\label{eq-Goedel-fusion}
	\end{equation}
	It uniquely determines the \emph{Gödel implication}
	$\to_{\mbox{\footnotesize G}}$ and the
	\emph{Gödel negation} $\neg_{\mbox{\footnotesize G}}$
	as
	\begin{equation}
		f \to_{\mathrm{G}} g = 
		\left\{ \begin{array}{ll}
			1 & \mbox{if $f \leqq g$,} \\
			g & \mbox{if $f > g$.}
		\end{array} \right.
		\qquad
		\neg_{\mbox{\footnotesize G}} f 
		= \delta(f) 
		= \left\{ \begin{array}{ll}
			1 & \mbox{if $f=0$,} \\
			0 & \mbox{if $f > 0$.}
		\end{array} \right.
		\label{eq-Goedel-implication-negation}
	\end{equation}
	Since for $0 \leqq f \leqq 1$ we have 
	$f \wedge \neg_{\mbox{\footnotesize G}} f = \min(f, \delta(f)) = 0$
	the fuzzy logic $(L, \neg_{\mbox{\footnotesize G}})$ is
	non-contradictory and therefore a logic.
	However, by
	$f \vee \neg_{\mbox{\footnotesize G}} f = \max(f, \delta(f)) < 1$
	for $0 < f < 1$, the tertium non datur does not hold in it.
	Moreover, with 
	$\neg_{\mbox{\footnotesize G}} (f \wedge g)$ $=$ $\delta(\min(f,g))$
	and
	$\neg_{\mbox{\footnotesize G}} f \vee \neg_{\mbox{\footnotesize G}} g$
	$=$ $\max(\delta(f), \delta(g))$
	we have either
	$\delta(\min(f,g)) = \delta(0) = 1 = \max(\delta(f), \delta(g))$
	if $fg = 0$, or
	$\delta(\min(f,g)) = 0 = \max(\delta(f), \delta(g))$
	if $fg > 0$, and therefore
	the conjunctive De Morgan law (\ref{conjunctive-De-Morgan}) does hold.
\end{beispiel}

\begin{beispiel}
	\textbf{(Product logic)}%
	\index{product logic}%
	\label{bsp-product}
	The \emph{product t-norm} $*_\Pi$ is defined as
	\begin{equation}
		f *_{\Pi} g = f \cdot g
		.
		\label{eq-product-fusion}
	\end{equation}
	It uniquely determines the \emph{Goguen implication}
	$\to_{\Pi}$ and the
	\emph{Gödel negation} $\neg_{\mbox{\footnotesize G}}$
	as
	\begin{equation}
		f \to_{\Pi} g = 
		\left\{ \begin{array}{ll}
			1 & \mbox{if $f \leqq g$,} \\
			g/f & \mbox{if $f > g$.}
		\end{array} \right.
		\qquad
		\neg_{\mbox{\footnotesize G}} f 
		= \delta(f) 
		= \left\{ \begin{array}{ll}
			1 & \mbox{if $f=0$,} \\
			0 & \mbox{if $f > 0$.}
		\end{array} \right.
		.
		\label{eq-product-implication-negation}
	\end{equation}
	Since the negation of the product logic coincides with the
	Gödel negation (\ref{eq-Goedel-implication-negation}),
	fuzzy logic $(L, \neg_{\mbox{\footnotesize G}})$ is
	non-contradictory with the De Morgan laws holding.
\end{beispiel}

By Examples \ref{bsp-Goedel} and \ref{bsp-product} we observe that
the same negation allows for different consistent residuum structures.

There are considerable generalizations of the notion of a t-norm, for instance
to the case of two dimensions leading to the intuitionistic fuzzy
logic defined in the next example. On a lattice structure of a subset of 
$\mathbb{R}^2$, however, a
t-norm does not imply a negation, 
any negation may be compatible with it
\cite[§8.3.3]{Deschrijver-Kerre-2005}.

\begin{beispiel}
	\label{bsp-L*}
	\emph{(Triangular intuitionistic logic $L^*$)} \cite{Deschrijver-Kerre-2005}
	Let
	\begin{equation}
		L^* = \{x = (x_1, x_2) \in [0,1]^2 : x_1 + x_2 \leqq 1\}
		\label{eq-L*-def}
	\end{equation}
	denote the triangular surface in the plane with vertices $(0,0)$, $(1,0)$,
	and $(0,1)$, and let $\leq_*$ be defined for all $x$, $y \in L^*$ by
	\begin{equation}
		(x_1, x_2) \leq_* (y_1, y_2) 
		\quad \mbox{if and only if} \qquad
		x_1 \leqq y_1 \mbox{ and } x_2 \geqq y_2
		\label{eq-L*-leq-def}
	\end{equation}
	Then $(L^*, \leq_*)$ is a poset.
	\begin{figure}[htp]
	\centering
	\begin{footnotesize}
	\unitlength1ex
	\begin{picture}(20,20)
		\put(0,0){\includegraphics[width=20ex]{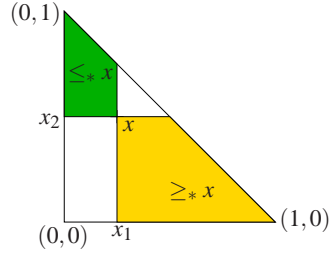}}
		\put(-0.2,19.0){\makebox(0,1)[r]{$(0,1)$}}
		\put( 0.0,-0.2){\makebox(0,0)[t]{$(0,0)$}}
		\put(20.2,-0.2){\makebox(0,1)[l]{$(1,0)$}}
		\put( 5.6, 9.6){\makebox(0,0)[tl]{$x$}}
		\put(-0.2, 9.6){\makebox(0,0)[r]{$x_2$}}
		\put( 5.6,-0.2){\makebox(0,0)[t]{$x_1$}}
		\put( 4.6,14.0){\makebox(0,0)[r]{$ \leq_* x$}}
		\put(12.0, 3.0){\makebox(0,0){$ \geq_* x$}}
	\end{picture}
	\end{footnotesize}
	\caption{\label{fig-triangular-lattice}\footnotesize
		The poset $(L^*,\leq_*)$
		of Example \ref{bsp-L*}, and the subsets of points comparable to an $x \in L^*$.
	}
	\end{figure}
	The partial order in $L^*$ is illustrated in Figure \ref{fig-triangular-lattice}.
	For a given point $x \in L^*$, the white regions contain the points
	$y \in L^*$ which cannot be compared to $x$ by the partial order
	$\leq_*$, whereas the light gray area consists of the points 
	$y \geq_* x$ and the upper dark gray area the points $y \leq_* x$.
	With the definitions
	\begin{equation}
		x \wedge y = (\min[x_1, y_1], \max[x_2, y_2]),
		\qquad
		x \vee y = (\max[x_1, y_1], \min[x_2, y_2])
	\end{equation}
	$L^*$ is a lattice with universal bounds 
	$0_* = (0,1)$ and $1_* = (1,0)$.
	Moreover, it is distributive since in each component the order is total
	\cite[§I.6]{Birkhoff-1973}.
	The standard negation on $L^*$ is given by
	\begin{equation}
		(x_1, x_2)' = (x_2, x_1).
	\end{equation}
	It is involutive and satisfies both De Morgan laws, as is directly verified.
	Since for $x \in L^*$, with $0 < x_i < 1$ for at least one $i \in \{1,2\}$,
	we have
	\begin{equation}
		x \wedge x' = (\min[x_1, x_2],\ \max[x_1, x_2]) >_* 0_*
	\end{equation}
	and
	\begin{equation}
		x \vee x' = (\max[x_1, x_2],\ \min[x_1, x_2]) <_* 1_*,
	\end{equation}
	$L^*$ does neither satisfy the law of non-contradiction 
	nor the \emph{tertium non datur}.
\end{beispiel}



\subsection{Examples of finite and discrete fuzzy logics} 

\begin{beispiel}
	\label{bsp-Lukasiewicz-3}
	\textbf{(Finite {\L}ukasiewicz logics)}
	\cite[§8.2]{Restall-2000}
	For $n \in \mathbb{N}$,
	the \emph{{\L}ukasiewicz logic}\index{Lukasiewicz logic}
	$\mbox{\L}_{n+1}$ is a well-known $(n+1)$-valued logic
	with the constant truth values
	\begin{equation}
		\textstyle
		\mbox{\L}_{n+1} = 
		\left(
		\left\{0, \frac{1}{n}, \frac{2}{n}, \ldots, \frac{n-1}{n}, 1 \right\}, 
		\neg_{\mbox{\footnotesize\L}} \right).
	\end{equation}
	Especially for $n=2$ we obtain $\mbox{\L}_3 = \{0, \frac12, 1\}$, i.e.,
	\begin{equation}
		\begin{array}{c|c}
			f & \neg_{\mbox{\footnotesize\L}} f
			\\ \hline
			0       & 1       \\
			\frac12 & \frac12 \\
			1       & 0
		\end{array}
		\quad
		\begin{array}{c|c@{\ \ }c@{\ \ }c}
			\wedge  &    0    & \frac12 & 1
			\\ \hline
			0       &    0    &     0   & 0 \\
			\frac12 &    0    & \frac12 & \frac12 \\
			1       &    0    & \frac12 & 1
		\end{array}
		\quad
		\begin{array}{c|c@{\ \ }c@{\ \ }c}
			\vee    &    0    & \frac12 & 1
			\\ \hline
			0       &    0    & \frac12 & 1 \\
			\frac12 & \frac12 & \frac12 & 1 \\
			1       &    1    &     1   & 1
		\end{array}
		\quad
		\begin{array}{c|c@{\ \ }c@{\ \ }c}
			*_{\mbox{\footnotesize\L}} &    0    & \frac12 & 1
			\\ \hline
			0       &    0    &    0    & 0 \\
			\frac12 &    0    &    0    & \frac12 \\
			1       &    0    & \frac12 & 1
		\end{array}
		\quad
		\begin{array}{c|c@{\ \ }c@{\ \ }c}
			\to_{\mbox{\footnotesize\L}}  &    0    & \frac12 & 1
			\\ \hline
			0       &    1    &     1   & 1 \\
			\frac12 & \frac12 &     1   & 1 \\
			1       &    0    & \frac12 & 1
		\end{array}
	\end{equation}
	Notably, by the fusion $\frac12 *_{\mbox{\footnotesize\L}} \frac12 = 0$ 
	we do not have
	$f \leqq f *_{\mbox{\footnotesize \L}} f$ in general, i.e., the 
	fusion in $\mbox{\L}_3$ does not obey the ``law of weak contraction.''
	Since this is one of the so-called structural laws, $\mbox{\L}_3$ is a 
	``substructural logic''\index{substructural logic}
	\cite{Restall-2000}.
	$\mbox{\L}_3$ is equivalent to the $(f, n, t)$ fragment of BN$_4$
	(Def.~\ref{bsp-BN4-MO1}).
	We can extend the domain of propositions
	to $\mbox{\L}_{\mathbb{Q}} = [0,1] \cap \mathbb{Q}$, or
	to $\mbox{\L}_{\mathbb{R}} = [0,1]$.
\end{beispiel}

\begin{beispiel}
	\label{bsp-Goedel-3}
	\textbf{(Finite Gödel logics)}
	For $n \in \mathbb{N}$,
	the \emph{Gödel logic}\index{Gödel logic}
	$\mbox{G}_{n+1}$ is an $(n+1)$-valued logic
	with the constant truth values
	\begin{equation}
		\textstyle
		\mbox{G}_{n+1} = 
		\left(
		\left\{0, \frac{1}{n}, \frac{2}{n}, \ldots, \frac{n-1}{n}, 1 \right\}, 
		\neg_{\mbox{\footnotesize G}} \right).
	\end{equation}
	Especially for $n=2$ we obtain $\mbox{G}_3 = \{0, \frac12, 1\}$, i.e.,
	\begin{equation}
		\begin{array}{c|c}
			f & \neg_{\mbox{\footnotesize G}} f
			\\ \hline
			0       & 1       \\
			\frac12 & 0 \\
			1       & 0
		\end{array}
		\quad
		\begin{array}{c|c@{\ \ }c@{\ \ }c}
			\wedge  &    0    & \frac12 & 1
			\\ \hline
			0       &    0    &     0   & 0 \\
			\frac12 &    0    & \frac12 & \frac12 \\
			1       &    0    & \frac12 & 1
		\end{array}
		\quad
		\begin{array}{c|c@{\ \ }c@{\ \ }c}
			\vee    &    0    & \frac12 & 1
			\\ \hline
			0       &    0    & \frac12 & 1 \\
			\frac12 & \frac12 & \frac12 & 1 \\
			1       &    1    &     1   & 1
		\end{array}
		\quad
		\begin{array}{c|c@{\ \ }c@{\ \ }c}
			\to_{\mbox{\footnotesize G}}\!\!  
			        &    0    & \frac12 & 1
			\\ \hline
			0       &    1    &     1   & 1 \\
			\frac12 &    0    &     1   & 1 \\
			1       &    0    & \frac12 & 1
		\end{array}
	\end{equation}
	Note that the residuum $*_{\mbox{\footnotesize G}}$ equals the
	meet operation $\wedge$ of the lattice.
	$\mbox{G}_3$ is also known as the
	Heyting lattice\index{Heyting lattice} $\mbox{H}_3$\index{H$_3$}
	\cite[§8.46]{Restall-2000}.
\end{beispiel}

\begin{beispiel}
	\textbf{(Sugihara models of RM$_{\mathbf{2n+1}}$)}
	For $n\in\mathbb{N}$, the logic RM$_{2n+1}$\index{RM$_{2n+1}$}\index{Sugihara model}
	(where RM stands for ``relevant logic with mingle'' \cite[§2.7]{Restall-2000})
	is given by
	\begin{equation}
		\mathrm{RM}_{2n+1}
		=
		\left(
			\{-n, -n-1, \ldots, -1, 0, 1, 2, \ldots, n\},
			\ -
		\right)
	\end{equation}
	with the negation $a'=-a$, and fusion and material implication defined as
	\begin{equation}
		a * b = 
		\left\{ \begin{array}{ll}
			a \wedge b & \mbox{if $a \leqq -b$,} \\
			a \vee   b & \mbox{if $a > -b$,}
		\end{array} \right.
		\qquad
		a \to b = 
		\left\{ \begin{array}{ll}
			-a \vee   b & \mbox{if $a \leqq b$,} \\
			-a \wedge b & \mbox{if $a < b$.} \\
		\end{array} \right.
	\end{equation}
	Fusion is communitative and associative, with identity 0, and
	we have $a * a = a$ for all $a$ $\in$ RM$_{2n+1}$.
	Moreover, the negation satisfies the De Morgan laws, but
	by $0 \vee (-0) = 0$,
	\emph{tertium non datur} does not hold.
	For instance, $n=1$ yields the three-valued logic
	RM$_3$ $=$ $(\{-1,0,1\}, -)$
	where
	\begin{equation}
		\begin{array}{r|r}
			a & - a
			\\ \hline
			-1 &  1 \\
			 0 &  0 \\
			 1 & -1
		\end{array}
		\qquad
		\begin{array}{r|rrr}
			\wedge & -1 &  0 &  1
			\\ \hline
			-1     & -1 & -1 & -1 \\
			 0     & -1 &  0 &  0 \\
			 1     & -1 &  0 &  1
		\end{array}
		\qquad
		\begin{array}{r|rrr}
			\vee    & -1 &  0 &  1
			\\ \hline
			-1      & -1 &  0 &  1 \\
			 0      &  0 &  0 &  1 \\
			 1      &  1 &  1 &  1
		\end{array}
	\end{equation}
	and
	\begin{equation}
		\begin{array}{r|rrr}
			* & -1 &  0 &  1
			\\ \hline
			-1    & -1 & -1 & -1 \\
			 0    & -1 &  0 &  1 \\
			 1    & -1 &  1 &  1
		\end{array}
		\qquad
		\begin{array}{r|rrr}
			\to   & -1 &  0 &  1
			\\ \hline
			-1    &  1 &  1 &  1 \\
			 0    & -1 &  0 &  1 \\
			 1    & -1 & -1 &  1
		\end{array}
	\end{equation}
	Hence RM$_3$ is equivalent to the $(f, b, t)$ fragment of BN$_4$
	(Def.~\ref{bsp-BN4-MO1}).
	The model can be extended to the infinite-valued logic\index{RM logic} 
	RM $=$ $(\mathbb{Z},-)$, however as a lattice it has 
	no universal lower and upper bounds $O=-\infty$, $I=\infty$.
\end{beispiel}

\begin{beispiel}
	\label{bsp-fuzzy-temperature}
	Consider the three propositions $a$, $b$, $c: X \to [0,1]$ on the universe
	$X=(-273,15; \infty) \subset \mathbb{R}$, 
	$a(x)$ $=$ ``$x$ is cold,''
	$b(x)$ $=$ ``$x$ is warm,''
	$c(x)$ $=$ ``$x$ is hot.''
\begin{figure}[htp]
\centering
\begin{footnotesize}
	\unitlength.8ex
	\begin{picture}(30,10)
		\put(0.2,0){\includegraphics[width=24ex]{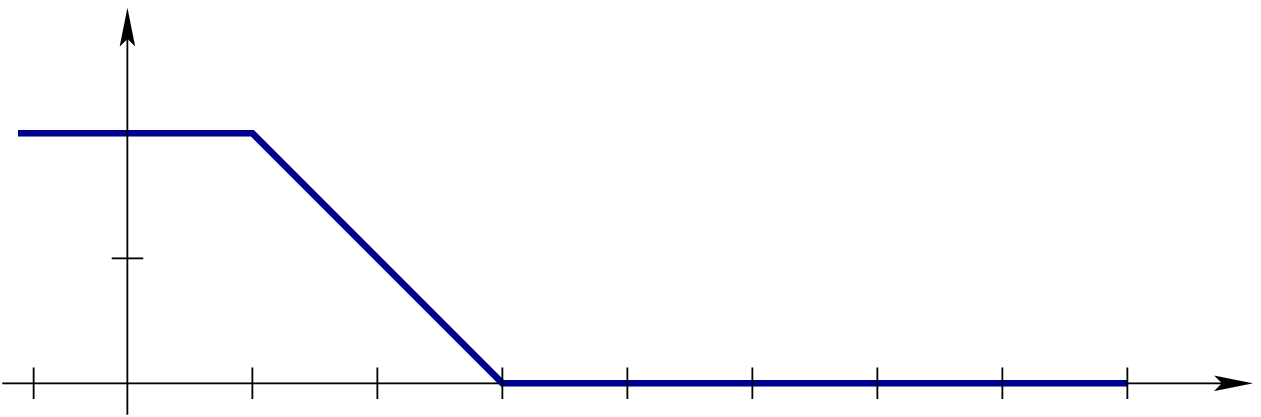}}
		\put( 2.5, 7.0){\makebox(0,0)[r]{$1$}}
		\put( 2.5, 4.1){\makebox(0,0)[r]{$0.5$}}
		\put(15.0, 6.0){\makebox(0,0)[b]{$a(x)$ ``cold''}}
		\put(29.0, 1.5){\makebox(0,0)[b]{$x$}}
		\put( 3.3, 0.0){\makebox(0,0)[t]{$0$}}
		\put( 9.2, 0.0){\makebox(0,0)[t]{$10$}}
		\put(15.4, 0.0){\makebox(0,0)[t]{$20$}}
		\put(21.0, 0.0){\makebox(0,0)[t]{$30$}}
		\put(29.0, 0.0){\makebox(0,0)[t]{${}^\circ$C}}
	\end{picture}
	\qquad
	\begin{picture}(30,10)
		\put(0.2,0){\includegraphics[width=24ex]{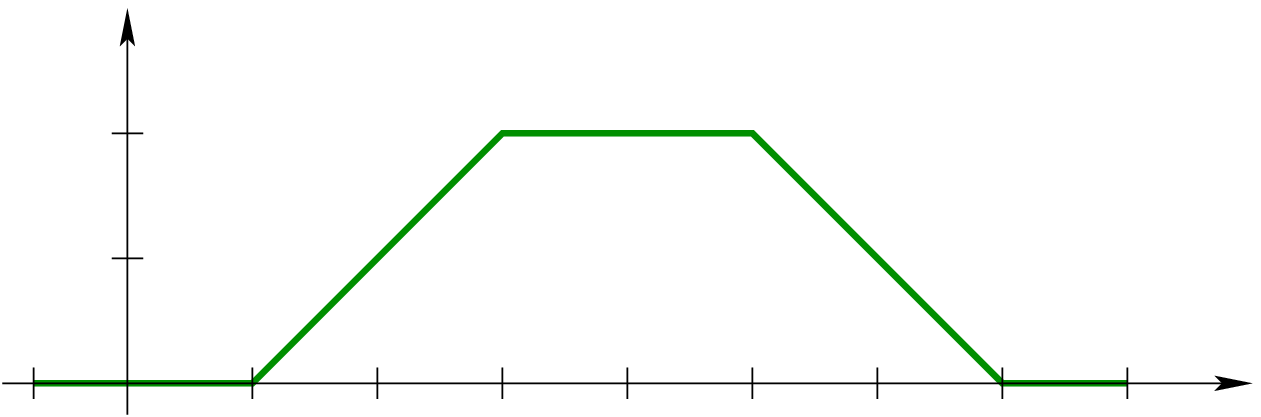}}
		\put( 2.5, 7.0){\makebox(0,0)[r]{$1$}}
		\put( 2.5, 4.1){\makebox(0,0)[r]{$0.5$}}
		\put(15.0, 7.0){\makebox(0,0)[b]{$b(x)$ ``warm''}}
		\put(29.0, 1.5){\makebox(0,0)[b]{$x$}}
		\put( 3.3, 0.0){\makebox(0,0)[t]{$0$}}
		\put( 9.2, 0.0){\makebox(0,0)[t]{$10$}}
		\put(15.4, 0.0){\makebox(0,0)[t]{$20$}}
		\put(21.0, 0.0){\makebox(0,0)[t]{$30$}}
		\put(29.0, 0.0){\makebox(0,0)[t]{${}^\circ$C}}
	\end{picture}
	\qquad
	\begin{picture}(30,10)
		\put(0.2,0){\includegraphics[width=24ex]{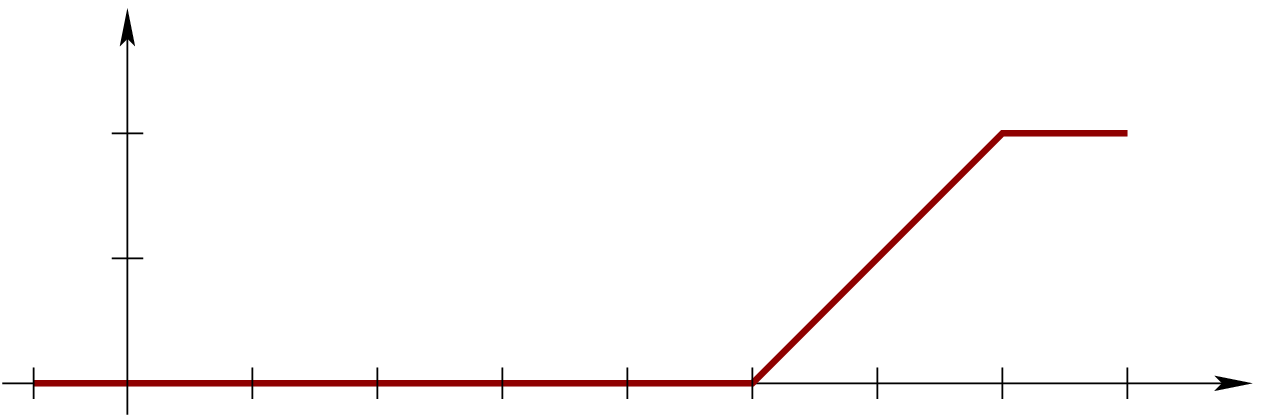}}
		\put( 2.5, 7.0){\makebox(0,0)[r]{$1$}}
		\put( 2.5, 4.1){\makebox(0,0)[r]{$0.5$}}
		\put(15.0, 6.0){\makebox(0,0)[b]{$c(x)$ ``hot''}}
		\put(29.0, 1.5){\makebox(0,0)[b]{$x$}}
		\put( 3.3, 0.0){\makebox(0,0)[t]{$0$}}
		\put( 9.2, 0.0){\makebox(0,0)[t]{$10$}}
		\put(15.4, 0.0){\makebox(0,0)[t]{$20$}}
		\put(21.0, 0.0){\makebox(0,0)[t]{$30$}}
		\put(29.0, 0.0){\makebox(0,0)[t]{${}^\circ$C}}
	\end{picture}
	\end{footnotesize}
	\caption{\label{fig-fuzzy-temperature}\footnotesize
		The membership functions 
		$a(x) = \mu_a(x)$,
		$b(x) = \mu_b(x)$,
		$c(x) = \mu_c(x)$.
	}
\end{figure}
	For a measurement outcome $x\in X$, they are defined as the functions
	(cf.\ Figure \ref{fig-fuzzy-temperature})
	\begin{eqnarray*}
		a(x) 
		& \hspace*{-1.5ex} = \hspace*{-1.5ex} &
		\left\{ \begin{array}{@{}cl}
			1 & \mbox{$x \leqq 5$,} \\
			\frac{15-x}{10} & \mbox{$5 < x \leqq 15$,} \\
			0 & \mbox{otherwise,}
		\end{array} \right.
		\quad 
		b(x)
		= 
		\left\{ \begin{array}{@{}cl}
			\frac{x-5}{10} & \mbox{$5 < x \leqq 15$,} \\
			1 & \mbox{$15 < x \leqq 25$,} \\
			\frac{35-x}{10} & \mbox{$25 < x \leqq 35$,} \\
			0 & \mbox{otherwise,}
		\end{array} \right.
		\\
		c(x)
		& \hspace*{-1.5ex} = \hspace*{-1.5ex} &
		\left\{ \begin{array}{@{}cl}
			0 & \mbox{$x \leqq 25$,} \\
			\frac{x-25}{10} & \mbox{$25 < x \leqq 35$,} \\
			1 & \mbox{otherwise.}
		\end{array} \right.
		\quad
	\end{eqnarray*}
	With the operations in (\ref{eq-fuzzy-operations}) 
	and the {\L}ukasiewicz negation (\ref{eq-Luksasiewicz-implication-negation})
	applied pointwise,
	i.e., $f'(x) = \neg_{\mbox{\footnotesize \L}} f(x)$,
	we have
	$a(x) \wedge b(x) \wedge c(x) = 0$ and
	$a(x) \vee b(x) \vee c(x) \geqq \frac12$
	for all $x\in X$, as is evident from
	Figure \ref{fig-fuzzy-temperature-2}.
\begin{figure}[htp]
\centering
\begin{footnotesize}
	\unitlength1.2ex
	\begin{picture}(30,10)
		\put(0.2,0){\includegraphics[width=36ex]{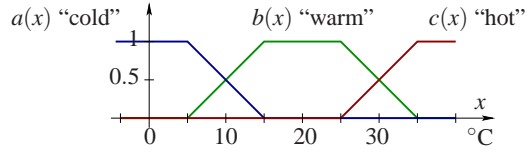}}
		\put( 2.5, 7.0){\makebox(0,0)[r]{$1$}}
		\put( 2.5, 4.1){\makebox(0,0)[r]{$0.5$}}
		\put( 1.0, 8.5){\makebox(0,0)[r]{$a(x)$ ``cold''}}
		\put(16.0, 8.5){\makebox(0,0){$b(x)$ ``warm''}}
		\put(25.0, 8.5){\makebox(0,0)[l]{$c(x)$ ``hot''}}
		\put(29.0, 1.5){\makebox(0,0)[b]{$x$}}
		\put( 3.3, 0.0){\makebox(0,0)[t]{$0$}}
		\put( 9.2, 0.0){\makebox(0,0)[t]{$10$}}
		\put(15.4, 0.0){\makebox(0,0)[t]{$20$}}
		\put(21.0, 0.0){\makebox(0,0)[t]{$30$}}
		\put(29.0, 0.0){\makebox(0,0)[t]{${}^\circ$C}}
	\end{picture}
	\end{footnotesize}
	\caption{\label{fig-fuzzy-temperature-2}\footnotesize
		The membership functions 
		$a(x)$, $b(x)$, $c(x)$, summarized into a single diagram.
	}
\end{figure}
	Moreover, the emerging fuzzy logic is contradictory since
	$f \wedge f' \ne 0$ for $f=a$, $b$, $c$ (Figure \ref{fig-fuzzy-temperature-atoms}).
\begin{figure}[htp]
\centering
\begin{footnotesize}
	\unitlength.8ex
	\begin{picture}(30,10)
		\put(0.2,0){\includegraphics[width=24ex]{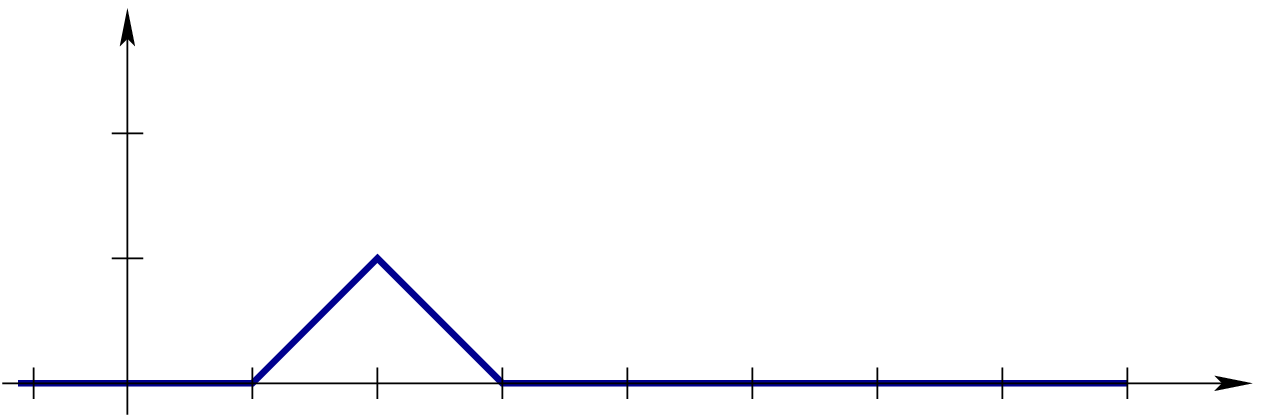}}
		\put( 2.5, 7.0){\makebox(0,0)[r]{$1$}}
		\put( 2.5, 4.1){\makebox(0,0)[r]{$0.5$}}
		\put(15.0, 7.0){\makebox(0,0){$(a\wedge a')(x)$}}
		\put(29.0, 1.5){\makebox(0,0)[b]{$x$}}
		\put( 3.3, 0.0){\makebox(0,0)[t]{$0$}}
		\put( 9.2, 0.0){\makebox(0,0)[t]{$10$}}
		\put(15.4, 0.0){\makebox(0,0)[t]{$20$}}
		\put(21.0, 0.0){\makebox(0,0)[t]{$30$}}
		\put(29.0, 0.0){\makebox(0,0)[t]{${}^\circ$C}}
	\end{picture}
	\qquad
	\begin{picture}(30,10)
		\put(0.2,0){\includegraphics[width=24ex]{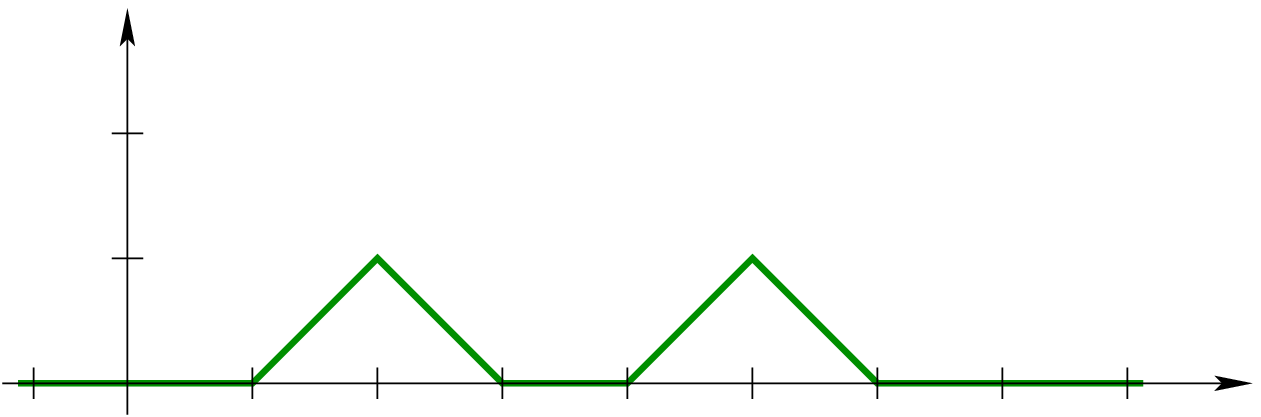}}
		\put( 2.5, 7.0){\makebox(0,0)[r]{$1$}}
		\put( 2.5, 4.1){\makebox(0,0)[r]{$0.5$}}
		\put(15.0, 7.0){\makebox(0,0){$(b\wedge b')(x)$}}
		\put(29.0, 1.5){\makebox(0,0)[b]{$x$}}
		\put( 3.3, 0.0){\makebox(0,0)[t]{$0$}}
		\put( 9.2, 0.0){\makebox(0,0)[t]{$10$}}
		\put(15.4, 0.0){\makebox(0,0)[t]{$20$}}
		\put(21.0, 0.0){\makebox(0,0)[t]{$30$}}
		\put(29.0, 0.0){\makebox(0,0)[t]{${}^\circ$C}}
	\end{picture}
	\qquad
	\begin{picture}(30,10)
		\put(0.2,0){\includegraphics[width=24ex]{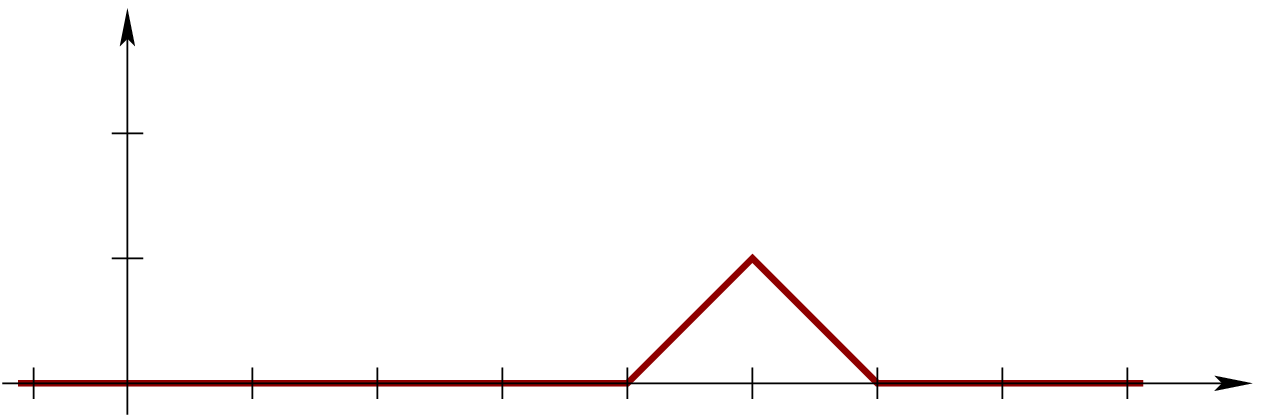}}
		\put( 2.5, 7.0){\makebox(0,0)[r]{$1$}}
		\put( 2.5, 4.1){\makebox(0,0)[r]{$0.5$}}
		\put(15.0, 7.0){\makebox(0,0){$(c \wedge c')(x)$}}
		\put(29.0, 1.5){\makebox(0,0)[b]{$x$}}
		\put( 3.3, 0.0){\makebox(0,0)[t]{$0$}}
		\put( 9.2, 0.0){\makebox(0,0)[t]{$10$}}
		\put(15.4, 0.0){\makebox(0,0)[t]{$20$}}
		\put(21.0, 0.0){\makebox(0,0)[t]{$30$}}
		\put(29.0, 0.0){\makebox(0,0)[t]{${}^\circ$C}}
	\end{picture}
	\end{footnotesize}
	\caption{\label{fig-fuzzy-temperature-atoms}\footnotesize
		The meets 
		$a\wedge a'$ (``cold and not cold''),
		$b\wedge b'$ (``warm and not warm''),
		$c\wedge c'$ (``hot and not hot'').
		They all are not identical to 0, thus the fuzzy logic they form is
		contradictory. Notice that $b\wedge b' = (a \wedge a') \vee (c \wedge c')$.
	}
\end{figure}
	For instance, a temperature of 10${}^\circ$C is ``cold and not cold'' 
	with a truth value $\frac{1}{2}$,
	and also ``warm and not warm'' with a truth value $\frac{1}{2}$.
	As a lattice, the fuzzy logic is depicted in Figure \ref{fig-fuzzy-lattice}.
	Notice that, e.g., $b\wedge b' = (a \wedge a') \vee (c \wedge c')$.
	It is non-orthomodular since, e.g.,
	$a \vee (a' \wedge b') = a' \wedge b' \ne b'$
	although $a < b'$.
	Note that it does not contain $O_6$ as a sublattice.
\begin{figure}[tb] 
\centering
	\begin{footnotesize}
	\unitlength1ex
	\begin{picture}(32,56)
		\put(3.0,3){\includegraphics[height=50ex]{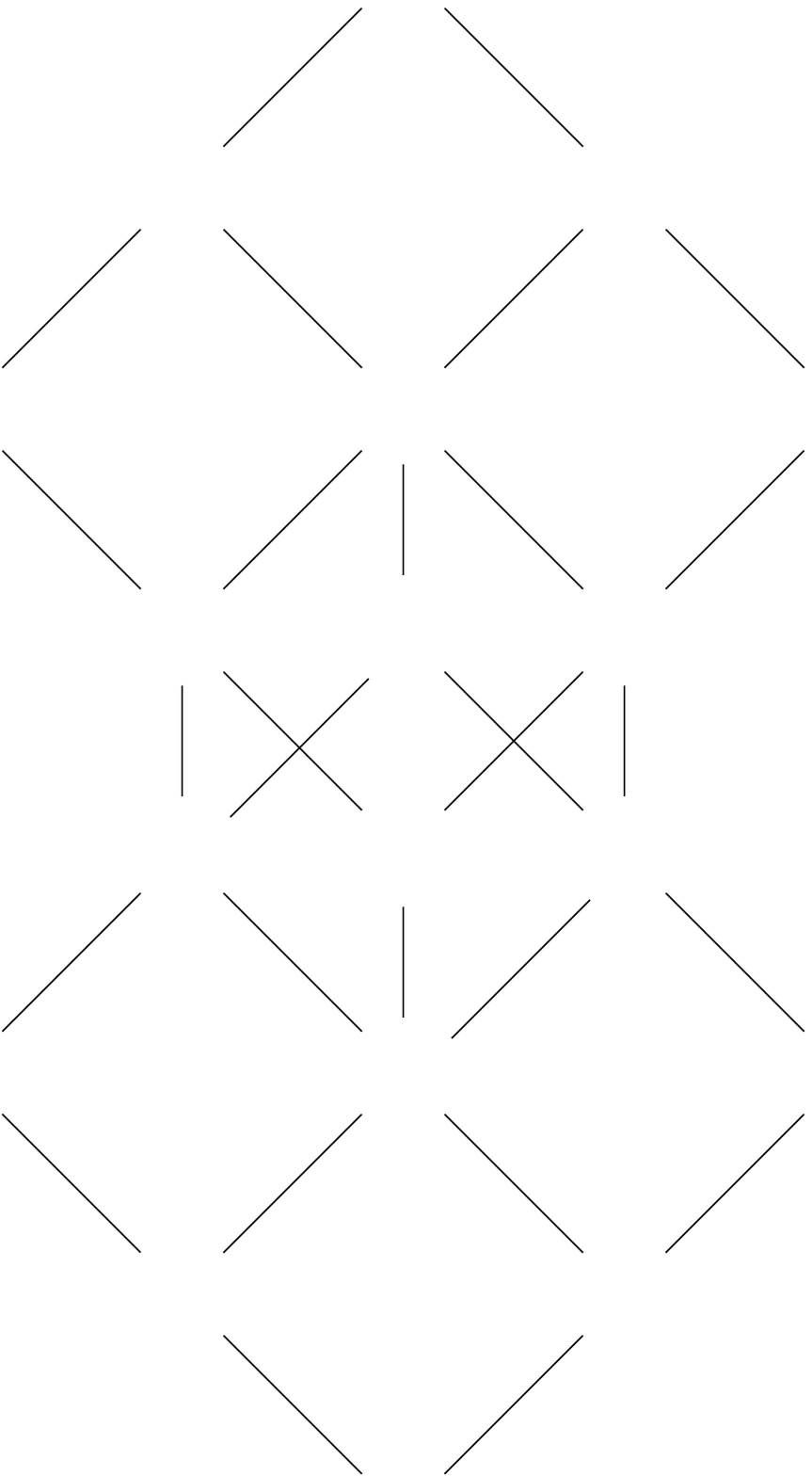}}
		\put(16.8,54.2){\makebox(0,0){\fbox{$1$}}}
		\put( 9.4,46.8){\makebox(0,0){\fbox{$c \vee c'$}}}
		\put(24.4,46.8){\makebox(0,0){\fbox{$a \vee a'$}}}
		\put( 1.8,39.2){\makebox(0,0){\fbox{$c'$}}}
		\put(17.0,39.2){\makebox(0,0){\fbox{$b \vee b'$}}}
		\put(31.8,39.2){\makebox(0,0){\fbox{$a'$}}}
		\put( 9.4,31.6){\makebox(0,0){\fbox{$a \vee b$}}}
		\put(16.8,31.6){\makebox(0,0){\fbox{$b'$}}}
		\put(24.2,31.6){\makebox(0,0){\fbox{$b \vee c$}}}
		\put( 9.2,24.2){\makebox(0,0){\fbox{$b' \wedge c'$}}}
		\put(16.8,24.2){\makebox(0,0){\fbox{$\phantom{'}\!\! b \phantom{'}\!\!$}}}
		\put(24.2,24.2){\makebox(0,0){\fbox{$a' \wedge b'$}}}
		\put( 1.8,16.7){\makebox(0,0){\fbox{$\phantom{'}\!\! a \phantom{'}\!\!$}}}
		\put(17.0,16.7){\makebox(0,0){\fbox{$b \wedge b'$}}}
		\put(31.8,16.7){\makebox(0,0){\fbox{$\phantom{'}\!\! c \phantom{'}\!\!$}}}
		\put( 9.4, 9.4){\makebox(0,0){\fbox{$a \wedge a'$}}}
		\put(24.4, 9.4){\makebox(0,0){\fbox{$c \wedge c'$}}}
		\put(16.8, 1.6){\makebox(0,0){\fbox{$0$}}}
	\end{picture}
	\qquad \qquad
	\begin{picture}(32,56)
		\put(3.0,3){\includegraphics[height=50ex]{./lattice-temperature-0}}
		\put(16.8,54.2){\makebox(0,0){{\includegraphics[width=7ex]{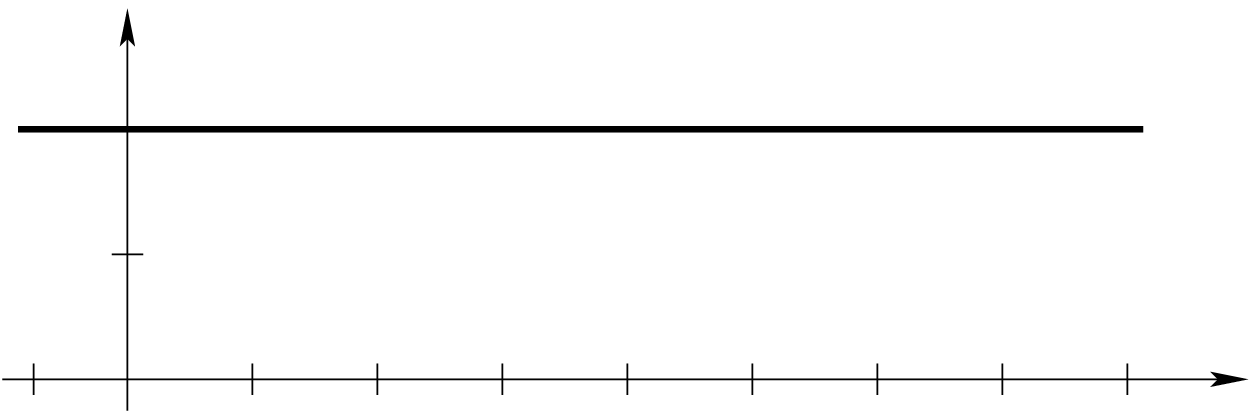}}}}
		\put( 9.4,46.8){\makebox(0,0){{\includegraphics[width=7ex]{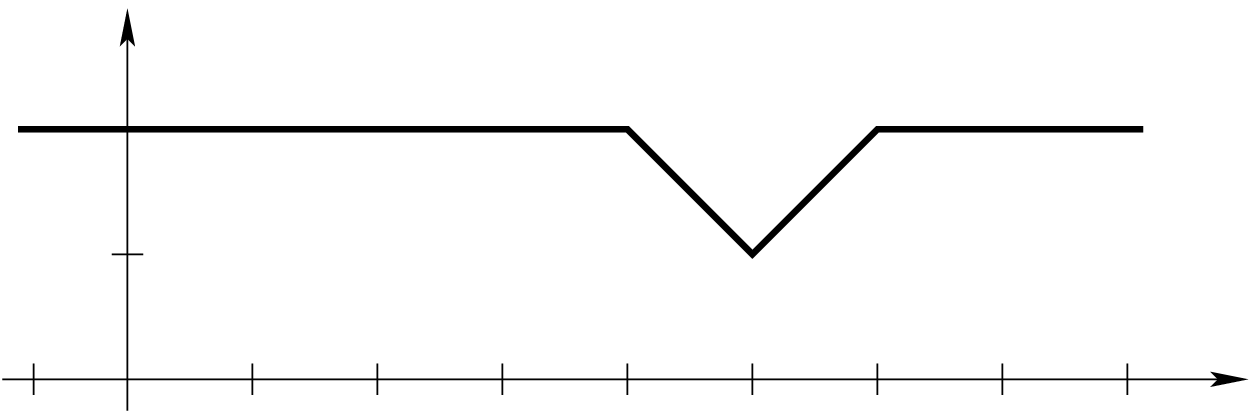}}}}
		\put(24.4,46.8){\makebox(0,0){{\includegraphics[width=7ex]{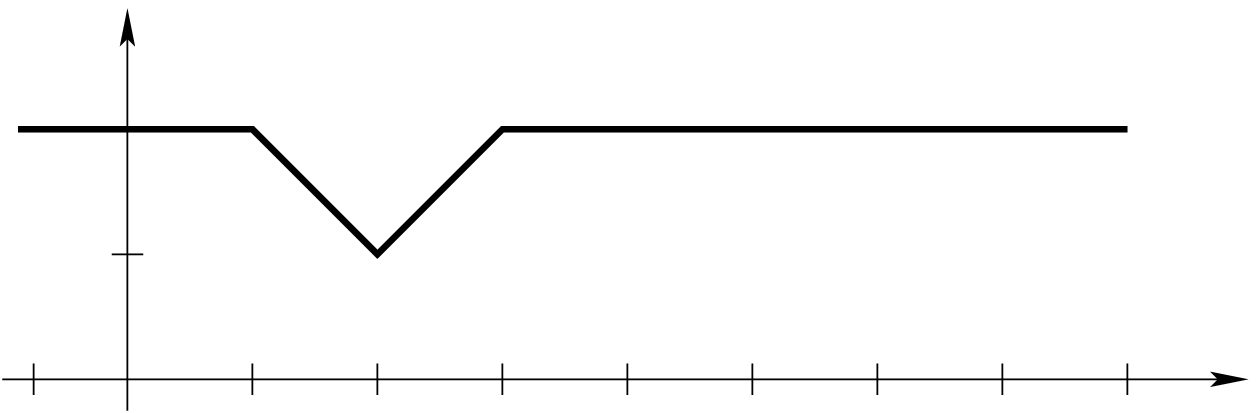}}}}
		\put( 1.8,39.2){\makebox(0,0){{\includegraphics[width=7ex]{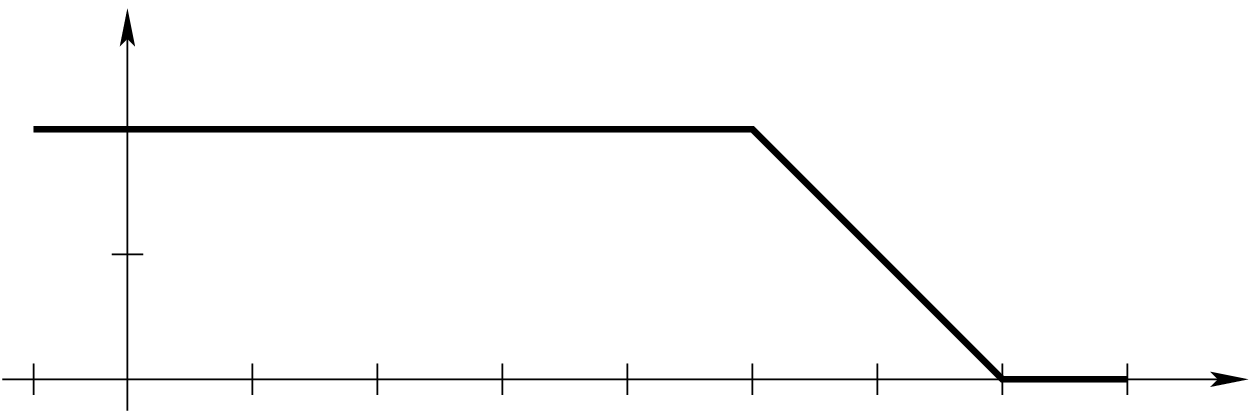}}}}
		\put(17.0,39.2){\makebox(0,0){{\includegraphics[width=7ex]{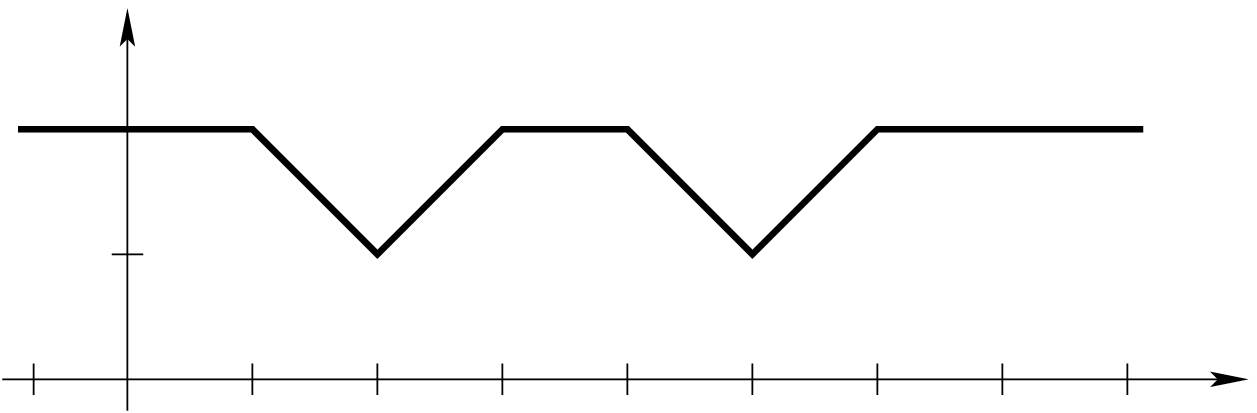}}}}
		\put(31.8,39.2){\makebox(0,0){{\includegraphics[width=7ex]{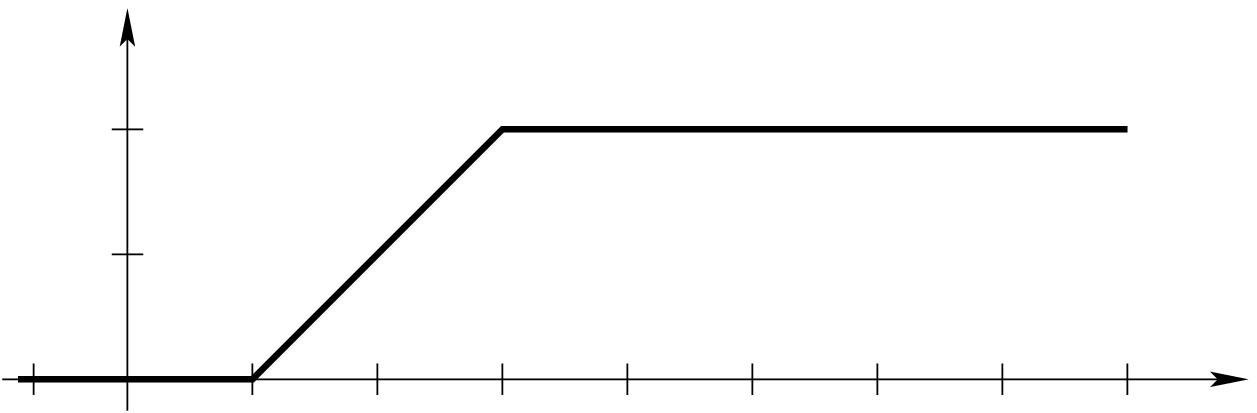}}}}
		\put( 9.4,31.6){\makebox(0,0){{\includegraphics[width=7ex]{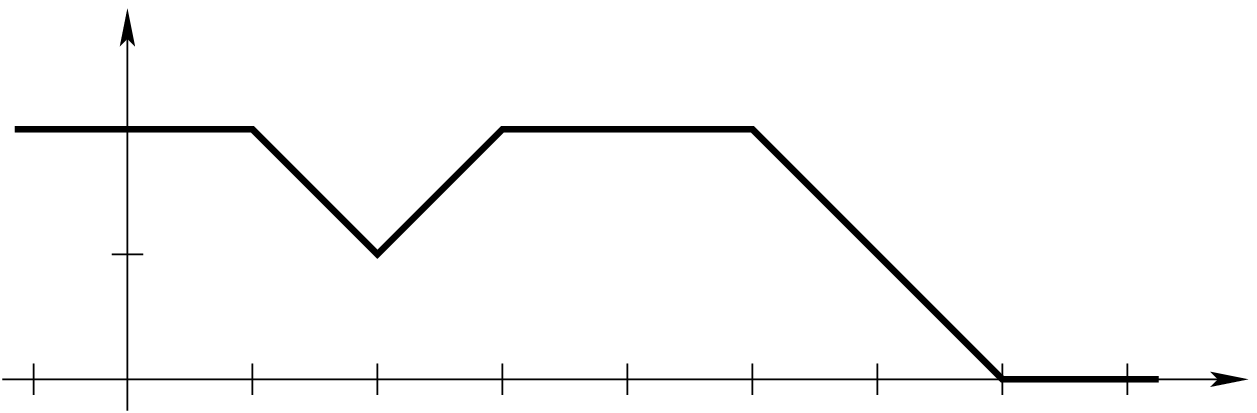}}}}
		\put(16.8,31.6){\makebox(0,0){{\includegraphics[width=7ex]{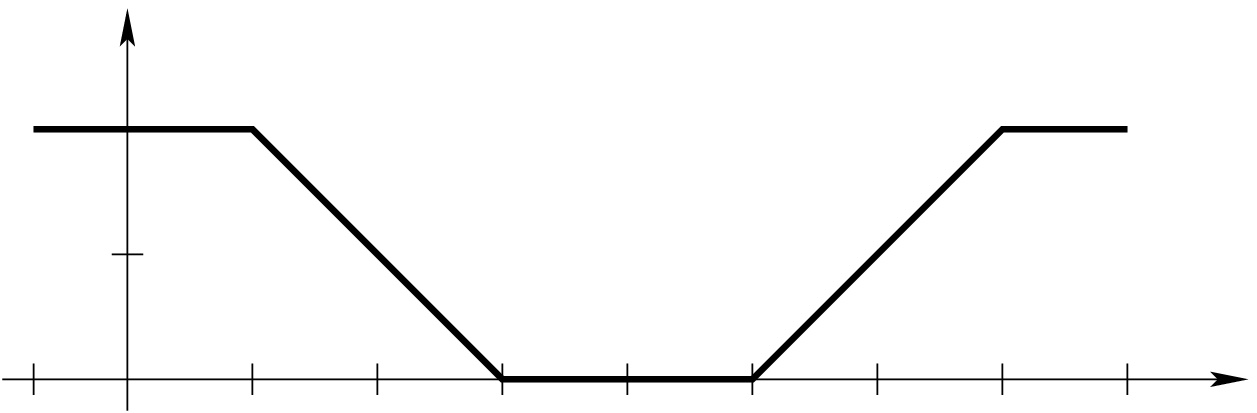}}}}
		\put(24.2,31.6){\makebox(0,0){{\includegraphics[width=7ex]{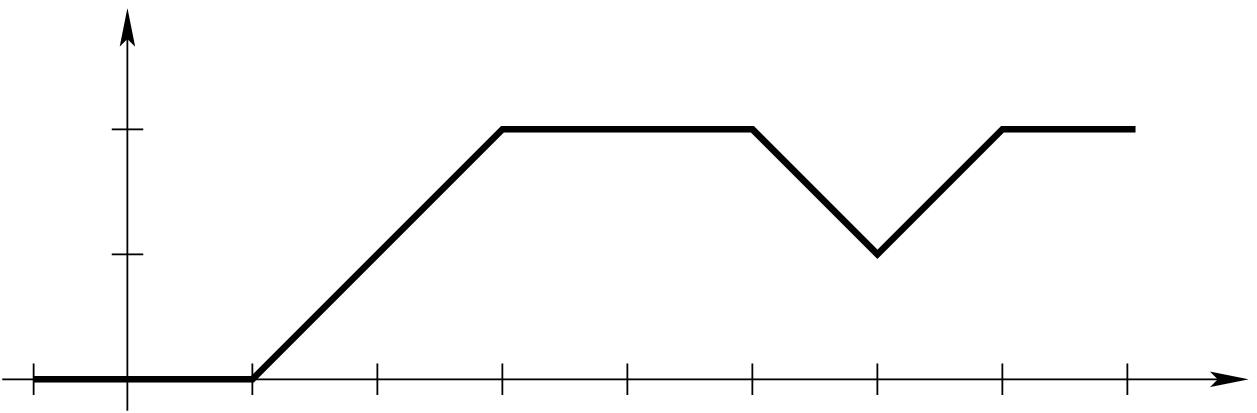}}}}
		\put( 9.2,24.2){\makebox(0,0){{\includegraphics[width=7ex]{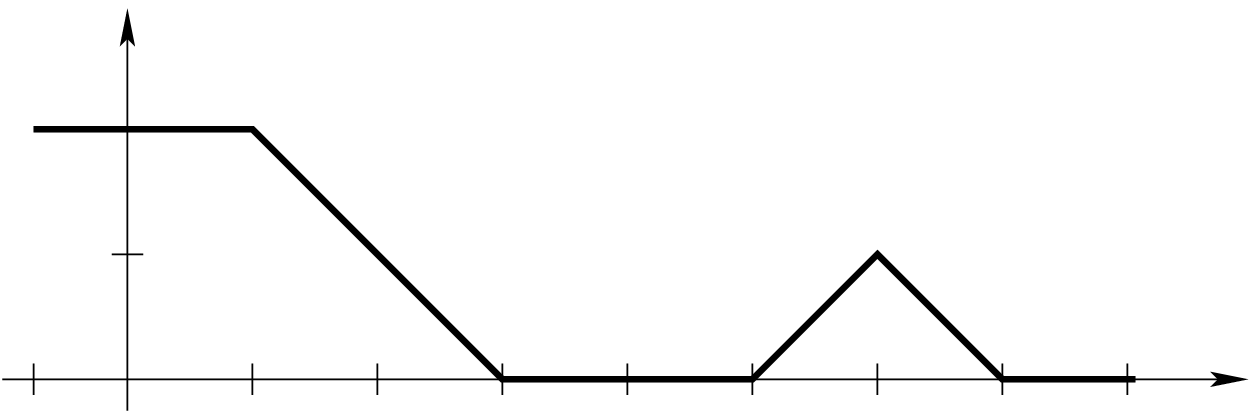}}}}
		\put(16.8,24.2){\makebox(0,0){{\includegraphics[width=7ex]{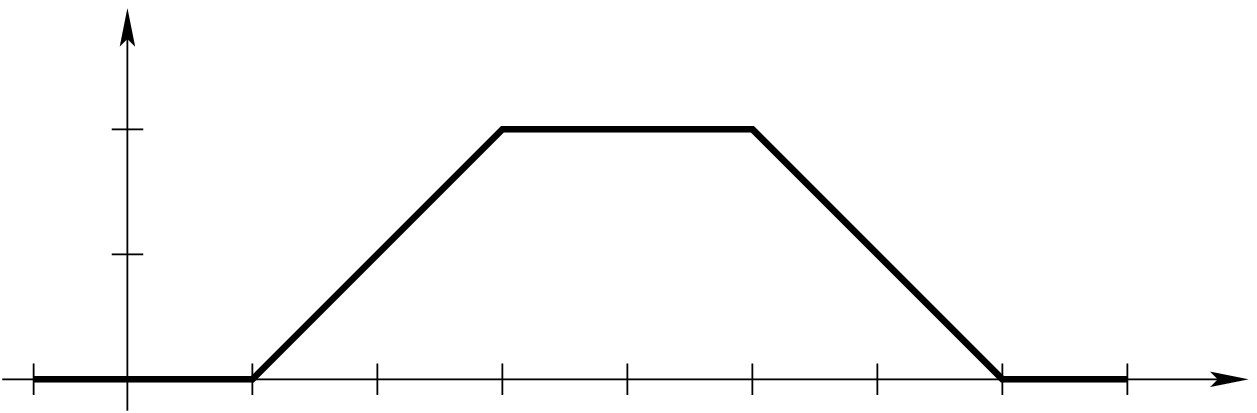}}}}
		\put(24.2,24.2){\makebox(0,0){{\includegraphics[width=7ex]{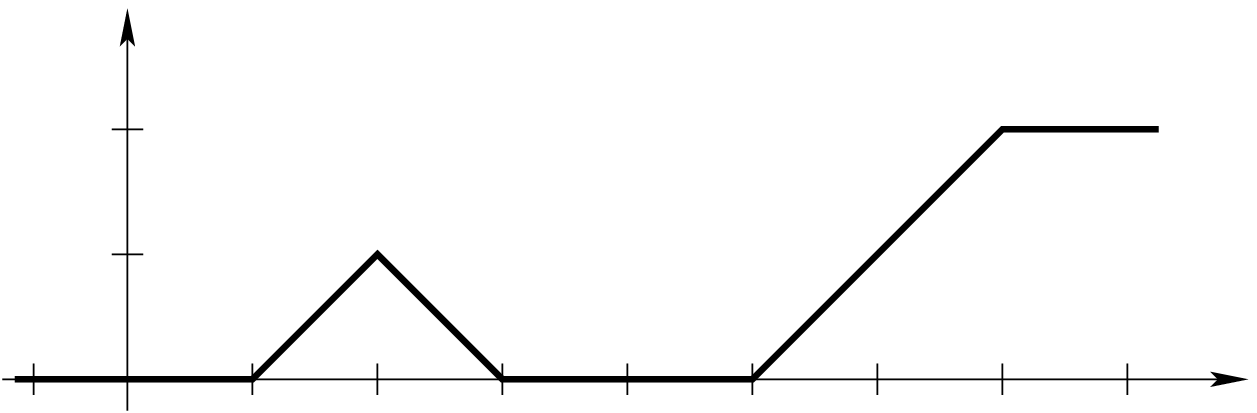}}}}
		\put( 1.8,16.7){\makebox(0,0){{\includegraphics[width=7ex]{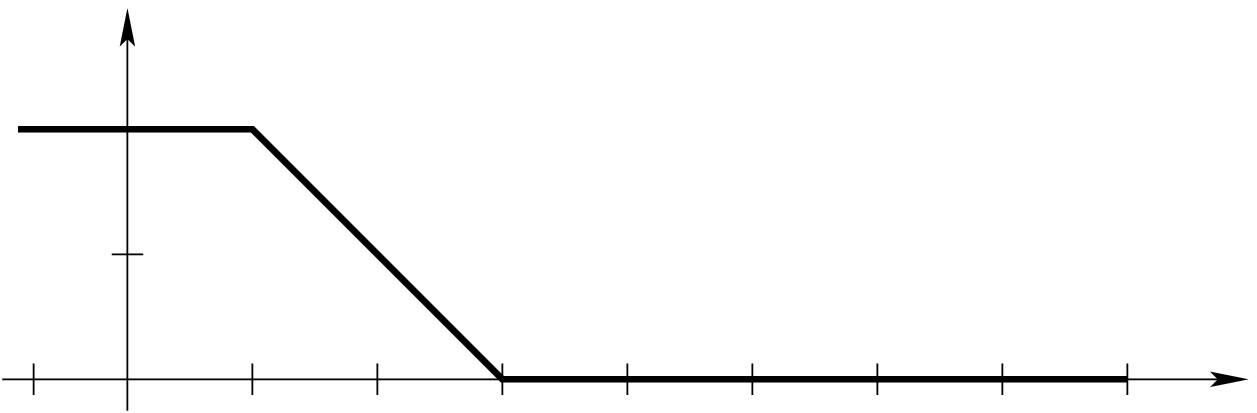}}}}
		\put(17.0,16.7){\makebox(0,0){{\includegraphics[width=7ex]{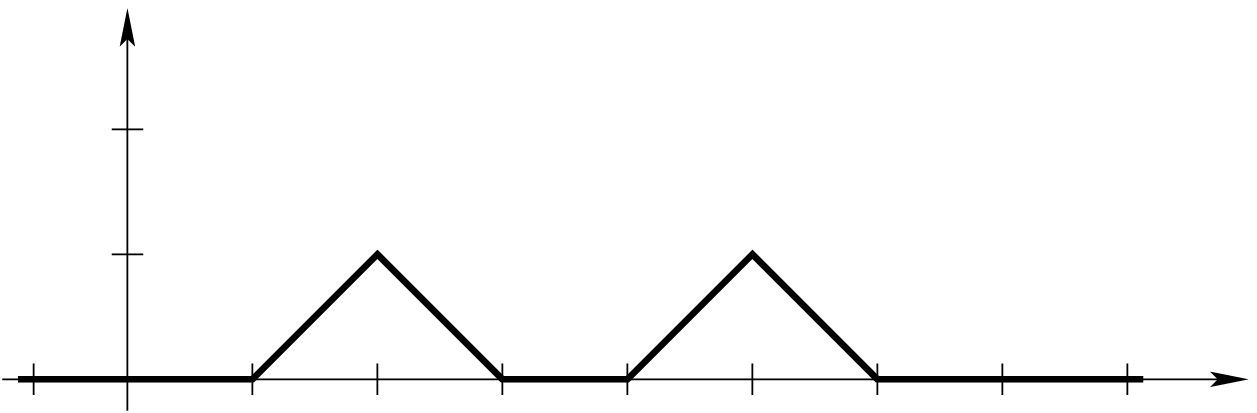}}}}
		\put(31.8,16.7){\makebox(0,0){{\includegraphics[width=7ex]{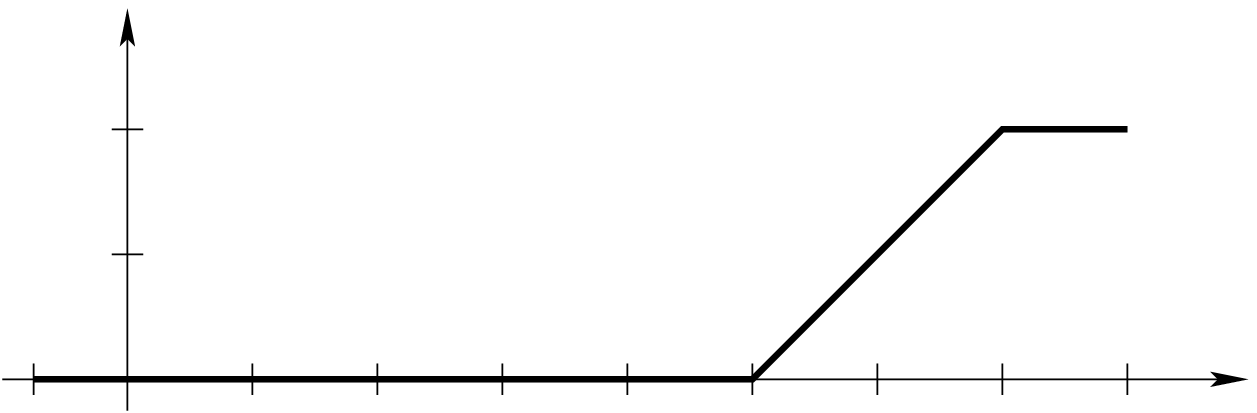}}}}
		\put( 9.4, 9.4){\makebox(0,0){{\includegraphics[width=7ex]{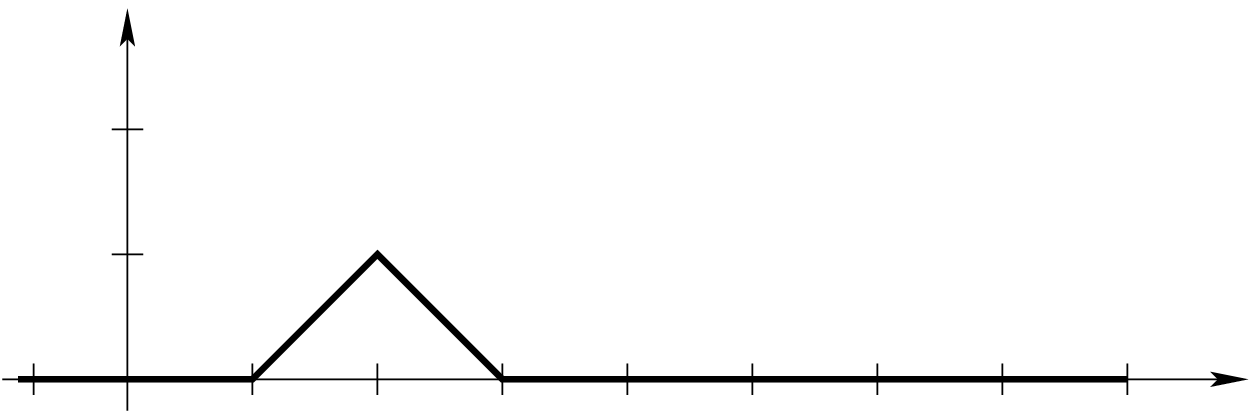}}}}
		\put(24.4, 9.4){\makebox(0,0){{\includegraphics[width=7ex]{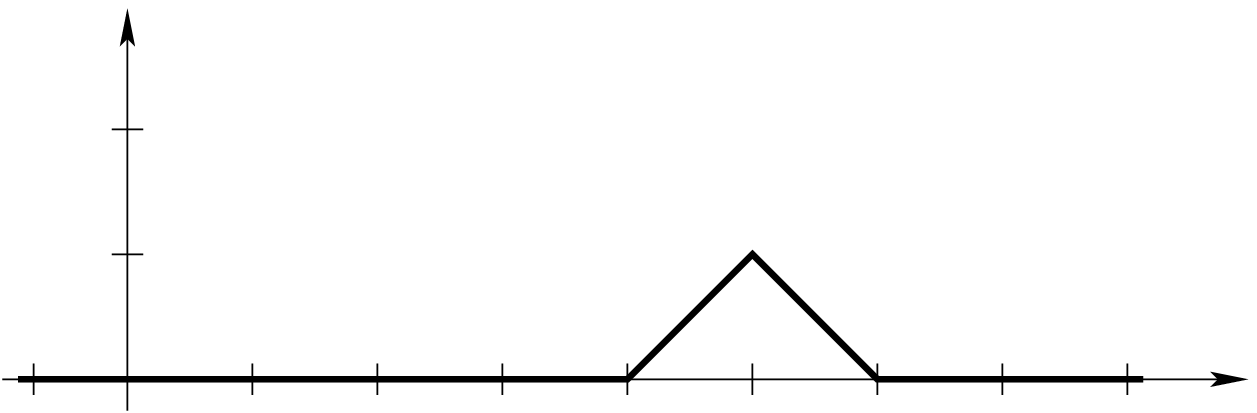}}}}
		\put(16.8, 1.6){\makebox(0,0){{\includegraphics[width=7ex]{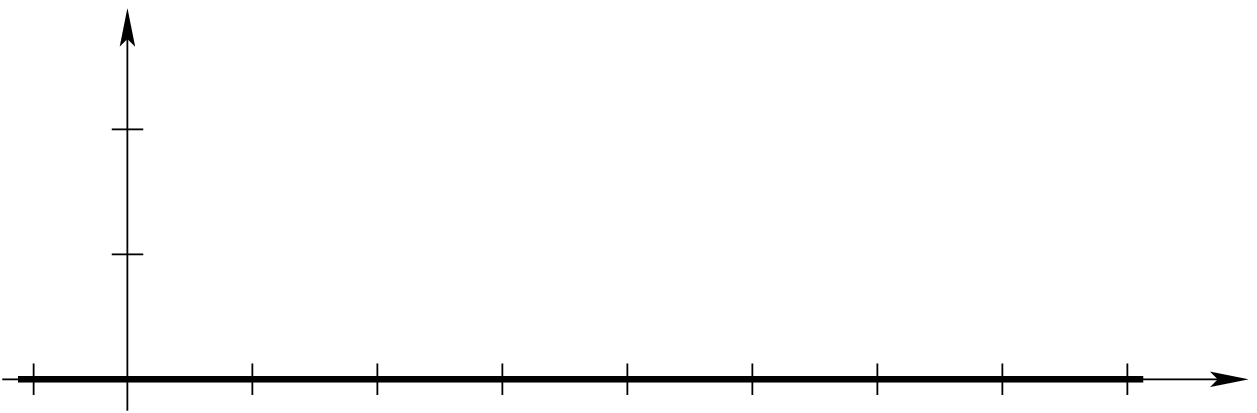}}}}
	\end{picture}
	\end{footnotesize}
	\caption{\label{fig-fuzzy-lattice} 
		The lattice of the fuzzy logic in Ex.~\ref{bsp-fuzzy-temperature},
		right hand the respective function graphs.
	}
\end{figure}
\end{beispiel}

\section{Quantum logics}
Historically, it was soon recognized that quantum mechanics involved departures from
classical Boolean logic, most strikingly
the Heisenberg uncertainty principle and the related principle of
the noncommutativity of physical observations.
Since the experimental verification of the distributive laws
for an 
algebra of attributes is based on the permutability
and repeatability of physical observations, distributivity
appeared inappropriate for a logic of projective quantum measurements. In addition, 
the Boolean concept of
negation had to be modified, since negating quantum attributes is based
on the orthogonality of subspaces of a Hilbert space and not on the
complements of subsets of a set.

In 1936, Birkhoff
and 
von Neumann
\cite{Birkhoff-von-Neumann-1936}
suggested to 
replace the distributivity condition by the weaker modularity condition. 
In general, however, the lattice of projections of a Hilbert space is 
not modular but only orthomodular.
It is modular\index{modular vs.\ orthomodular}
if and only if the Hilbert space is finite-dimensional,
and it is distributive\index{distributive vs.\ modular}
if and only if it is one-dimensional. 
Thus the requirement that the lattice of propositions be modular rules out 
the prototype quantum mechanical 
infinite-dimensional Hilbert space of Schrödinger's wavefunctions.
%
%
Nonetheless, the idea remained 
and became explicit
in Mackey's programme \cite{Mackey-1963} as a set of 
axioms, framing a conservative generalized probability theory, 
that underwrite the construction of a logic of experimental propositions, 
or, in his terminology, ``questions,'' having the structure of a 
$\sigma$-orthomodular poset. 

The observables of quantum mechanics are assumed to be Hermitian
operators acting on a Hilbert space $\mathscr{H}$. Therefore, if one
defines an observable \emph{attribute} of a quantum-mechanical state
$|\psi\rangle$ as a statement of the form that the observation $M$ on 
$|\psi\rangle$ will give a value $\lambda$ in a Borel set 
$S\subset\mathbb{R}$, then $|\psi\rangle$ has the attribute $(M,S)$
with probability 1 if and only if $|\psi\rangle$ is in a 
\emph{closed subspace} $X$ of $\mathscr{H}$.
If $S'$ is the complement of $S$ (as the set $S'=\mathbb{R} \setminus S$),
then the \emph{negation} of $(M,S)$ is certain if and only if
$|\psi\rangle$ $\in$ $X^\bot$, where $X^\bot$ denotes the orthogonal complement of $X$.

If every such pair $(M,S)$ is \emph{assumed} to correspond to an observable
attribute, then it follows that any closed subspace $\mathscr{H}$ of Hilbert
space represents an attribute observable ``with certainty'' if and only if
$|\psi\rangle$ $\in$ $\mathscr{H}$. 


\begin{axiom}[Mackey's Axiom VII]
The partially ordered set of all propositions \emph{(``questions'')} 
about a quantum system is 
isomorphic to the partially ordered set $L(\mathscr{H})$ of all closed subspaces of a 
separable, infinite dimensional Hilbert space.
\end{axiom}

\noindent
For Mackey, the outstanding problem was to 
explain \emph{why} the poset of questions
ought to be isomorphic to $L(\mathscr{H})$.

\begin{quote}
``This axiom has rather a different character from Axioms I through VI. 
These all had some degree of physical naturalness and plausibility. 
Axiom VII seems entirely ad hoc. Why do we make it? Can we justify making it?
Ideally, one would like to have a list of physically plausible assumptions 
from which one could deduce Axiom VII. Short of this one would like a list 
from which one could deduce a set of possibilities for the structure \ldots\
all but one of which could be shown to be inconsistent with suitably planned 
experiments.''
\cite[pp. 71--72]{Mackey-1963}
\end{quote}

This problem is still unresolved.
Even worse, rather natural examples of composite quantum systems
\cite{Aerts-1982, Foulis-Randall-1981} are known 
to violate orthomodularity or an equivalent structure, although 
in each of their single subsystems it does hold true.
These examples put serious doubts on the universality of the original Mackey programme.
The question of how to determine logics of general \emph{composite} quantum systems
is under current research 
\cite{Aerts-van-Steirteghem-2000,Ptak-Weber-2001,Soler-1995,Wilce-2006}.

\subsection{Subspaces in Hilbert space}
Observe that $|\psi\rangle \in X_1 \cap X_2$ for two closed subspaces 
$X_1$, $X_2$ $\subset$ $\mathscr{H}$ means that successive measurements of a system
in the state $|\psi\rangle$ will \emph{certainly} verify the predictions
$(H_1, S_1)$ and $(H_2,S_2)$ corresponding to $X_1$ and $X_2$, respectively.
Thus, like the orthocomplement $X^\bot$, $X_1 \wedge X_2$ has a simple and direct
physical meaning.

The lattice of projections of a Hilbert space is in fact a complete lattice, 
i.e., arbitrary sets of propositions, whether mutually orthogonal or not, 
countable or not, have well-defined suprema and infima.

Throughout this section, $\mathscr{H}$ will denote an arbitrary separable
Hilbert space.
\begin{definition}
	A set $\mathscr{M} \subset \mathscr{H}$ is a 
	\emph{linear manifold}\index{linear manifold}
	in $\mathscr{H}$ if for all $x$, $y$ $\in$ $\mathscr{M}$ and all
	constants $z\in \mathbb{C}$, both
	$x + y$ and $zx$ belong to $\mathscr{M}$.
\end{definition}
If $\mathscr{H}$ is finite-dimensional, a linear manifold is just a
``vector subspace.'' In the infinite-dimensional case, however, the
term ``subspace'' is reserved for certain kinds of manifolds
(Def.~\ref{def-subspace}).

\begin{definition}
	A set $S\subset\mathscr{H}$ is \emph{closed}\index{closed} if every
	Cauchy sequence in $S$ converges in norm to a vector in $S$.
\end{definition}

Let, for a general collection of sets $\mathscr{S}$, denote
$\cap \mathscr{S}$ the intersection of all sets of $\mathscr{S}$,
i.e., $\cap \mathscr{S}$ $=$ \{$x$: $x\in S$ for every $S\in\mathscr{S}$\}.

\begin{satz}
	\label{theo-closed-collection}
	If $\mathscr{S}$ is a collection of closed subsets of $\mathscr{H}$,
	then $\cap \mathscr{S}$ is a closed set in $\mathscr{H}$.
\end{satz}
\begin{proof}
	Suppose \{$x_k$\}$_k$ is a Cauchy sequence in $\cap \mathscr{S}$, which
	converges to a $x\in \mathscr{H}$. Then for all sets $C\in\mathscr{S}$,
	\{$x_k$\}$_k$ is a Cauchy sequence in $C$, and since $C$ is closed,
	$x$ must belong to $C$. Hence $x \in \cap\mathscr{S}$.
\end{proof}

\begin{definition}
	For a set $S\subset\mathscr{H}$, we define the \emph{closure of}\index{closure} 
	$S$ by
	\begin{equation}
		\mathrm{clos}(S) = \cap \{ C \subset\mathscr{H}: 
		\ S\subset C \mbox{, and $C$ closed}\}
	\end{equation}
\end{definition}
Observe from Theorem \ref{theo-closed-collection} that clos$(S)$ is a closed set.
From the definition of intersection, it is a subset of every closed set containing
$S$.

\begin{beispiel}
	\emph{(Counterexample)}
	Let $\mathscr{H}$ be a two-dimensional Hilbert space and
	$S=\{x\in\mathscr{H}:$ $\|x\| < 1\}$ its inner unit disc.
	Suppose 0 $\ne$ $x$ $\in$ $S$, and consider the sequence defined by
	\[ 
		x_k = \Big( \frac{1}{\|x\|} - \frac{1}{k} \Big) x
	\]
	for all $k\in \mathbb{N}$.
	Since $\|(1/\|x\| - 1/k)x\|$ $<$ $(1/\|x\|)\|x\|$ $=$ 1, we have
	$x_k \in S$ for all $k\in\mathbb{N}$. 
	Moreover, $x_k - x/\|x\|$ $=$ $-x/k$, i.e.,
	$\|(x_k - x/\|x\|) \|$ $=$ $\|x\|/k$ $\to$ 0 as $k \to \infty$.
	Hence $x_k$ is a Cauchy sequence converging to $\tilde x=x/\|x\|$, but
	$\tilde x \notin S$. That means, $S$ is not closed.
\end{beispiel}

\begin{definition}
	\label{def-subspace}
	A \emph{subspace}\index{subspace} in a Hilbert space $\mathscr{H}$
	is a closed linear manifold in $\mathscr{H}$.
\end{definition}

\begin{beispiel}
	\label{bsp-l-2}
	Let $\mathscr{H} = l_2$ be the vector space
	\begin{equation}
	\label{eq-l-2}
		l_2 = \{ \{x_k\}_k: \ x_k \in \mathbb{C} \mbox{ such that }
			\sum_{k=1}^\infty |x_k|^2 < \infty \},
	\end{equation}
	together with the inner product 
	$\langle x, y\rangle = \sum_{k=1}^\infty x_k y_k^*$.
	($l_2$ is a complex vector space, since if $x$, $y$ $\in$ $l_2$, then also
	$ax+by\in l_2$ for any constants $a$, $b$ $\in$ $\mathbb{C}$;
	moreover, the series $\sum_{k=1}^\infty x_k y_k^*$ converges, since
	$0 \leqq (|x_k - |y_k|)^2 = |x_k|^2 - 2|x_k||y_k| + |y_k|^2$, i.e.,
	$2|x_k y_k^*| = 2|x_k||y_k| \leqq |x_k|^2 + |y_k|^2$. Thus,
	$\sum_k |x_k||y_k^*| \leqq \frac12 (\sum_k |x_k|^2 + \sum_k |y_k|^2)$,
	and both series on the right converge because $x$, $y$ $\in$ $l_2$.)
	Define now the subset
	\begin{equation}
		S = \{ \{x_k\}_k \in l_2:\ x_k=0 \mbox{ for all but a finite number of } k\}.
	\end{equation}
	Then $S$ is a linear manifold, since for 
	$x=\{x_k\}_k$, $y=\{y_k\}_k$ $\in$ $S$, also the
	series $ax + by$ $=$ $\{ax_k + by_k\}$ for $a$, $b$ $\in$ $\mathbb{C}$
	are in $l_2$, since there still are
	only finitely many series member non-vanishing.
	However, $S$ is not closed, since for instance the series
	$\{x_k\}_k$ where each $x_k$ again is a series $x_k=\{x_k^j\}_j$ defined by
	\[ 
		x_k^j = \left\{ \begin{array}{ll}
			\sqrt{1/2^j} & \mbox{if $j\leqq k$,}\\
			0 & \mbox{if $j > k$}
		\end{array} \right.
	\]
	is a Cauchy sequence in $S$, $\{x_k\}_k \in S$, which converges in norm to
	$x=\{\sqrt{1/2^j}\}_j$, but $x\notin S$.
\end{beispiel}

\begin{beispiel}
	Suppose $[a,b] \subset \mathbb{R}$ is a closed real interval, and let
	$L^2(a,b)$ denote the Hilbert space of (Le\-bes\-gue-) square-integrable
	functions on $(a,b)$. Then the following sets are linear manifolds in
	$L^2(a,b)$.
	\begin{eqnarray}
		C[a,b]
		& \hspace*{-.5em} = \hspace*{-.5em} &
		\{ f:[a,b] \to \mathbb{C}:\ f \mbox{ is continuous}\}
		\\
		C^\infty[a,b]
		& \hspace*{-.5em} = \hspace*{-.5em} &
		\{ f:[a,b] \to \mathbb{C}:\ f \mbox{ is infinitely often differentiable}\}
	\end{eqnarray}
	If $f$ is continuous on $[a,b]$, then $|f|^2$ is also continuous and hence
	integrable, i.e., $C[a,b] \subset L^2(a,b)$. From this and from the fact that
	functions are continuous at all points where they are differentiable, we have
	\begin{equation}
		C^\infty[a,b] \subset C[a,b] \subset L^2(a,b).
	\end{equation}
	That $C^\infty[a,b]$ and $C[a,b]$ are linear manifolds in $L^2(a,b)$ follows
	from standard theorems in calculus. However, neither of these manifolds is
	closed in $L^2(a,b)$. On the other hand, the closure of $C[a,b]$ is
	$L^2(a,b)$.
\end{beispiel}

\begin{definition}
	If $S \subseteq \mathscr{H}$, we define the \emph{span}\index{span} of $S$
	by
	\begin{equation}
		\mathrm{span} S = \cap \{K \subseteq \mathscr{H}:\ 
		K \mbox{ is a subspace in $\mathscr{H}$ with } S \subseteq K\}.
	\end{equation}
	Two vectors $x$, $y$ $\in$ $\mathscr{H}$ of a Hilbert space $\mathscr{H}$
	are called \emph{orthogonal}\index{orthogonal}, in symbols
	$x \bot y$, if $\langle x, y\rangle = 0$.
	If $S\subset \mathscr{H}$, we define the \emph{orthogonal complement}%
	\index{complement! orthogonal -}\index{orthogonal complement}
	$S^\bot$ of $S$ as the subset 
	$S^\bot = \{x \in \mathscr{H}:$ $x\bot s$ for all $s\in S$\}.
	If $\mathscr{S}$ is a collection of subsets of $\mathscr{H}$, we write
	\begin{equation}
		\bigvee_{S \in \mathscr{S}} S := \mathrm{span}( \bigcup_{S\in\mathscr{S}} S )
		\mbox{,\qquad and \qquad}
		\bigwedge_{S \in \mathscr{S}} S := \mathrm{span}( \bigcap_{S\in\mathscr{S}} S )
		.
	\end{equation}
\end{definition}
For the intersection of two subspaces 
$X$, $Y$ in $\mathscr{H}$ we have simply
$X \wedge Y = X \cap Y$.
Cf.\ the finite-dimensional analog in Example \ref{bsp-vector-space}.

\begin{satz}
	\label{theo-ortho}
	Let be $S\subset \mathscr{H}$. Then the following hold:
	
	\begin{itemize}
	\item[(i)] $S^\bot \cap S = \{0\}$;
	
	\item[(ii)] $S^\bot$ is a subspace in $\mathscr{H}$ (even if $S$ is not);
	
	\item[(iii)] $S \subseteq T \subseteq \mathscr{H}$ $\Rightarrow$
	$T^\bot \subseteq S$;
	
	\item[(iv)] $S \subseteq (S^\bot)^\bot$; 
	
	\item[(v)] $S$ is a subspace in $\mathscr{H}$ $\Rightarrow$
	$(S^\bot)^\bot = S$;
	
	\item[(vi)] If $\mathscr{S}$ is a collection of subspaces in $\mathscr{H}$
	then 
	$(\bigvee_{S\in\mathscr{S}} S)^\bot$ $=$ $\bigcap_{S\in\mathscr{S}} S^\bot$,
	and
	$\bigvee_{S\in\mathscr{S}} S^\bot$ $=$ $(\bigcap_{S\in\mathscr{S}} S)^\bot$.
	\end{itemize}
\end{satz}
\begin{proof}
	(i) If $x \in S \cap S^\bot$, then $x \bot x$, which implies $x=0$.

	(ii) If $x$, $y\in S^\bot$ and $a$, $b\in\mathbb{C}$, then for all $s\in S$,
	$\langle ax+by, s\rangle$ $=$ $a\langle x,s\rangle$ $+$ $b\langle y.s\rangle$
	$=$ 0 $+$ 0 $=$ 0. Thus $S^\bot$ is a linear manifold.
	Suppose that $\{x_k\}$ is a sequence in $S^\bot$ that converges in norm to
	$x\in\mathscr{H}$. Then for all $s\in S$,
	$\langle x,s \rangle$ $=$ $\lim_{k\to \infty} \langle x_k, s\rangle$ by
	a fundamental Hilbert space property \cite[Theor.2.21B(i)]{Cohen-1989}.
	Since $\langle x_k, s\rangle$ $=$ 0 for all $k\in\mathbb{N}$, we have
	that $\langle x, s\rangle = 0$. So, $x\in S^\bot$, which shows that 
	$S^\bot$ is closed.
	
	(iii) Suppose $x\in T^\bot$. For every $s\in S$ we have $s\in T$, so that
	$x \bot s$. Hence, $x\in S^\bot$.
	
	(iv) Suppose $x \in S$. If $t \in S^\bot$, then $x \bot t$. Thus,
	$x\in (S^\bot)^\bot$.
	
	(v) -- (vi) see \cite[p.~123]{Cohen-1989}.
\end{proof}

\begin{satz}
	If $\mathscr{H}$ is a Hilbert space and $L(\mathscr{H})$ is the collection
	of all subspaces in $\mathscr{H}$, then
	$L(\mathscr{H})$ together with the set-inclusion $\subseteq$ and the
	complementation ${}^\bot$ is a 
	quantum logic with $0 = \{0\}$ and $1 = \mathscr{H}$.
\end{satz}
\begin{proof}
	If $K_1$, $K_2$ $\in$ $L(\mathscr{H})$, then $\{K_1$, $K_2\}$ has supremum
	$K_1 \vee K_2$ and infimum $K_1 \cap K_2$. So $L(\mathscr{H}, \subseteq)$ is
	a lattice. Clearly, \{0\} and $\mathscr{H}$ are, respectively, the least and
	greatest members of the lattice. By Theorem \ref{theo-ortho} (i), (v) and (vi), 
	$^\bot$ is a non-contradictory negation. So there remains to establish the 
	orthomodular identity (\ref{orthomodular-identity}).
	Suppose $J$, $K$ $\in$ $L(\mathscr{H})$ with $J \subseteq K$.
	We wish to show that $K=J\vee (J^\bot \cap K)$. Then there exist
	bases $B$ and $B_0$ for $J$ and $K$, respectively, with $B\subseteq B_0$
	\cite[Theor.~4.8]{Cohen-1989}.
	
	First, suppose $x\in K$. Then
	$x$ $=$ $\sum_{b\in B_0} \langle x, b\rangle b$
	$=$ $\sum_{b\in B} \langle x, b\rangle b$
	$+$ $\sum_{b\in B_0\setminus B} \langle x, b\rangle b$.
	Since $B_0\setminus B$ $\subseteq$ $J^\bot$, the second term belongs to 
	$J^\bot \cap K$. Since the first term belongs to $J$, we have that
	$x$ is a linear combination of vectors in
	$J \cup (J^\bot \cap K)$, so
	$x$ $\in$ $J \vee (J^\bot \cap K)$. This establishes that
	$K$ $\subseteq$ $J \vee (J^\bot \cap K)$.
	
	On the other hand, observe that since both $J$ and $J^\bot \cap K$ are
	subspaces of $K$, so is the span of their join.
\end{proof}

\begin{beispiel}
	Let $\mathscr{H}$ be a Hilbert space of dimension
	$\mathrm{dim}_{\mathbb{C}} \mathscr{H} \geqq 2$, and
	$x$, $y$ $\in$ $\mathscr{H}$ two nonzero orthogonal vectors in
	$\mathscr{H}$. Denote $X$, $Y$, $Z$ $\subset$ $\mathscr{H}$
	the one-dimensional subspaces spanned by them,
	$X=\mathrm{span}(\{x\})$, $Y=\mathrm{span}(\{y\})$, and
	$Z=\mathrm{span}(\{x+y\})$.
	In fact, $X$, $Y$, $Z$ $\in$ $L(\mathscr{H})$.
	Then we directly verify that
	$X \vee Y = 1$, $X \vee Z = 1$, 
	$X \wedge Y = 0$, $X \wedge Z = 0$, and
	$Z \wedge Y = 0$.
	Hence we have
	\begin{equation}
		X \vee (Y \wedge Z) = X,
		\qquad
		X \wedge (Y \vee Z) = X,
	\end{equation}
	and on the other hand,
	\begin{equation}
		(X \vee Y) \wedge (X \vee Z) = 1,
		\qquad
		(X \wedge Y) \vee (Y \wedge Z) = 0.
	\end{equation}
	Thus, the distributive laws (\ref{L6'}) and (\ref{L6''}) 
	are \emph{not} satisfied in the logic $L(\mathscr{H})$.
\end{beispiel}

\subsection{Quantum mechanics constructed from quantum logic}
\begin{beispiel}
	\label{bsp-projection-logic}
	Let $\mathscr{H}$ be a Hilbert space, and
	$\fett{P}$: $\mathscr{H} \to \mathscr{H}$ a self-adjoint operator with spectrum
	$\sigma_{\fett{\scriptsize P}} \subset \{0,1\}$. 
	Then $\fett{P}$ must be a \emph{projection}\index{projection},
	i.e., $\fett{P}^2 = \fett{P}$.
	Projections are in a bijective correspondence with the closed subspaces of
	$\mathscr{H}$: if $\fett{P}$ is a projection, its range ran($\fett{P}$) is closed, and
	any closed subspace is the range of a unique projection.
	If $u \in \mathscr{H}$ is a unit vector, then 
	$\langle \fett{P}u, u\rangle$ $=$ $\|\fett{P}u\|^2$ is the expected value of the 
	corresponding observable in the state represented by $u$.
	Since this is 0-1 valued, we can interpret $\|\fett{P}u\|^2$ as the probability
	that a measurement of the observable will produce the ``affirmative''
	answer 1. In particular, this affirmative answer will have probability
	1 if and only if $\fett{P}u=u$, i.e., $u$ $\in$ ran($\fett{P}$).
	
	We thus can impose on the set $L(\mathscr{H})$ of projections on $\mathscr{H}$
	the structure of a complete uniquely complemented lattice, defining
	\begin{equation}
		\fett{P} \leqq \fett{Q}
		\quad \mbox{if} \quad
		\textrm{ran}(\fett{P}) \subset \textrm{ran}(\fett{Q}),
		\qquad \mbox{and} \qquad
		\fett{P}' = 1 - \fett{P}
	\end{equation}
	(such that ran($\fett{P}'$) $=$ ran($\fett{P}$)$^\bot$). It is straightforward that
	$\fett{P} \leqq \fett{Q}$ just in the case $\fett{P}\fett{Q}=\fett{Q}\fett{P}=\fett{P}$.
	More generally, if $\fett{P}\fett{Q}$ $=$ $\fett{Q}\fett{P}$, then $\fett{P}\fett{Q}=\fett{P}\wedge \fett{Q}$; also in this case
	their join is given by $\fett{P}\vee \fett{Q}$ $=$ $\fett{P}+\fett{Q}-\fett{P}\fett{Q}$.
	Then, $L(\mathscr{H})$ is a quantum logic.
\end{beispiel}

Example \ref{bsp-projection-logic} motivates the following.
Call two projections $\fett{P}$, $\fett{Q}$ $\in$ $L(\mathscr{H})$ \emph{orthogonal},
in symbols $\fett{P}\bot \fett{Q}$, if $\fett{P} \leqq \fett{Q}'$. It follows that $\fett{P}\bot \fett{Q}$ if and only if
$\fett{P}\fett{Q}$ $=$ $\fett{Q}\fett{P}$ $=$ 0. If $\fett{P}$ and $\fett{Q}$ are orthogonal projections, then their
join is simply their sum, denoted traditionally $\fett{P} \oplus \fett{Q}$; in other words,
$\fett{P}\vee \fett{Q}$ $=$ $\fett{P} \oplus \fett{Q}$ if $\fett{P}\bot \fett{Q}$.
We denote the identity mapping on $\mathscr{H}$ by $1_{\mathscr{H}}$.

\begin{definition}
	A \emph{probability measure}\index{probability measure} on 
	$L$ $=$ $L(\mathscr{H})$ is a mapping $\mu$: $L\to [0,1]$ such that
	$\mu(1_{\mathscr{H}}) = 1$ and, for any sequence of pairwise orthogonal
	projections $\fett{P}_j$ $\in$ $L$, $j$ $=$ 1, 2, \ldots,
	\begin{equation}
		\textstyle
		\mu ( \oplus_j \fett{P}_j) = \sum_j \mu(\fett{P}_j).
	\end{equation}
\end{definition}
A way to construct a probability measure on $L(\mathscr{H})$ is,
for any unit vector $u\in \mathscr{H}$, to set 
$\mu_u (\fett{P}) = \langle \fett{P}u | u \rangle$.
This gives the orthodox Copenhagen interpretation\index{Copenhagen interpretation}
of quantum mechanics for the probability that $\fett{P}$ will yield the value 1 if the
physical system is in the state $u$.
Another way to express this fact is to write $\mu_u(\fett{P})$ $=$ Tr$(\fett{P}\fett{P}_u)$,
where $\fett{P}_u$ is the projection on the  one-dimensional subspace generated by 
the unit vector $u$.

More generally, probability measures $\mu_j$, $j$ $=$ 1, 2, \ldots, on
$L(\mathscr{H})$ form a {mixture}\index{mixture} 
$\mu$ $=$ $\sum_j t_j \mu_j$ where $0\leqq t_j \leqq 1$ and $\sum_j t_j = 1$
(the convex combination of the $\mu_j$'s).
Given any sequence $u_1$, $u_2$, \ldots of unit vectors in $\mathscr{H}$,
let be $\mu_j$ $=$ $\mu_{u_j}$ and $\fett{P}_j$ $=$ $\fett{P}_{u_j}$. Then for the
{density operator}\index{density operator} $\rho$ of the mixture,
\begin{equation}
	\rho = t_1 \fett{P}_1 + t_2 \fett{P}_2 + \ldots,
\end{equation}
we have
\begin{equation}
	\mu(\fett{P}) = t_1 \textrm{Tr}(\fett{P}\fett{P}_1) + t_2 \textrm{Tr}(\fett{P}\fett{P}_2) + \ldots
	       = \textrm{Tr}( \rho \fett{P}).
\end{equation}
Therefore, every density operator $\rho$ gives rise to a probability measure
$\mu$ on $L(\mathscr{H})$.
The remarkable theorem of Gleason 
shows the converse, i.e., that
to every probability measure there exists a density operator $\rho$.

\begin{satz}[Gleason (1957)]
	\label{satz-Gleason}
	Let $\mathscr{H}$ be a separable Hilbert space with dimension
	$\mathrm{dim}_{\mathbb{C}}\mathscr{H}$ $\geqq$ $3$. Then every probability measure on
	the space $L(\mathscr{H})$ of projections has the form
	$\mu(\fett{P})$ $=$ $\mathrm{tr}(\rho \fett{P})$ for a density operator $\rho$ on
	$\mathscr{H}$.
\end{satz}

An important direct consequence of Gleason's Theorem is that $L(\mathscr{H})$
does not permit any probability measures having only the values 0 and 1.
To see this, note that for any density operator $W$, the mapping 
$u \mapsto \langle Wu|u\rangle$ is continuous on the unit sphere of
$\mathscr{H}$. But since it is connected, no continuous function on it
can take only the discrete values 0 and 1.
Thus this result rules out the the possibility of ``hidden variables,''
an issue which had been subject of a long debate and was decided in the
experiments of Aspect in the early 1980's \cite{Aspect-et-al-1982}.

From the single premise that the ``experimental propositions'' associated with
a physical system are encoded by projections as in Example \ref{bsp-projection-logic}, 
one can reconstruct the formal apparatus of quantum mechanics. The first step is
Gleason's Theorem, which tells us that probability measures on the quantum logic
$L(\mathscr{H})$ correspond to density operators.
With the fundamental ``spectral theorem''
of functional analysis, stating that to every observable there exists a
certain family of projection operators, the observables are derivable 
by quantum logic. Even the dynamics of quantum mechanics, i.e., the unitary evolution,
can be deduced with the aid of a deep theorem of Wigner on the projective
representations of groups. For details see \cite{Varadarajan-1985}.

\subsection{Compatible propositions}
\begin{definition}
	\label{def-compatible}
	Two propositions $x$ and $y$ in a quantum logic $L$ are called 
	\emph{orthogonal}\index{orthogonal propositions}, symbolically
	$x \bot y$, if the relations
	\begin{equation}
		x \wedge y = 0, \qquad x \leqq y', \qquad y \leqq x',
		\label{eq-orthogonal-propositions}
	\end{equation}
	are satisfied. 
	Two propositions $x,$ $y \in  L$ are called 
	\emph{compatible}\index{compatible},
	or \emph{comeasurable}\index{comeasurable},
	if there exist $u$, $v$, $w$ $\in$ $L$ such that
	\begin{itemize}
		\item[(i)]
		$u$, $v$, $w$ are pairwise orthogonal;
		
		\item[(ii)]
		$u \vee v$ $=$ $x$ and $v \vee w = y$;
		
		\item[(iii)]
		The sublattice $B$ formed by $\{u, v, w, x, y, u', v', w', x', y'\}$ is orthomodular,
		i.e., each pair of propositions in $B$ satisfy (\ref{orthomodular-identity}).
	\end{itemize}
	We then call \{$u$, $v$, $w$\} a 
	\emph{compatible decomposition} 
	for $x$ and $y$.
\end{definition}

\begin{figure}[tb] 
\centering
	\begin{scriptsize}
	\unitlength1ex
	(a)\hspace*{-3ex}%
	\begin{picture}(14,13)
		\put(0.2,0){\includegraphics[height=14ex]{./F2}}
		\put( 8.9,12.5){\makebox(0,0)[l]{$= x \vee y$}}
		\put( 7.2,12.7){\makebox(0,0){$1$}}
		\put( 1.6, 7.0){\makebox(0,0){$x$}}
		\put(12.8, 6.8){\makebox(0,0){$y$}}
		\put( 7.2, 1.5){\makebox(0,0){$0$}}
	\end{picture}
	\qquad \quad
	(b)\hspace*{-3ex}%
	\begin{picture}(14,20)
		\put(0,0){\includegraphics[height=20ex]{./Z_2-3}}
		\put( 8.8,18.5){\makebox(0,0)[l]{$= x\vee y$}}
		\put( 7.2,18.6){\makebox(0,0){$1$}}
		\put( 1.5,12.8){\makebox(0,0){$x$}}
		\put( 7.2,13.0){\makebox(0,0){$v'$}}
		\put(12.8,12.6){\makebox(0,0){$y$}}
		\put( 1.5, 7.0){\makebox(0,0){$u$}}
		\put( 7.2, 7.0){\makebox(0,0){$v$}}
		\put(13.0, 7.0){\makebox(0,0){$w$}}
		\put( 7.2, 1.5){\makebox(0,0){$0$}}
	\end{picture}
	\qquad \quad
	(c)\hspace*{-3ex}%
	\begin{picture}(30,25.7)
		\put(0,0){\includegraphics[height=25.7ex]{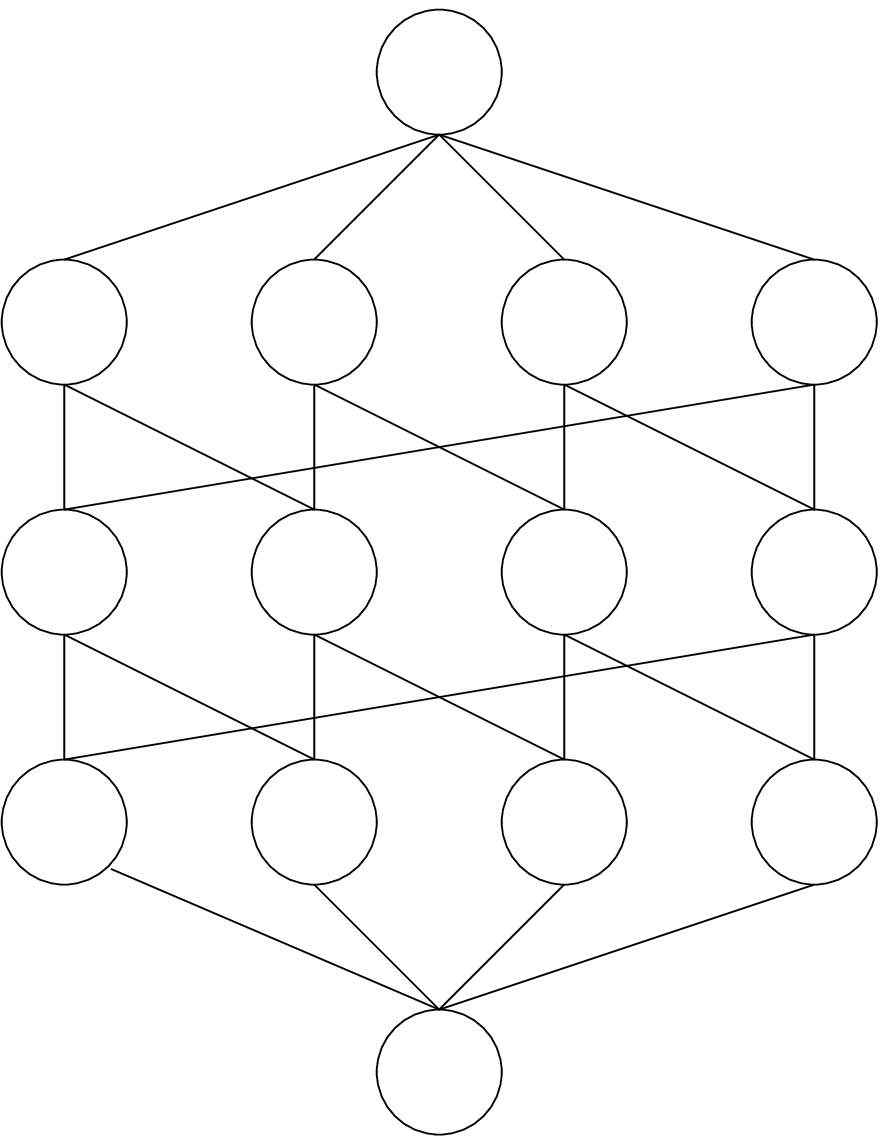}}
		\put(10.0,24.3){\makebox(0,0){$1$}}
		\put( 1.4,20.5){\makebox(0,0)[b]{$x \vee y$}}
		\put( 7.2,18.6){\makebox(0,0){$u'$}}
		\put(12.9,18.6){\makebox(0,0){$v'$}}
		\put(18.5,18.6){\makebox(0,0){$w'$}}
		\put( 1.3,12.8){\makebox(0,0){$x$}}
		\put( 7.2,12.6){\makebox(0,0){$y$}}
		\put(12.9,13.0){\makebox(0,0){$x'$}}
		\put(18.7,12.8){\makebox(0,0){$y'$}}
		\put( 1.4, 7.2){\makebox(0,0){$u$}}
		\put( 7.2, 7.2){\makebox(0,0){$v$}}
		\put(12.9, 7.2){\makebox(0,0){$w$}}
		\put(18.6, 5.1){\makebox(0,0)[t]{$x'\wedge y'$}}
		\put(10.0, 1.5){\makebox(0,0){$0$}}
	\end{picture}
	\end{scriptsize}
\caption{\label{fig-compatible} 
	Compatible decompositions $\{u, v, w\}$ for $x$ and $y$:
	(a) $F_2$ $=$ 
	$\mathbf{2}^2$ for $x \bot y$, i.e., $v = 0$, $w=x$, $u=y$;
	(b) $\mathbf{2}^3$ for $v \ne 0$, $w'=x$, $u'=y$;
	(c) for $v\ne 0$, $w' \ne x$, $u' \ne y$.
}
\end{figure}

By (\ref{disjunctive-De-Morgan}) two propositions $x$, $y$ with the
compatible decomposition $\{u, v, w\}$ satisfy
\begin{equation}
	x' = u' \wedge v' \geqq w,
	\qquad
	y' = v' \wedge w' \geqq u,
	\qquad
	x \vee y = u \vee v \vee w.
\end{equation}
Clearly, two orthogonal propositions\index{orthogonal proposition}
$x$, $y \in L$ of a logic 
are always compatible, since
we may simply identify $u=x$, $v=0$, and $w=y$.
The notion ``orthogonal'' originally refers to the geometrical relationship 
between two propositions of a quantum logic, which are subspaces of a Hilbert space.
Figure \ref{fig-compatible} depicts three compatible decompositions.
They are distributive, (a) and (b) are Boolean and (c) is a sublattice of
the Boolean lattice $\mathbf{2}^4$.

\begin{satz}
	\label{theo-u-v-w}
	Let $L$ be a quantum logic. 
	If $\{u$, $v$, $w\}$ is a compatible decomposition for $x$ and $y$
	in $L$, then
	\begin{equation}
		\label{compatible-1}
		u = x \wedge y',
	\end{equation}
	\begin{equation}
		\label{compatible-2}
		w = x' \wedge y,
	\end{equation}
	\begin{equation}
		\label{compatible-3}
		v = x \wedge y
		= (x' \vee y) \wedge x = (x \vee y') \wedge y.
	\end{equation}
	Therefore, $u$, $v$, and $w$ are uniquely determined by $x$ and $y$.
\end{satz}
\begin{proof}
	By Theorem \ref{theo-quantum-logic} and Theorem \ref{theo-De-Morgan-fuzzy},
	in a quantum logic De Morgan's laws
	(\ref{disjunctive-De-Morgan}) and (\ref{conjunctive-De-Morgan}) hold.

	Proof of (\ref{compatible-1}): 
	Since by definition we have $u\leqq w'$ and $u \leqq v'$, we have
	$u \leqq v' \wedge w' = (v \vee w)' = y'$, where the first equality
	follows from (\ref{conjunctive-De-Morgan}). But also $u \leqq x$ (by definition), so
	\begin{equation}
	\label{compatible-*}
		u \leqq x \wedge y'.
	\end{equation}
	Now $x \wedge y' \leqq y'$, thus $u \leqq y'$, implying that $y\leqq u'$,
	so by the orthomodular identity (\ref{orthomodular-identity}),
	$u' = y \vee (y'\wedge u')$, or by (\ref{conjunctive-De-Morgan})
	$u = y' \wedge (y \vee u)$. Now $y\vee u = y \vee x$. Now if $r \leqq x$
	and $r \leqq y'$, then $r \leqq y'\wedge x \leqq y' \wedge (y \vee x) = u$.
	This together with (\ref{compatible-*}) establishes that 
	$u = \inf\{x,y'\} = x \wedge y'.$
	
	Proof of (\ref{compatible-2}): 
	This is proved with an argument analogous to (\ref{compatible-1}).

	Proof of (\ref{compatible-3}): 
	First observe that $w$, $u$ $\leqq$ $v'$, so $w \vee u \leqq v'$.
	By (\ref{orthomodular-identity}) thus 
	$v' = w \vee u \vee ((w\vee u)' \wedge v')$, so that
	$v$ $=$ $(w \vee u)'$ $\wedge$ $(w \vee u \vee v)$ 
	$=$ $((x' \wedge y)$ $\vee$ $(x \wedge y'))$ $\wedge$ $(x\vee y)$
	$=$ $((x \vee y')$ $\wedge$ $(x' \vee y))$ $\wedge$ $(x \vee y)$.
	Now $v$ $\leqq$ $x$, $y$, and we shall immediately establish that
	$v=\inf\{x,y\} = x \wedge y$. Let $r$ $\leqq$ $x$, $y$. Then
	$r$ $\leqq$ $((x \vee y')$ $\wedge$ $(x' \vee y))$ $\wedge$ $(x \vee y)$ $=$ $v$.
	This establishes the first equality in (\ref{compatible-3}).
	For the second equality, observe that $u \leqq v'$ so that by the
	orthomodular identity (\ref{orthomodular-identity}),
	$v'$ $=$ $u \vee (u' \wedge v')$, or
	$v$ $=$ $u' \wedge (u \vee v)$ $=$ $(x' \vee y) \wedge x$.
	The third equality is proved similarly.
\end{proof}

\begin{satz}
	\label{theo-compatible}
	Let $L$ be a quantum logic. 
	Then two propositions $x$ and $y$ $\in$ $L$ are compatible if and only if
	\begin{equation}
		\label{theo-compatible-1}
		x = (x \wedge y) \vee (x \wedge y')
	\end{equation}
	and
	\begin{equation}
		\label{theo-compatible-2}
		y = (y \wedge x) \vee (y \wedge x').
	\end{equation}
\end{satz}
\begin{proof}
	Since $x \wedge y \leqq x$ and $x \wedge y' \leqq x$, we deduce
	$(x \wedge y)$ $\vee$ $(x \wedge y')$ $\leqq$ $x$, whether or not 
	$x$ and $y$ are compatible.
	The reverse inequality follows from Theorem \ref{theo-u-v-w};
	for if \{$u$, $v$, $w$\} is a compatibility decomposition for
	$x$ and $y$, then $x \wedge y = v$, $x \wedge y' = u$,
	and $u \vee v = x$. This yields (\ref{theo-compatible-1}).
	Similarly, (\ref{theo-compatible-2}) follows by an analogous argument.
	Conversely, if (\ref{theo-compatible-1}) and (\ref{theo-compatible-2})
	hold, we have that \{$x \wedge y'$, $x \wedge y$, $x' \wedge y'$\}
	form a compatibility decomposition of $x$ and $y$.
\end{proof}

In \cite{Mittelstaedt-1978} propositions satisfying 
(\ref{theo-compatible-1}) and (\ref{theo-compatible-2})
are called ``commensurable''\index{commensurable},
in \cite{Cohen-1989} ``commuting.''
These conditions reveal the intuitive meaning of the 
orthomodularity condition: \emph{comparable elements are compatible}. 
Notice that by the first distributive inequality (\ref{5}), with $z=y'$,
we have always $x \geqq (x \wedge y) \vee (x \wedge y')$. Theorem
\ref{theo-compatible} then states that we generally do not have the reverse
$x \leqq (x \wedge y) \vee (x \wedge y')$ unless $x$ and $y$ are compatible.
Thus, if $x$ and $y$ are not compatible,
then knowing that proposition $x$ is true \emph{is not sufficient}
for concluding that at least one of the following ist true:
\begin{itemize}
	\item[(i)]
	$x$ and $y$ are simultaneously true (i.e., $u$ is true), or
	
	\item[(ii)]
	$x$ and ``not-$y$'' are simultaneously true (i.e., $v$ is true).
\end{itemize}
This logical structure captures the idea of Heisenberg's uncertainty principle
that for some physical systems there might exist a pair of propositions
whose truth values simply cannot be simultaneously determined.
The algebraic structure which guarantees compatibility is the
distributivity.

\subsection{The logic of quantum registers}

\subsubsection{The Boolean logic of a single qubit}
\begin{figure}[htb] 
\centering
	\begin{footnotesize}
	\unitlength1ex
	\begin{picture}(18,18)
		\put(0,0){\includegraphics[width=18ex]{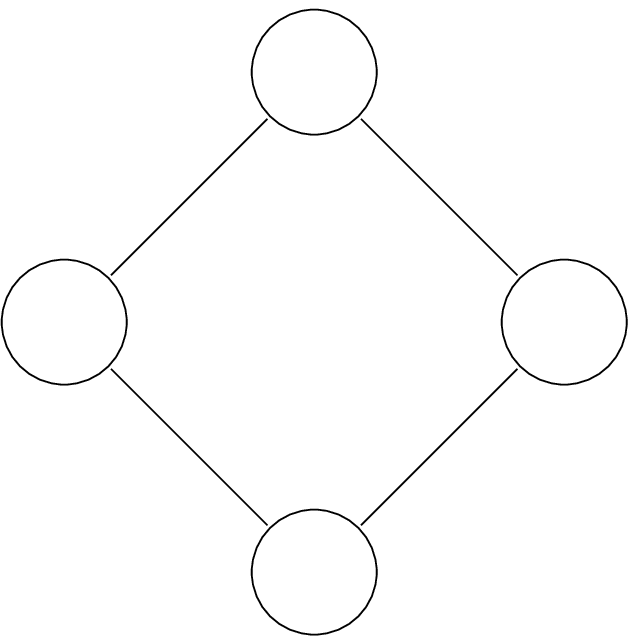}}
		\put( 9.0,17.0){\makebox(0,0)[t]{$1$}}
		\put( 1.9, 9.2){\makebox(0,0){$p^-$}}
		\put(16.3, 9.2){\makebox(0,0){$p^+$}}
		\put( 9.0, 1.2){\makebox(0,0)[b]{$0$}}
	\end{picture}
	\end{footnotesize}
	\caption{\label{fig-MO-1} 
	The Boolean lattice $MO_1$ $\simeq$ $\mathbf{2}^2$ of a single qubit,
	with the (orthogonal) states $p^-$ and $p^+$.
	(Cf.\ Figs.\ \ref{fig-F2}) 
	}
\end{figure}
A two-state quantum system, such as a spin-$\frac12$ particle, is described
by the two-di\-men\-si\-o\-nal Hilbert space $\mathscr{H} = \mathbb{C}^2$. 
Consider measurements of the spin-component along a particular direction, say
along the $z$-axis (any direction will do as well). This can be
operationalized by a Stern-Gerlach type experiment 
\cite{Goswami-1997}
using an inhomogeneous magnetic field.
There are two possible spin components of the particle, namely
spin $-\frac12 \hbar$ and $+\frac12 \hbar$; we will shortly say
that the particle is in state ``$+$'' or $|0\rangle$ 
if it has spin $+\frac12 \hbar$,
and in state ``$-$'' or $|1\rangle$ if it has spin $-\frac12 \hbar$.
This corresponds to the following propositions.
\begin{itemize}
	\item[$p^-$:]
	``The particle is in state `$-$'\,''
	$=$ 1-dimensional subspace $\mathrm{span}\,\{|1\rangle\}$

	\item[$p^+$:]
	``The particle is in state `$+$'\,''
	$=$ 1-dimensional subspace $\mathrm{span}\,\{|0\rangle\}$

	\item[1:]
	``The particle is in some state''
	$=$ whole space $\mathscr{H} = \mathbb{C}^2$

	\item[0:]
	``The particle is in no state''
	$=$ zero-dimensional subspace $\{(0,0)\}$
\end{itemize}
The proposition 1 is the tautology\index{tautology}, 0 is the
absurd statement. The propositions $p^-$ and $p^+$ are the
atoms. Since they are complements of each other, i.e.,
$p^+ = (p^-)'$ and $p^- = (p^+)'$, they form a Boolean logic
(Figure \ref{fig-MO-1}).

\subsubsection{Quantum register of size $n$}
\begin{figure}[tb] 
\centering
\begin{footnotesize}
	\unitlength1ex
	\begin{picture}(54,18)
		\put(0,0){\includegraphics[width=54ex]{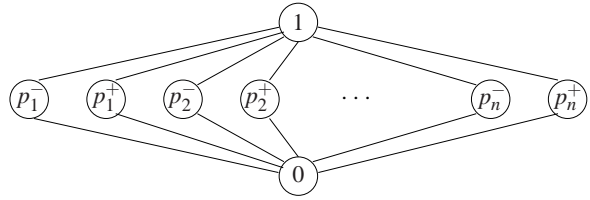}}
		\put(27.0,17.0){\makebox(0,0)[t]{$1$}}
		\put( 2.0, 9.2){\makebox(0,0){$p_1^-$}}
		\put( 9.1, 9.2){\makebox(0,0){$p_1^+$}}
		\put(16.2, 9.2){\makebox(0,0){$p_2^-$}}
		\put(23.3, 9.2){\makebox(0,0){$p_2^+$}}
		\put(45.1, 9.2){\makebox(0,0){$p_n^-$}}
		\put(52.2, 9.2){\makebox(0,0){$p_n^+$}}
		\put(27.0, 1.2){\makebox(0,0)[b]{$0$}}
	\end{picture}
\end{footnotesize}
\caption{\label{fig-MO-n} 
	The modular lattice $MO_n$ of a quantum register of $n$ qubits,
	if only single-qubit measurements are considered.
}
\end{figure}

A quantum register of size $n$ is a physical system of $n$ qubits.
Mathematically, it is a composition of several single systems 
(a `tensor product'). A general lattice representing 
such a composition seems to be impossible \cite[p.~51]{Svozil-1998},
but for the special case of a 
a composite system without entanglement, for instance a 
measurement of a single qubit of a quantum register,
a lattice is given as follows.
According to the ``pasting construction'' \cite[\S3.2]{Svozil-1998}
each qubit is considered as a Boolean ``block,'' and identical
propositions in different blocks are identified such that the
logical structure in each block remains intact.
This yields the lattice $MO_n$ (for ``modular orthocomplemented'')
in Fig.\ \ref{fig-MO-n}, consisting of $2n$ atoms $p_n^\pm$
satisfying
\begin{equation}
	p_j^+ = (p_j^-)',
	\qquad
	p_j^- = (p_j^+)'
	\quad
	\mbox{for } j = 1, 2, \ldots, n.
\end{equation}
Formally, $MO_n = \bigoplus_{j=1}^n L_j$ with
$L_j = \{0,p_j^+,p_j^-,1\}$ where 0 and 1 in each $L_j$ are identified.
Therefore, $p_j^+ \wedge p_j^- = 0$, and $p_j^+ \vee p_j^- = 1$. 
For $n>1$, the lattices $MO_n$ are not distributive, since for
$j \ne k$, with $j$, $k$ $=$ 1, \ldots, $n$,  we have
$
	p_j^+ \vee (p_k^+ \wedge p_k^-)
	= p_j^+ \vee 0
	= p_j^+
$, 
but
$
	(p_j^+ \vee p_k^+) \wedge (p_j^+ \vee p_k^-)
	= 1 \wedge 1
	= 1.
$

\subsubsection{Entire register measurements without entanglement}
In case of a measurement of the entire quantum register of size $n$,
the atoms are given by the four propositions
$p_{00}$, $p_{01}$, $p_{10}$, and $p_{11}$, where $p_{ij}$ 
corresponds to the basis vector
$|ij\rangle$ of the four-dimensional Hilbert space 
$\mathscr{H} = \mathbb{C}^4$.
They are composed to the following propositions:
\begin{equation}
\begin{array}{r@{\ =\ }l@{\mbox{, }\quad}r@{\ =\ }l@{\mbox{, }\quad}r@{\ =\ }l@{\mbox{,}}}
	u & p_{00} \vee p_{01} &
	v & p_{00} \vee p_{10} & 
	w & p_{00} \vee p_{11} \\
	x & p_{01} \vee p_{10} &
	y & p_{01} \vee p_{11} &
	z & p_{10} \vee p_{11}
\end{array}
\end{equation}
\begin{figure}[tb] 
\centering
\begin{footnotesize}
	\unitlength1ex
	\begin{picture}(40,32)
		\put(0,-.5){\includegraphics[width=40ex]{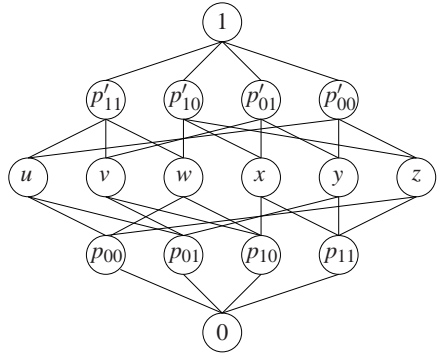}}
		\put(20.0,30.5){\makebox(0,0){$1$}}
		\put( 9.2,23.4){\makebox(0,0){$p_{11}'$}}
		\put(16.5,23.4){\makebox(0,0){$p_{10}'$}}
		\put(23.7,23.4){\makebox(0,0){$p_{01}'$}}
		\put(31.0,23.4){\makebox(0,0){$p_{00}'$}}
		\put( 1.7,16.0){\makebox(0,0){$u$}}
		\put( 9.0,16.0){\makebox(0,0){$v$}}
		\put(16.4,16.0){\makebox(0,0){$w$}}
		\put(23.5,16.0){\makebox(0,0){$x$}}
		\put(30.8,15.8){\makebox(0,0){$y$}}
		\put(38.1,16.0){\makebox(0,0){$z$}}
		\put( 9.2, 8.6){\makebox(0,0){$p_{00}$}}
		\put(16.5, 8.6){\makebox(0,0){$p_{01}$}}
		\put(23.7, 8.6){\makebox(0,0){$p_{10}$}}
		\put(31.0, 8.6){\makebox(0,0){$p_{11}$}}
		\put(20.0, 1.4){\makebox(0,0){$0$}}
	\end{picture}
\end{footnotesize}
	\caption{\label{fig-2-4} 
	The Boolean lattice $\mathbf{2}^4$ of a quantum register of 2 qubits,
	if entire register measurements are considered.
	}
\end{figure}
This yields the Boolean lattice $\mathbf{2}^4$, cf.\ Fig.\ \ref{fig-2-4}
and \cite[§5.2]{Svozil-1998}.

\section{Paraconsistent logics and effect structures}
The (orthodox) quantum logic based on the proposals of Birhhoff and von Neumann
as described above, is both a ``total'' and a ``sharp'' logic.
A logic is \emph{total}\index{total} if the set of ``meaningful propositions'' is
closed under the basic logical operations, i.e., the conjunction as well as the 
disjunction of two meaningful propositions is again a meaningful proposition.
A logic is \emph{sharp}\index{sharp} if each proposition corresponds to
exact physical properties of the corresponding physical system.

There has been deep criticism in different contexts on both the total and the sharp 
character of quantum logic.
One main objection to orthodox quantum logic, for instance, has been that by the 
identification of propositions with the projections of the Hilbert space given
the quantum system,
propositions in fact are identified with the \emph{physical properties} of the system.
This one-to-one correspondence implies the identification of the extensional notion
of a proposition (according to the standard tradition of semantics) with the
collapse of the empirical and intensional concepts of ``experimental proposition,''
``physical property.''
Although quite convenient mathematically, 
this collapse has been called a ``metaphysical disaster'' 
\cite{Randall-Foulis-1983}.
It stimulated the investigation about more and more general quantum 
structures which, however, do not yield lattice structures in a direct way.
The main goal is to find some algebraic structure on the semantic level
and to derive a quantum logic.

An important extension of the mathematical representation is the notion of an ``effect''
as a representative of an experimental proposition
\cite{Dalla-Chiara-et-al-2004, Dvurecenskij-Pulmannova-2000}.
The main idea is to define a quantum physical observable operationally in terms of an
experimental procedure or a class of them,
and the ensuing measurement statistics are to be described by probability measures
depending on the input states.
In this way, any observable is a normalized 
positive-operator-valued measure (POVM)\index{POVM} $A$,
assigning to each measurement outcome $x$ or outcome range $X$ 
(Borel set, usually a real set) its
measurement probability $A(\{x\})$ or $A(X)$.
If the measurement outcomes are finite, e.g., \{$x_1$, \ldots, $x_n$\},
it is custom to write the POVM as the set \{$A_1$, \ldots, $A_n$\}
such that $\langle \psi| A_j |\psi\rangle$ yields the probability to
measure the state $|\psi\rangle$.
Necessarily, $0 \leqq \langle \psi| A_j |\psi\rangle \leqq 1$,
and $\sum_j A_j = I$.
In terms of density operator, we have $\rho = |\psi\rangle \langle \psi|$,
and $\mathrm{tr}\,(A \rho) = \langle \psi| A |\psi\rangle$.
\cite[§2.2.6]{Nielsen-Chuang-2000}
Any such operator $A$ is an effect.

\begin{definition}
	Let $\mathscr{H}$ be a given Hilbert space describing a quantum system.
	Then the set $E(\mathscr{H})$ of all \emph{effects}\index{effect} 
	of $\mathscr{H}$ is defined as the set of all linear bounded 
	self-adjoint operators
	$A$ on $\mathscr{H}$ such that
	$
		\mathrm{tr}\,(A \rho) \in [0,1]
	$ 
	for an arbitrary density operator $\rho$.
\end{definition}

An effect $A$ is a projection if and only if $A^2 = A$, that is, projections
are exactly the idempotent effects.
An important difference between effects and projections is that effects may
represent fuzzy propositions like ``the value for the observable $A$ lies in the
fuzzy Borel set $\mu_A$.''
In particular, there exists effects $A$ different from the null projection
$O$ 
such that no state $\rho$ can verify $A$ with probability 1, i.e.,
$\mathrm{tr}\,(A \rho) < 1$ for all statistical operators $\rho$.
A limiting case is the \emph{semitransparent effect}\index{semitransparent effect}
$\frac12 I$, to which any statistical operator $\rho$ assigns the
probability $\frac12$. It represents the prototypical ambiguity.

The set $E(\mathscr{H})$ of effects can be naturally structured as a 
poset
$(E(\mathscr{H}), \sqsubseteq)$
with minimal bound $O$ and maximal bound $I$
where $A \sqsubseteq B$ means by definition that 
$\mathrm{tr}\,(A \rho) \leqq \mathrm{tr}\,(B \rho)$
for any statistical operator $\rho$.
It is easily checked that the null projection $O$ and the identity
$I$ are minimum and maximum, respectively, with respect to $\sqsubseteq$.
However,
$(E(\mathscr{H}), \sqsubseteq)$
is not a lattice.
Contrary to projections, there exist pairs of effects which have no uniquely
defined infimum or supremum.

\begin{beispiel}
	Consider the following effects on the Hilbert space $\mathscr{H} = \mathbb{C}^2$,
	given in matrix-representation:
	\[ 
	A = \left( \begin{array}{cc}
		\frac12 &       0 \\[1.0ex]
		0       & \frac12
	\end{array} \right),
	\quad
	B = \left( \begin{array}{cc}
		\frac34 &       0 \\[1.0ex]
		0       & \frac14
	\end{array} \right),
	\quad
	C = \left( \begin{array}{cc}
		\frac12 &       0 \\[1.0ex]
		0       & \frac14
	\end{array} \right),
	\quad
	D = \left( \begin{array}{cc}
		\frac7{16} & \frac18 \\[1.0ex]
		\frac18    & \frac3{16}
	\end{array} \right).
	\] 
	Then $C \sqsubseteq A$, $B$, as well as $D \sqsubseteq A$, $B$.
	However, $C \not\sqsubseteq D$ and $D \not\sqsubseteq C$,
	i.e., there does not exist an infimum of $A$ and $B$.
\end{beispiel}

A possibility to obtain a lattice structure from $(E(\mathscr{H}), \sqsubseteq)$
is to embed it into its Mac Neille completion.

\begin{definition}
	Let $(\mathscr{B}, \sqsubseteq)$ 
	be a poset with universal bounds $O$ and $I$
	and an involution $'$ satisfying
	(i) $A=A''$ and (ii) $A \sqsubseteq B$ $\Rightarrow$ $B' \sqsubseteq A'$.
	Then the \emph{Mac Neille completion}\index{Mac Neille completion}
	is defined as the tuple 
	$(\mathrm{MC}\,(\mathscr{B}), \subseteq, ', \{0\}, \mathscr{B})$
	where
	\begin{equation}
		\mathrm{MC}(\mathscr{B})=\{X \subseteq \mathscr{B}: \ X = u(l(X))\}
	\end{equation}
	with $l(X)$ and $u(X)$ denoting all lower bounds and all upper bounds of $X$, 
	respectively, and
	\begin{equation}
		X' := 
		\{A \in \mathscr{B}: \
		A \sqsubseteq B' \ \forall B \in X\}
		.
	\end{equation}	
\end{definition}

It turns out that $X \in \mathrm{MC}\,(\mathscr{B})$ if and only if
$X = X''$. Moreover, $\mathrm{MC}(\mathscr{B})$
is universally bounded lattice where
$X \sqcap Y = X \cap Y$ and $X \sqcup Y = (X \cup Y)''$.
The negation is fuzzy and does not satisfy neither the non-contradiction law
nor the \emph{tertium non datur}
\cite{Dalla-Chiara-Giuntini-2002}.

Another common possibility to derive a logic structure from effects is
provided by the fact that effects naturally form an algebra
(``effect algebra'')
\cite{Dvurecenskij-Pulmannova-2000} which is equivalent to
a quasilinear QMV algebra \cite{Dalla-Chiara-Giuntini-2002}.
In this approach, the operator
\begin{equation}
	A \oplus B
	= \left\{ \begin{array}{ll}
		A+B & \mbox{if $A+B \in E(\mathscr{H})$,} \\
		I & \mbox{otherwise,}
	\end{array}\right.
\end{equation}
and the negation
\begin{equation}
	A' = I - A
\end{equation}
are introduced.
The semitransparent effect $\frac12 I$ is a fixed point of the negation,
i.e., it is its own negation.
Moreover, $\frac12 I \oplus \frac12 I = I$.

\subsection{Paraconsistent quantum logics (PQL)}
\begin{beispiel}
	\label{bsp-pql}%
	\emph{(The lattices ${G}_{6}$, ${G}_{8}$, ${G}_{14}$)}
	Among the simplest finite paraconsistent logics are the
	the lattices 
	${G}_{6}$\index{$\mathscr{G}_{8}$},
	${G}_{8}$\index{$\mathscr{G}_{8}$},
	and
	${G}_{14}$\index{$\mathscr{G}_{8}$}
	depicted in Figure \ref{fig-G-8}.
\begin{figure}[htb] 
\centering
	\begin{footnotesize}
	\unitlength1ex
	\begin{picture}(16,28)
		\put(0,0){\includegraphics[width=16ex]{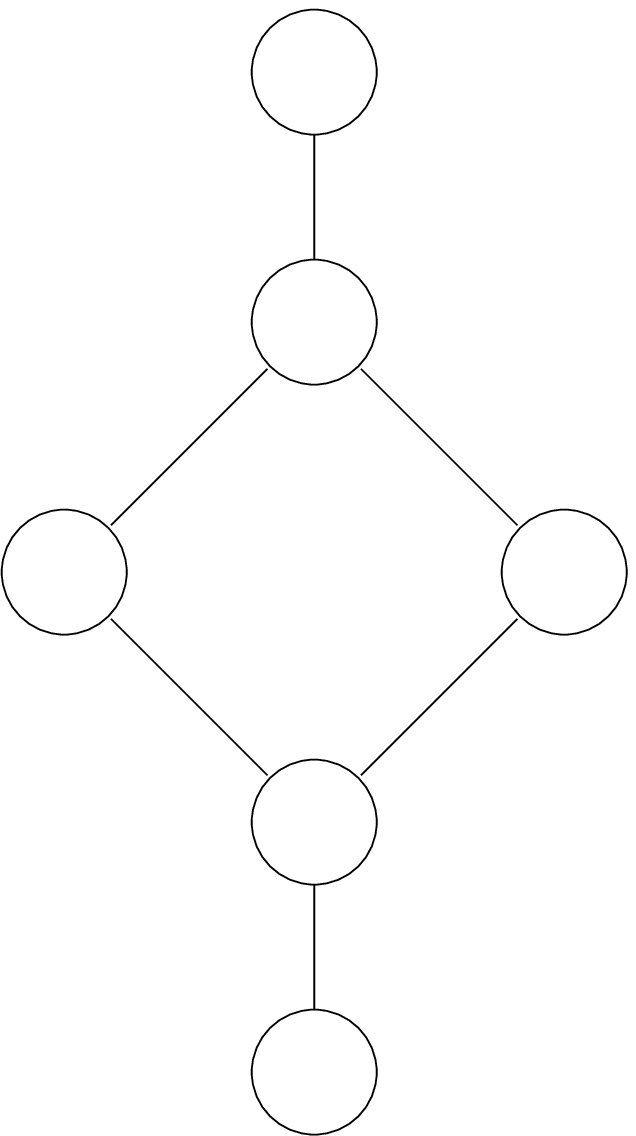}}
		\put( 8.0,27.2){\makebox(0,0){$1$}}
		\put( 8.0,21.0){\makebox(0,0){$x'$}}
		\put( 1.8,14.3){\makebox(0,0){$y$}}
		\put(14.4,14.5){\makebox(0,0){$y'$}}
		\put( 8.0, 7.8){\makebox(0,0){$x$}}
		\put( 8.0, 1.0){\makebox(0,0)[b]{$0$}}
	\end{picture}
	\qquad
	\begin{picture}(15.5,34)
		\put(0,0){\includegraphics[width=15.5ex]{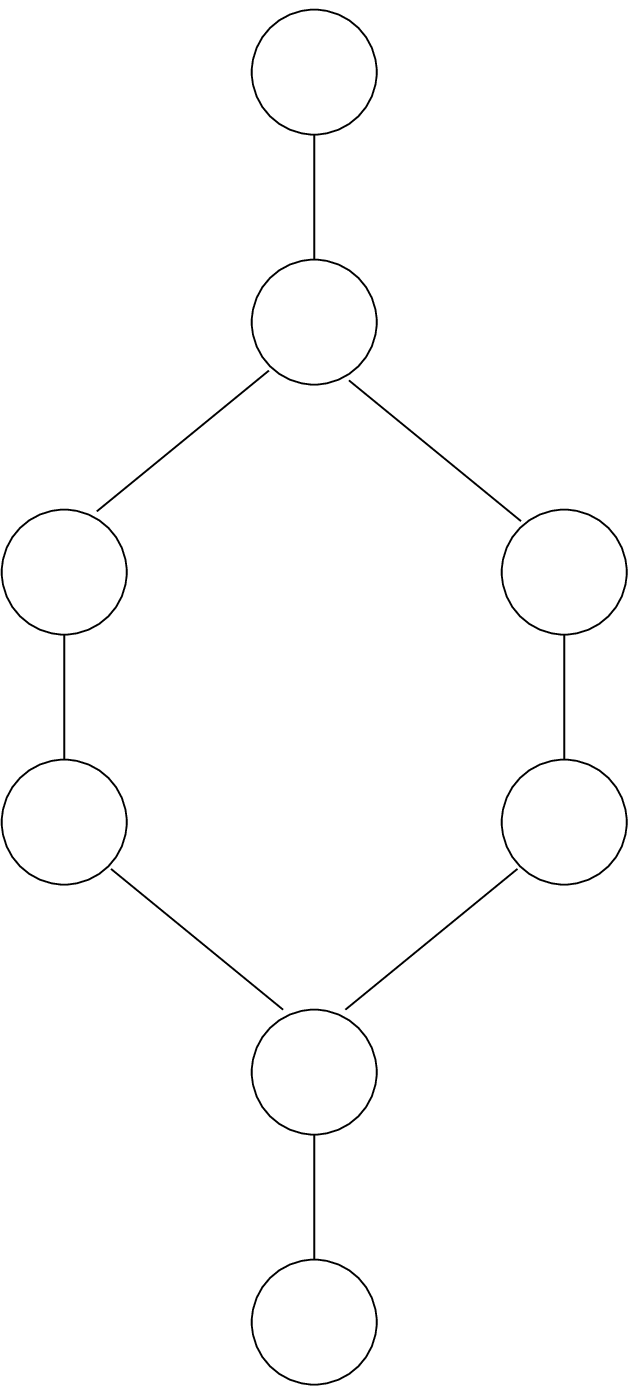}}
		\put( 7.7,32.7){\makebox(0,0){$1$}}
		\put( 7.7,26.3){\makebox(0,0){$f'$}}
		\put( 1.8,20.0){\makebox(0,0){$y$}}
		\put(14.0,20.3){\makebox(0,0){$x'$}}
		\put( 1.6,14.0){\makebox(0,0){$x$}}
		\put(14.0,14.0){\makebox(0,0){$y'$}}
		\put( 7.7, 7.8){\makebox(0,0){$f$}}
		\put( 7.7, 1.0){\makebox(0,0)[b]{$0$}}
	\end{picture}
	\qquad
	\begin{picture}(28,34)
		\put(0,0){\includegraphics[width=28ex]{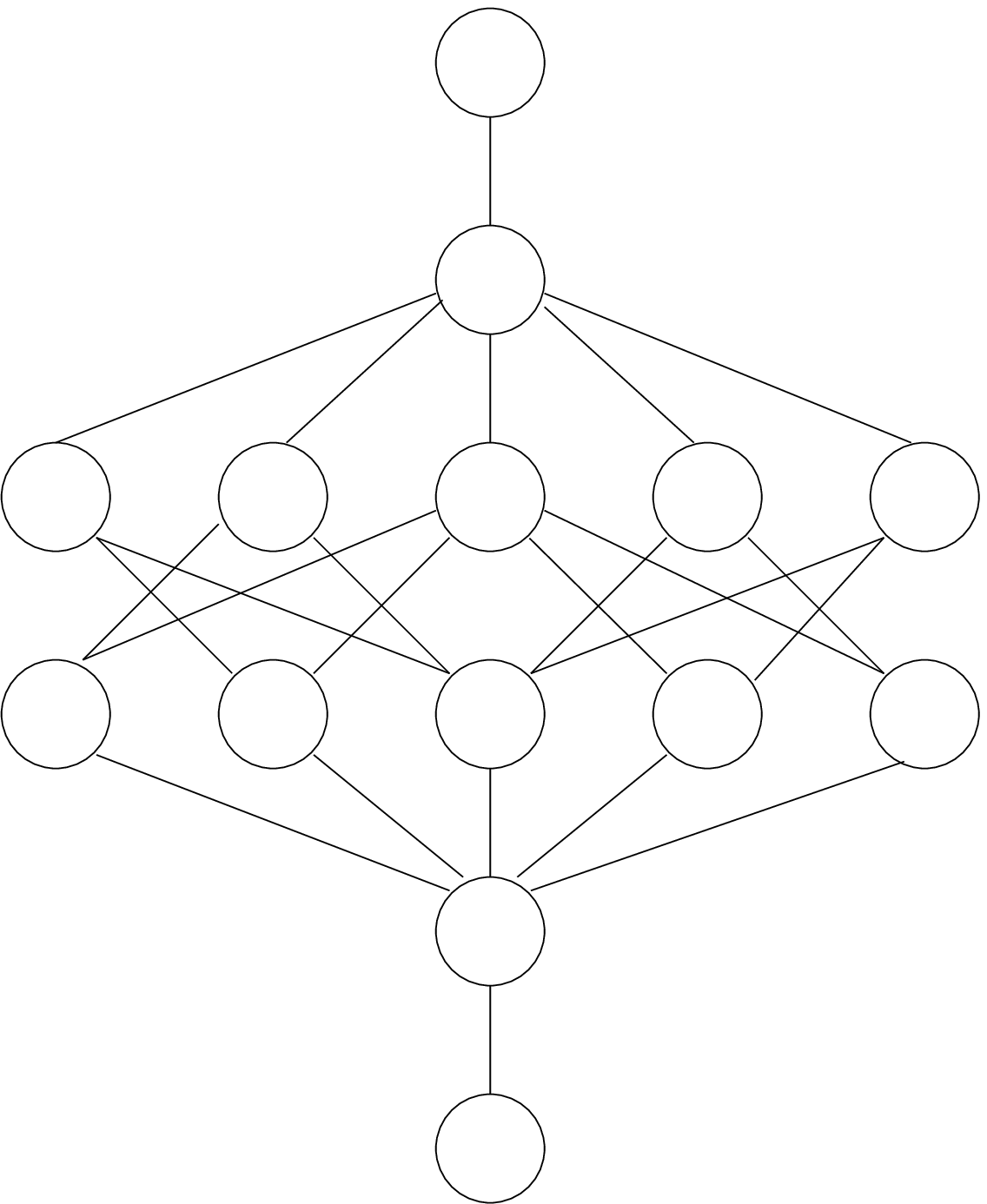}}
		\put(14.0,32.7){\makebox(0,0){$1$}}
		\put(14.0,26.5){\makebox(0,0){$f'$}}
		\put( 1.6,20.3){\makebox(0,0){$a'$}}
		\put( 7.8,20.3){\makebox(0,0){$b'$}}
		\put(14.0,20.3){\makebox(0,0){$c'$}}
		\put(20.2,20.3){\makebox(0,0){$d'$}}
		\put(26.4,20.3){\makebox(0,0){$e'$}}
		\put( 1.6,14.0){\makebox(0,0){$a$}}
		\put( 7.8,14.2){\makebox(0,0){$b$}}
		\put(14.0,14.0){\makebox(0,0){$c$}}
		\put(20.2,14.2){\makebox(0,0){$d$}}
		\put(26.4,14.0){\makebox(0,0){$e$}}
		\put(14.0, 7.8){\makebox(0,0){$f$}}
		\put(14.0, 0.8){\makebox(0,0)[b]{$0$}}
	\end{picture}
	\end{footnotesize}
	\caption{\label{fig-G-8} 
		The paraconsistent quantum logics 
		$G_6$, $G_{8}$, and $G_{14}$.
	}
\end{figure}
	They are not orthomodular, 
	in $G_{6}$ and $G_{8}$,
	for instance,
	$
		x \vee (x' \wedge y) = x \ne y 
	$
	although $x < y$. However, they satisfy (\ref{paraconsistency})
	since the premise of (\ref{paraconsistency}) is satisfied 
	only if both $x$ and $y$
	are either 0 or 1.
	The lattice $G_{14}$ is given in \cite{Dalla-Chiara-Giuntini-2002}.
\end{beispiel}

\begin{satz}
	Orthomodularity implies paraconsistency, but not vice versa.
\end{satz}
\begin{proof}
	By orthomodularity (\ref{orthomodular-identity}) we have
	$y = x \vee (x' \wedge y)$ for $x \leqq y$, implying $x=y$ if additionally
	$x' \wedge y = 0$. 
	A paraconsistent logic which is not a quantum logic, on the other hand,
	is $G_{6}$, cf.\ Example \ref{bsp-pql}.
\end{proof}

In other words: A quantum logic is always paraconsistent,
but there exist paraconsistent logics which are not quantum logics.

\begin{satz}
	A non-contradictory logic is not necessarily paraconsistent.
\end{satz}
\begin{proof}
	A non-contradictory logic which is not paraconsistent is
	$O_6$ with $x \vee y = 1$ in Figure \ref{fig-O6}, since
	$x \ne y$ although $x \leqq y$ and $x' \wedge y = 0$. 
\end{proof}

The Mac Neille completion of the effect poset 
$E(\mathscr{H})$
is a paraconsistent logic.
Paraconsistent quantum logic is the most obvious fuzzy weakening of quantum logic
\cite{Dalla-Chiara-Giuntini-2002}.
From the physical point of view, propositions in PQL (possible worlds)
represent pieces of information about the physical system under investigation.
Any information may be either maximal (a pure state) or non-maximal (a mixed state),
either sharp (a projection) or fuzzy (a proper effect).
Violations of the non-con\-tra\-dic\-tion law are caused by fuzzy (ambiguous)
pieces of knowledge.

Another important property of the paraconsistent quantum logic is that it is
a sublogic of a wide class of logics, such as Girard's linear logic,
{\L}ukasiewicz infinite many-valued logic $L_{\mathbb{R}}$ and some relevant logics
\cite{Dalla-Chiara-Giuntini-2002}.



\section{Implicative lattices and intuitionistic logic}
\label{sec-intuitionistic-logic}
Roughly speaking, ``intuitionism'' holds that logic and mathematics are 
``constructive'' mental activities, i.e.,
theorems are not \emph{discovered}, but \emph{invented}.
Thus, logic and mathematics are the application of internally 
consistent methods to realize more complex mental constructs.
Intuitionistic logic\index{intuitionistic logic}\index{logic! intuitionistic -}
is the logic used in mathematical constructivism,
introduced by
Brouwer %
and Heyting. 

In classical logic, a formula, say $P$, asserts that $P$ is true.  
In intuitionistic logic a formula is only considered to be true if it can be 
\emph{proved}.
Intuitionistic logic substitutes ``provability'' for ``truth'' in its logical calculus.
The logical calculus preserves provability, rather than truth, across transformations 
yielding derived propositions.

The essential difference to 
propositional logic 
is the interpretation of negation.
In propositional logic, $\neg P$ asserts that $P$ is false; in intuitionistic logic, 
$\neg P$ asserts that a proof of $P$ is impossible.  
The asymmetry is apparent: 
If $P$ is provable, then it is certainly impossible to prove that there is no proof of $P$;
however, 
we cannot conclude 
that there \emph{is} a proof of $P$, only
from an absence of a proof for a proof of the impossibility of $P$.

A direct consequence
is that many tautologies can no longer be proved within intuitionistic logic. 
An example is the \emph{tertium non datur} 
(Theorem \ref{theo-tertium-non-datur}, Example \ref{bsp-tertium}).
In intuitionistic logic, $P \vee \neg P=1$ says that at least one of 
$P$ or $\neg P$ can be proved, which is stronger than saying that their 
disjunction is true.
Similarly, \emph{Peirce's Law}\index{Pierce's law}
$((P \to Q) \to P) \to P = 1$ does not hold in intuitionistic logic.

From a practical point of view, there is a strong motivation for using 
intuitionistic logic. 
Solving problems in computer science,
one obviously is not interested in mere statements of existence.
A computer program is assumed to compute an answer, not to state that there is one. 
In applications one usually looks for a ``witness'' for a given existence 
assertion. In addition, one may have concerns about a proof system which has a proof 
for $\exists x: P(x)$, but which fails to prove $P(b)$ for any concrete $b$ it considers.

The observation that many classically valid tautologies are not theorems of 
intuitionistic logic leads to the idea of weakening the proof theory of
classical logic. This has for example been done by Gentzen with his 
``sequent calculus LK,'' obtaining a weaker version
that he called LJ \cite[§6.A]{Curry-1977}, \cite[§6]{Restall-2000}. 
A model theory can be given by Heyting algebras or, equivalently, by 
Kripke semantics, developed in the 1950's and 1960's.
The most natural mathematical framework for intuitionistic logic is the
algebraic concept of an implicative lattice, i.e.,
a lattice which is relatively pseudocomplemented and turns
out to be distributive. 

\subsection{Implicative lattices}
In connection with the foundations of logic, Brou\-wer and Heyting characterized an
important generalization of Boolean algebra, motivated by the following
consideration.
 In a Boolean algebra $A$, the complement $a'$ of an element $a\in A$
is the greatest element $x$ such that $a \wedge x = 0$,
i.e., such that $a$ and $x$ are ``disjoint''; more generally, $a \wedge x \leqq b$
if and only if $a \wedge x \wedge b' = 0$, that is
$(a \wedge b') \wedge x = 0$ or $x \leqq (a \wedge b')' = b \vee a' = a \to b$.
Hence, given $a$, $b\in A$, there exists a greatest element $c$ such that 
$a \wedge c \leqq b$.

\begin{definition}
	\label{def-implicative-lattice}
	A lattice $L$ is called
	\emph{implicative}\index{implicative lattice}, or
	\emph{Brouwerian}\index{Brouwerian lattice},
	if for any $a$, $b \in L$ the set $\{x \in L: a \wedge x \leqq b\}$ contains a
	greatest element, the 
	\emph{residuum}\index{residuum}, 
	\emph{relative pseudo-complement}\index{relative pseudo-complement}%
	\index{pseudo-complement}
	or \emph{material implication}\index{implication, material}%
	\index{material implication}
	$a \to b$ of $a$ in $b$.
	In a Brouwerian lattice with a universal bound $O$, 
	the element $a \to O$ is called a
	\emph{pseudocomplement} of $a$.
\end{definition}

Therefore, in an implicative lattice we have:
\begin{equation}
	\label{IP1}
	a \wedge (a \to b) \leqq b
	\qquad
	\mbox{(``modus ponens'')\index{modus ponens},}
\end{equation}
\begin{equation}
	\label{IP2}
	a \wedge c \leqq b 
	\Longleftrightarrow
	c \leqq a \to b,
	\qquad
	\mbox{(``residuation'')\index{residuation},}
\end{equation}
for any $a$, $b$, $c\in L$.\footnote{%
	The modus ponens follows directly from the definition of $a \to b$,
	as well as the $\Rightarrow$-direction of (\ref{IP2});
	the $\Leftarrow$-direction of (\ref{IP2}) follows immediately by modus ponens.
}
The material implication can be considered as
an operation $\to: L^2 \to L$.

\begin{lem}
	Let $L$ be an implicative lattice. Then the operation $\to: L^2 \to L$
	is inversely monotone with respect to its left argument and directly
	monotone with respect to its right argument. Further the following relations hold
	for all $a$, $b$, $c \in L$:
	\begin{equation}
		\label{eq-qi-1}
		b \leqq a \to b,
	\end{equation}
	\begin{equation}
		\label{eq-qi-2}
		a \to ( b \to c) \ =\ (a \wedge b) \to c \ =\ b \to (a \to c),
	\end{equation}
	\begin{equation}
		\label{eq-qi-3}
		a \to (b \to c) \leqq (a \to b) \to (a \to c),
	\end{equation}
	\begin{equation}
		\label{eq-qi-4}
		a \to (b \wedge c) = (a \to b) \wedge (a \to c),
	\end{equation}
	\begin{equation}
		a \to ( a \to b) = a \to b
		.
		\label{eq-qi-contraction}
	\end{equation}
	Eq.~(\ref{eq-qi-contraction}) is called 
	\emph{law of contraction}\index{contraction}\index{law of contraction}.
	If $L$ has a universal upper bound $1$, then
	\begin{equation}
		\label{eq-qi-5}
		a = 1 \to a,
	\end{equation}
	\begin{equation}
		\label{eq-qi-6}
		a \leqq b \mbox{ \ if and only if \ } 1 \leqq a \to b.
	\end{equation}
	If $L$ has a universal lower bound 0, then with the negation
	$\neg a := a \to 0$ we have
	\begin{equation}
		\label{eq-qi-7}
		\neg a \vee b \leqq a \to b,
		\qquad
		(\neg a \vee a) \wedge (a \to b) \leqq \neg a \vee b.
	\end{equation}
\end{lem}
\begin{proof}
	Assume $a \leqq b$. Then
	$ a \wedge ( b \to c) \leqq b \wedge (b \to c) \leqq c$
	by (\ref{L-B}) and (\ref{IP1}), hence $b \to c \leqq a \to c$ by
	(\ref{IP2}), proving the left monotonicity. Analogously,
	for $a \leqq b$ the right monotony follows from
	$ c \wedge ( c \to a) \leqq a \wedge b$
	by (\ref{IP1}), hence $c \to a \leqq c \to b$ by
	(\ref{IP2}).
	
	\emph{Proof of (\ref{eq-qi-1}):}
	By definition (\ref{def-wedge-vee}), we have $a \wedge b \leqq b$,
	hence $b \leqq a \to b$ by (\ref{IP2}).
	
	\emph{Proof of (\ref{eq-qi-2}):}
	By (\ref{IP1}), 
	$a \wedge (a \to ( b \to c)) \leqq b \to c$,
	hence 
	\[
		(a \wedge b) \wedge ( a \to (b \to c )) \leqq b \wedge (b \to c) \leqq c
	\]
	again by (\ref{IP1}) and the monotony. Then
	$a \to (b \to c) \leqq (a \wedge b) \to c$ by (\ref{IP2}). 
	This proves the left equation. The right equation follows by interchanging
	$a$ and $b$ and using the commutativity law (\ref{L2}) and the monotony.
	
	\emph{Proof of (\ref{eq-qi-3}):} We have
	\begin{eqnarray*}
		a \wedge (a \to b) \wedge (a \to (b \to c))
		& \leqq &
		a \wedge (a \to b) \wedge a \wedge (a \to (b \to c))
		\\
		& \leqq &
		b \wedge (b \to c) 
		\\
		& \leqq & c . 
	\end{eqnarray*}
	Here the first inequality follows by
	(\ref{L1}), (\ref{L2}), the second one follows by (\ref{IP1}) and monotony,
	and the last one by (\ref{IP1}).
	Hence
	$(a \to b) \wedge ( a \to (b \to c)) \leqq a \to c$ by (\ref{IP2}), and
	$a \to (b \to c) \leqq (a \to b) \wedge (a \to c)$ by (\ref{IP2}).
	
	\emph{Proof of (\ref{eq-qi-4}):} We have
	$a \to (b \wedge c) \leqq a \to b$ by (\ref{def-wedge-vee}) and monotony,
	and also
	$a \to (b \wedge c) \leqq a \to c$.
	Hence $a \to (b \wedge c) \leqq (a \to b) \wedge (a \to c)$ by (\ref{consistency}).
	Conversely, $a \wedge (a \to b) \wedge (a \to c) \leqq b \wedge c$ by (\ref{IP1})
	as in the proof for (\ref{eq-qi-3}). Hence 
	$(a \to b) \wedge (a \to c) \leqq b \wedge c$ by (\ref{IP2}).
	
	\emph{Proof of (\ref{eq-qi-contraction}):} This is the first equation in
	(\ref{eq-qi-2}) with $a=b$.
	
	\emph{Proof of (\ref{eq-qi-5}):} By (\ref{eq-qi-1}) we have $a \leqq 1 \to a$,
	and conversely $1 \to a = 1 \wedge (1 \to a)$, since by 
	definition of the universal bound, $1\wedge x = x$. Thus,
	$a \leqq 1 \to a \leqq a$.
	
	\emph{Proof of (\ref{eq-qi-6}):} 
	Suppose $a \leqq b$. Then $1 \wedge a \leqq b$ and
	$1 \leqq a \to b$ by (\ref{IP2}).
	Conversely, assume $1 \leqq a \to b$. 
	Then $a = a \wedge 1 \leqq b$ by (\ref{IP1}).
\end{proof}

\begin{satz}
	\label{satz-implikativ-distributiv}
	Any implicative lattice is distributive and a relative pseu\-do-com\-ple\-ment 
	in it is unique.
\end{satz}
\begin{proof}
	Given $a$, $b$, $c\in L$, form $d = (a \wedge b) \vee (a \wedge c)$ and consider
	$a \to d$. Since $a \wedge b \leqq d$, we have $b \leqq a \to d$ and 
	$c \leqq a \to d$.
	Hence $b \wedge c \leqq a \to d$, and so 
	$a \wedge (b \vee c) \leqq a \wedge (a \to d) \leqq d = (a \wedge b) \vee (a \wedge c)$.
	But this implies distributivity, by the distributive inequality
	(\ref{5}) and 	Theorem \ref{satz-equivalence-L6'-L6''}.
	Assume that $x$ and $y$ are relative pseudocomplements of $a$ in $b$. Then
	$x \wedge a = y \wedge a \leqq b$, since both are the greatest elements,
	hence by Theorem \ref{theo-distributive-uniqueness}, $x = y$.
\end{proof}

\begin{cor}
	An implicative lattice with universal bounds $0$ and $1$ is
	a distributive logic with the unique negation 
	$\neg a := a' = a \to 0$ for $a \in L$.
\end{cor}
\begin{proof}
	Let $L$ be an implicative lattice.
	With Theorem \ref{satz-implikativ-distributiv}, $L$ is distributive and has
	a unique pseudocomplement of $a$ relative to $b$ for all $a$, $b \in L$.
	In particular, $a'$ is unique, and 
	by the definition of material implication
	we have
	$0' = 0 \to 0 = 1$.
	Thus,
	\begin{equation}
		a \to (a')' 
		= a \to (a' \to 0) 
		= (a \wedge a') \to 0
		= 0 \to 0 = 1,
	\end{equation}
	where the second equality follows from (\ref{eq-qi-2}).
	With (\ref{eq-qi-6}) we thus have $a \leqq (a')'$.
	Moreover, for all $a$, $b$, $c \in L$ we have the relation
	\begin{equation}
		\label{eq-4C7}
		(a \vee b) \to c = (a \to c) \wedge (b \to c),
	\end{equation}
	since on the one hand,
	$(a \vee b) \to c \leqq a \to c$ by the antitony,
	and analogously
	$(a \vee b) \leqq b \to c$, i.e.,
	$(a \vee b) \leqq (a \to c) \wedge (b \to c)$,
	and on the other hand by distributivity,
	\begin{eqnarray*}
		\lefteqn{(a \vee b ) \wedge (a \to c) \wedge (b \to c)}
		\\
		& = & 
		a \wedge (a \to c) \wedge (b \to c)
		\vee b \wedge (a \to c) \wedge (b \to c)
		\\
		& \leqq & 
		(c \wedge (b \to c)) \vee (c \wedge (a \to c))
		\qquad \mbox{(by (\ref{IP1})} 
		\\ 
		& = & 
		c \wedge ((b \to c) \vee (a \to c))
		\leqq c,
	\end{eqnarray*}
	thus $(a \to c) \wedge (b \to c) \leqq (a \vee b) \to c$ by (\ref{IP2}).
	By (\ref{eq-4C7}) for $c=0$ we have
	$(a \vee b)' = (a \to b) \to c = (a \to c) \wedge (b \to c) = a' \wedge b'$,
	i.e., the disjunctive De Morgan law (\ref{disjunctive-De-Morgan}).
	Hence $L$ is a logic.
\end{proof}

\begin{satz}
	\label{theo-implicative-strong-double-negation}
	In an implicative lattice $L$ with universal bounds $0$ and $1$
	and the pseudocomplement defined by $\neg a := a' = a \to 0$,
	the following relations are equivalent for all $a$, $b\in L$.
	\begin{eqnarray}
		\label{eq-impl-3}
		\index{stable negation}\index{negation! stable -}%
		\index{negation! strong double -}\index{strong double negation}%
		\index{double negation! strong -}%
		\mbox{(Strong double negation, stability law)} & &
		a = \neg \neg a
		\\
		\label{eq-impl-1}
		\mbox{(tertium non datur)}\index{tertium non datur} & &
		1 = a \vee \neg a, 
		\\
		\label{eq-impl-4}
		\mbox{(Peirce's law)}\index{Pierce's law} & &
		a = (a \to b) \to a,
		\qquad
		\\
		\label{eq-impl-2}
		& &
		1 = a \vee (a \to b).
	\end{eqnarray}
\end{satz}
\begin{proof}	
	(\ref{eq-impl-3}) $\Rightarrow$ (\ref{eq-impl-1}):
	Since $L$ is a logic, this follows by Theorem \ref{theo-tertium-non-datur}.

	(\ref{eq-impl-1}) $\Rightarrow$ (\ref{eq-impl-3}):
	By the distributivity law and (\ref{double-negation}) it follows
	\[
		\neg \neg a
		= \neg \neg a \wedge (a \vee \neg a) 
		= (\neg \neg a \wedge a) \vee (\neg \neg a \wedge \neg a)
		= a \vee 0
		= a.
	\]

	(\ref{eq-impl-1}) $\Rightarrow$ (\ref{eq-impl-4}):
	By (\ref{eq-qi-1}) we have
	$(a \to b) \to a \leqq \neg a \to ((a \to b) \to a)$,
	i.e, by (\ref{eq-qi-2}) and (\ref{eq-qi-4})
	\[
		(a \to b ) \to a
		\leqq ( \neg a \wedge (a \to b)) \to a
		= ((a \to 0 ) \wedge (a \to b)) \to a
		= \neg a \to a,
	\]
	hence by (\ref{eq-qi-1})
	\begin{equation}
		a \leqq (a \to b) \to a \leqq \neg a \to a.
		\label{eq-pre-Peirce}
	\end{equation}
	Since by (\ref{eq-impl-1}) we deduce from (\ref{eq-qi-7}) that
	$a \to b = \neg a \vee b$, we have
	$\neg a \to a = a$, and by (\ref{eq-pre-Peirce}) we have (\ref{eq-impl-4}).
	
	(\ref{eq-impl-4}) $\Rightarrow$ (\ref{eq-impl-2}):
	First we see that 
	\begin{equation}
		\label{eq-pre-tertium}
		((a \vee (a \to b)) \to a ) \to (a \vee (a \to b)) 
		= a \vee (a \to b)
	\end{equation}
	by Peirce's law (\ref{eq-impl-4}). Moreover,
	by (\ref{eq-qi-1}) we have $a \leqq (a \vee (a \to b)) \to a$,
	but by the antitony of $\to$ with respect to its left argument,
	$(a \vee (a \to b)) \to a \leqq (a \to b) \to a = a$,
	applying Peirce's law again, hence
	$(a \vee (a \to b)) \to a = a$, and insertion in (\ref{eq-pre-tertium}) yields
	\begin{equation}
		a \to (a \vee (a \to b)) = a \vee (a \to b). 
	\end{equation}
	By (\ref{eq-qi-6}) we therefore have
	$1 \leqq a \to a = a \to (a \vee (a \to b)) = a \vee (a \to b).$
	
	(\ref{eq-impl-2}) $\Rightarrow$ (\ref{eq-impl-1}):
	Set $b = 0$ in (\ref{eq-impl-2}).
\end{proof}

It is easily verified that any Boolean algebra is an implicative lattice, in which
$a \to b = a' \vee b$ is the relative complement of $a$ in $[a \wedge b, 1]$.
Likewise, any finite distributive lattice is implicative, since the join
$u = \bigvee x_j$ of the $x_j$ such that $a \wedge x_j \leqq b$ satisfies
$a \wedge u = a\wedge \bigvee x_j = \bigvee ( a \wedge x_j) \leqq b$.
In \cite[§4C]{Curry-1977}, implicative lattices and the dual
``subtractive'' lattices are summarized under the notion 
``Skolem lattice''\index{Skolem lattice}.
Since analogously to the implicative case, every finite distributive lattice
is subtractive, we conclude that any finite distributive lattice is a Skolem
lattice.

The simplest non-Boolean implicative lattice is the Gödel logic $\mbox{G}_3$
given in Example \ref{bsp-Goedel-3}.


\begin{beispiel}
\textbf{(Heyting algebra of open sets of $\mathbb{R}^2$)}
Let $\mathscr{O}(\mathbb{R}^2, \subseteq)$ be the po\-set 
of open sets of the plane $\mathbb{R}^2$, 
and let
$\wedge$ and $\vee$ denote the 
operations correspond to set intersection and union, respectively.
Then $A\to  B$ $=$ $(A^C)^\circ \cup B$
is the union of the set $B$ and the inner set $(A^C)^\circ$ 
of the complement $A^C = \mathbb{R}^2 \setminus A$ of the set $A$.
Therefore, $\mathscr{O}(\mathbb{R}^2, \subseteq)$ is a lattice, with
universal bounds $0=\emptyset$ and $1=\mathbb{R}^2$.
With the pseudocomplement
$\neg A := (A^C)^\circ$ as negation, 
$(\mathscr{O}(\mathbb{R}^2, \subseteq), \neg)$ is even a distributive logic.
Especially, the formula $A \wedge \neg A = 0$ is valid, since
$A \cap (A^C)^\circ = \emptyset$.
However, \emph{tertium non datur} is not valid, as is seen
by letting $A = \{(x,y): y > 0 \}$ be the upper half plane:
$\neg A = \{(x,y) : y \leqq 0 \}^\circ = \{(x,y) : y < 0 \}$, 
and  $A\vee \neg A = \{(x,y)\in \mathbb{R}^2 : y \ne 0 \} \ne \mathbb{R}^2$, 
i.e.,  $A\vee \neg A < 1$.
%
\end{beispiel}

In general, the complete distributive lattice of all open subsets of any topological
space is implicative. However, the complete distributive lattice of all
\emph{closed} subsets of the line is not implicative: there is no greatest closed
set satisfying $p \wedge x = \emptyset$.
Therefore, not all distributive lattices are implicative,
distributivity is necessary but not sufficient for a lattice to
be implicative.
Nondistributive logics like quantum logics thus are not implicative lattices.

\begin{satz}[\textbf{Curry's Paradox \cite{Curry-1942}\index{Curry's paradox}}]
	\label{theo-Curry}
	Let $L$ be an implicative lattice with a universal upper bound $1\in L$,
	and let $y \in L$ be an arbitrary proposition. Suppose moreover
	that the mapping $f_y:L \to L$, 
	\begin{equation}
		f_y(x) = (x \to y)
		\label{eq-mapping-implication}
	\end{equation}
	has a fixed point $x_* \in L$, i.e., that there exists a proposition 
	$x_* \in L$ such that $x_* = f_y(x_*)$.
	Then both $x_*$ and $y$ are true, i.e., $x_*=y=1$.
\end{satz}
\begin{proof}
	First we notice that by the residuation condition (\ref{IP2}) we have
	$1 = (x_* \to x_*)$, setting $a=b=x_*$ and $c=1$. Since $x_*$ is a fixed point,
	we have $x_* \to f_y(x_*)$, hence
	\begin{equation}
		1 = x_* \to (x_* \to y) 
		\stackrel{\mathrm{(\ref{eq-qi-3})}}{=}
		(x_* \wedge x_*) \to y
		= x_* \to y.
		\label{eq-Curry-1}
	\end{equation}
	Thus $f_y(x_*) = 1$, i.e., $x_*=f_y(x_*)=1$. By \emph{modus ponens}
	(\ref{IP1}), with $a=1$ and $b=y$, this implies $1 \leqq y$.
\end{proof}

A paradox arises in an implicative lattice with a universal lower bound 0
since $y$ can be \emph{any} proposition, especially a false proposition. 
A natural language version
of Curry's paradox reads ``If this sentence is true, then Santa Claus exists.''
Here $x$ $=$ ``this sentence is true'' and $y$ $=$ ``Santa Claus exists,''
and thus $y$ is true.
However, in the version above the paradox is resolved since it is formulated
completely on the object language level:

\begin{cor}
	In an implicative lattice with universal upper bound $1$, the mapping
	$f_y$ for $y<1$ in Eq.~(\ref{eq-mapping-implication}) has no
	fixed point.
\end{cor}


Curry's paradox is a serious challenge in naive truth theory 
with an unrestricted $T$-schema $T[x] \leftrightarrow x$.
In a similar fashion as Theorem \ref{theo-Curry} it then can be
proved that any sentence $y$ can be derived from this schema.
In natural language, $T[x]$ $=$ ``This sentence $x$ is derivable,''
and $T[x] \leftrightarrow x$. Concerning provability instead of truth,
the paradox is known as \emph{Löb's Theorem} \cite[§18]{Boolos-et-al-2002}.

\section{Boolean algebras}
A complemented distributive
lattice is called
\emph{Boolean lattice}\index{Boolean lattice},
or
\emph{Boolean algebra}\index{Boolean algebra}.
In a Boolean algebra $L$, the notation is usually modified.
The 
complement $x'$ is replaced by $\neg x$\index{$\neg$} and called 
\emph{negation}\index{negation}.
The universal bounds are denoted simply by 0 and 1.

The fact that it can be regarded as an algebra is justified
by the first part of following statement which asserts
that in a Boolean lattice complements are unique.

\begin{satz}
	In a Boolean lattice $L$ each element $x$ $\in$ $L$
	has one and only one complement, and the 
	complementation is an involutive negation%
	\index{involutive negation}\index{negation! involutive -},
	i.e.,
	\begin{equation}
		\label{BL10}
		(x')' = x.
	\end{equation}
\end{satz}
\begin{proof}
	Assume that for an element $x\in L$ in the lattice there exist
	two complements $x'$, $y$. Then we have
	$x \wedge x'$ $=$ $x \wedge y$ $=$ $O$, and 
	$x \vee x'$ $=$ $x \vee y$ $=$ $I$. 
	With Theorem \ref{theo-distributive-uniqueness}
	this implies $x' = y$. Therefore, $x \mapsto x'$ is single-valued.
	But by the symmetry of the definition \ref{def-complement}
	of complement, $x$ is a complement of $x'$, hence
	$x=(x')'$ by uniqueness, proving (\ref{BL10}). Thus the correspondence
	$x \mapsto x'$ is one-one.
\end{proof}

By the definition \ref{def-complement} of complement, we have
$x = x \wedge I = x \wedge (y \vee y')$, i.e., 
in a Boolean algebra by the distributive law
\begin{equation}
	x = ( x \wedge y) \vee ( x \wedge y')
	\qquad
	\mbox{for all } x, y \in L.
\end{equation}
In a Boolean algebra there can be defined the addition $\oplus$ 
(i.e., addition modulo 2)
and multiplication $\cdot$ by
\begin{equation}
	x \cdot y := x \wedge y,
	\qquad
	x \oplus y := (x \wedge y') \vee (x' \wedge y)
\end{equation}
Then the $\vee$-operation and the negation are related to the addition by
\begin{equation}
	x \vee y = x \oplus y \oplus xy,
	\qquad
	x' = 1 \oplus x,
\end{equation}
where we assign $O$ by 0 and $I$ by 1.
With multiplication and addition defined in this way, 
the Boolean algebra is a ``Boolean ring with unit.''

\begin{satz}
	Any complete atomic Boolean lattice is isomorphic to
	${\mathbf{2}}^\aleph$, where $\aleph$ is the cardinality of the set of
	its atoms.
\end{satz}
\begin{proof}
	\cite[\S\,VIII.9]{Birkhoff-1973}.
\end{proof}

\begin{satz}
	A necessary and sufficient condition that 
	$x = y$
	hold in a Boolean algebra is that $x$ and $y$ have the same value
	in every evaluation of $0$-$1$ truth tables.
	Analogously, a necessary and sufficient condition for
	$x \leqq y$
	is that $y$ has the value $1$ in every evaluation of $0$-$1$ tables
	in which $x$ has the value $1$.
\end{satz}
\begin{proof}
	\cite[Theorem 6 in \S6.D]{Curry-1977}.
\end{proof}

For a Boolean algebra, this theorem gives a theoretical solution to the 
decision problem\index{decision problem} wheth\-er a given ``well-formed expression'' 
$x$ is provable. 
However, it is not always the fastest method.
If the number of indeterminates is $n$, then there are $2^n$ possibilities
to be considered.
The method of reduction by translating into a Boolean ring with  unit,
multiplying out, and using addition modulo 2 is faster.

\subsection{Propositional logic}
Boolean algebra especially applies to propositions.
Thus, for any two propositions $x$ and $y$, one can denote the
propositions ``$x$ and $y$,'' ``$x$ or $y$,'' and ``not $x$''
by $x\wedge y$, $x \vee y$, and $x'$, respectively.
With respect to this interpretation we conclude:
\emph{Propositions form a Boolean algebra}\index{proposition}
(``Boole's third law'').
Hence propositional logic is a Boolean logic.
In addition to the properties of a Boolean algebra as discussed above,
in classical, ``two-valued'', logic\index{logic}, all propositions
are either true or false, and never both.
Moreover, $x \wedge y$ is true if and only if $x$ and $y$ are both true;
$x\vee y$ is true when $x$ or $y$ is true; of $x$ and $x'$, one is true
and the other false.

Under these assumptions, the compound proposition ``$x$ implies $y$''
(i.e., ``if $x$ then $y$''), denoted $x \to y$, has a special meaning.
It is true or false according to ``$y$ or not-$x$'' is true or
false. Hence one can interpret
$x \to y$ meaning $x' \vee y$, i.e.,
\begin{equation}
	x \leqq y \ = \ x' \vee y.
\end{equation}
In two-valued logics, one can similarly symbolize the statement
``$x$ is equivalent to $y$'' by $x$ $\leftrightarrow$ $y$, and replace
it by ``$x$ implies $y$ and $y$ implies $x$,'' i.e., with the
symmetric addition operation of \emph{complemented} Boolean ring:
\begin{equation}
	x \leftrightarrow y \ = \ (x' \vee y ) \vee (x \wedge y')
	 \ = \
	 (x \oplus y )'.
\end{equation}
Many compound propositions $x$ are ``tautologies''\index{tautology},
that is, true by virtue of their logical structure alone. This
amounts algebraically to saying that $x \to 1$.
The simplest tautology is $x \vee x'$ (``$x$ or not $x$'').
It is a simple exercise in Boolean algebra to show that the following
propositions are all tautologies of two-valued logic \cite[§XII.2]{Birkhoff-1973}:
\begin{equation}
	\begin{array}{c}
		x \vee x', \quad 0 \to x, \quad x \to 1, \quad x \to x, \quad
		x \leftrightarrow x, \quad x \to ( y \to x),
		\\
		(x \leftrightarrow y) \leftrightarrow (y \leftrightarrow x), \quad
		(x \leftrightarrow 0) \vee (x \leftrightarrow 1), \quad
		(x \leftrightarrow y) \vee (y \to x).
	\end{array}
\end{equation}

\paragraph{Critique.}
\label{sec-propositional-calculus}
There often has been raised intuitive objections to the preceding logical
rules, pushed to extremes. For instance, the tautology
$0 \to x$ asserts that ``a false proposition implies every proposition''
(``ex falso quodlibet'').
But what does this really mean\footnote
{%
	Bertrand Russell\index{Russell} 
	is reputed to have been challenged
	to prove that the false hypothesis $2+2=5$ implied that he was the
	Pope. Russell replied: ``You admit $2+2=5$; but I can prove $2+2=4$;
	therefore 5 $=$ 4. Taking two away from both sides, we have 3 $=$ 2;
	taking one more, 2 $=$ 1. But you will admit that I and the Pope are
	two. Therefore, I and the Pope are one. Q.E.D.''
}%
?
Likewise, one may question the validity of the tautology
$(x \to y) \vee (y \to x)$, which asserts that ``of any two propositions
$x$ and $y$, either $x$ implies $y$ or $y$ implies $x$.''

Most important, intuition can make one skeptical of the validity of
``proofs by contradiction,'' or ``reductio ad absurdum,'' 
on which many mathematical proofs
are based upon. The idea of this principle is simple. Since
$(x \to y)$  $=$ $x' \vee y$, and $(y' \to x')$ $=$ $y \vee x'$,
by the symmetry of the $\vee$-operation 
we have 
\begin{equation}
	(x \to y) \ = \ (y' \to x'), 
\end{equation}
i.e.,
``$x$ implies $y$'' is the same as
``not-$y$ implies not-$x$.'' 
So if one does not succeed in proving the left hand side of this equation, 
one can try to prove the right hand side (which is often simpler,
since it is not \emph{constructive}).
But why should the disproof ``not-$x$'' imply the truth of $x$?
See Example \ref{bsp-tertium}.

However, propositional logic is ``complete,''\index{complete}
i.e., in propositional logic any Boolean formula $\phi$ 
is derivable from a set $\Delta$ of Boolean formulas, if and only if,
for any variable assignment,
$\phi$ is true whenever each formula in $\Delta$ is true,
in symbols $\Delta \vdash \phi$ if and only if $\Delta \models \phi$
\cite[§1.9]{Hedman-2004}.

\subsection{First-order logic}
First-order logic\index{first-order logic}\index{predicate logic},
sometimes also called \emph{predicate logic},
is an extension of propositional logic containing the 
\emph{universal quantifier}\index{existential quantifier}\index{quantifier}
$\forall$ (``for all'') and the
\emph{existential quantifier}\index{existential quantifier}
$\exists$ (``there exists'').

The basis of first-order logic is given by a
\emph{term}\index{term} which is recursively defined by the following rules:
\begin{enumerate}
\item
Any constant symbol $a$, $b$, $c$, \ldots is a term with no free variables.
\item
Any variable $x$, $y$, $z$, \ldots is a term whose only free variable is itself.
\item
Any expression $f(t_1$, \ldots, $t_n)$ of $n \geqq 1$ arguments is a term,
where each argument $t_i$ is a term and $f$ is a function symbol of arity $n$,
whose free variables are the free variables of any of the terms $t_i$.
\end{enumerate}
A \emph{relation}\index{relation}, or \emph{predicate variable}\index{predicate},
is then an $n$-ary relation $S=S(t_1,$ \dots, $t_n)$ of terms.
A relation $P(t_1$, \ldots, $t_n)$ then forms an 
\emph{atomic formula}. If in addition it contains no free variable 
it is an \emph{atomic sentence}. Atomic sentences
play the role of the basic propositions of the first-order logic.
Its \emph{vocabulary}\index{vocabulary} $\Sigma$ then is the
set of relations, functions, and constant symbols.

\begin{beispiel}
	The vocabulary of ordered abelian groups has a 
	constant 0, a unary function $-$, a binary function $+$, 
	and two binary relations $=$ and $\leqq$, so
	$\Sigma$ $=$ \{0, $-$, $+$, $=$, $\leqq$\}, and
	\begin{itemize}
	\item
	$0$, $x$, $y$ are atomic terms;
	\item
	$+(x, y)$, $+(x, +(y, -(z)))$ are terms, usually written as 
	$x + y$, $x + y - z$;
	\item
	$=(+(x, y), 0)$, $\leqq (+(x, +(y, -(z))), +(x, y))$ are atomic formulas, 
	usually written as $x + y = 0$, $x + y - z \leqq x + y$.
	\end{itemize}
	The relations, or predicates, $=$ and $\leqq$ are Boolean-valued
	if constants are inserted. 
\end{beispiel}

Given a vocabulary $\Sigma$ of constant symbols, functions, and relation,
the 
\emph{universal quantifier}\index{universal quantifier}\index{quantifier}
$\forall$ (``for all'') and the
\emph{existential quantifier}\index{existential quantifier}
$\exists$ (``there exists'') can be introduced, 
acting on formulas $P(\ldots, x, \ldots)$ with a free variable $x$
by writing $\forall x\, P(x)$ or $\exists x\, P(x)$, respectively.
They have to satisfy the axioms
\begin{itemize}
	\item[\textbf{P1}:] $(\forall x\, P(x)) \to P(t)$ for any term $t$
	without free variables.
	\item[\textbf{P2}:] $P(t) \to (\exists x\, P(x))$ for any term $t$
	without free variables.
	\item[\textbf{P3}:] $(\forall x\, (\phi \to P(x))) \ \to \ (\phi \to \forall x\, P(x))$
	for any atomic sentence $\phi$.
	\item[\textbf{P4}:] $(\forall x\, (P(x) \to \phi)) \ \to \ (\exists x\, P(x) \to \phi)$
	for any atomic sentence $\phi$.
\end{itemize}

In fact, the two quantifiers $\forall$ and $\exists$ are not independent. Since
\begin{equation}
	\fbox{$
	\exists = \neg \forall \neg,
	\qquad
	\forall = \neg \exists \neg,
	$}
\end{equation}
both quantifiers are dual to each another. Thus in fact, the introduction
of only one of them suffices to extend propositional logic to first-order logic.

A \emph{sentence}\index{sentence} of first-order logic is a formula having no
free variables, i.e., only quantified variables.
For instance, $\exists x$ $\forall y$ $x\cdot y = x$ is a sentence of first order logic,
whereas $\exists x$ $x\cdot y = x$ is a formula but not a sentence,
because $y$ is a free variable.
As a primality formula for a natural number, we could define
\begin{equation}
	\pi(x) = 
	\left( 
		\left(1 < x \right)
		\ \wedge \
		\forall u \forall v\ 
		\left( \left(x = u \cdot v\right) \ \rightarrow \ 
		\left( u = 1 \vee v = 1 \right) \right)
	\right)
\end{equation}
with the free variable $x$. The formula $\pi(x)$ is true if and only if $x$ is a prime.
This is not a sentence. However, the assertion
$
	\forall z \exists x \ (x > z \wedge \pi(x))
$
that there are infinitely many primes is a sentence.

\begin{beispiel}
	The vocabulary of ordered abelian groups has a 
	constant 0, a unary function $-$, a binary function $+$, 
	and two binary relations $=$ and $\leqq$, so
	$\Sigma$ $=$ \{0, $-$, $+$, $=$, $\leqq$\}, and
	\begin{itemize}
	\item
	$0$, $x$, $y$ are atomic terms;
	\item
	$+(x, y)$, $+(x, +(y, -(z)))$ are terms, usually written as 
	$x + y$, $x + y - z$;
	\item
	$=(+(x, y), 0)$, $\leqq (+(x, +(y, -(z))), +(x, y))$ are atomic formulas, 
	usually written as $x + y = 0$, $x + y - z \leqq x + y$.
	\item
	$(\forall x$ $\exists y$ $\leqq(+(x,y)), z)$
	$\wedge$ $(\exists x$ $=(x,y) = 0))$
	is a formula, usually written as
	$(\forall x$ $\exists y$ $x + y$ $\leqq$ $z)$ $\wedge$ $\exists$ $y$ $x+y = 0)$.
	\end{itemize}
\end{beispiel}

To summarize, first-order logic contains the nine fixed symbols
\begin{equation}
	\wedge, \vee, \neg, \rightarrow, \leftrightarrow, (, ), \exists, \forall.
\end{equation}

Like propositional logic, first-order logic is complete
\cite[§V.4]{Ebbinghaus-et-al-1996}, \cite[§4]{Hedman-2004}. 
But first-order logic is, in some sense, weak,
although it is much richer than propositional logic.
Whereas it can express ``$\exists$ $n$ elements'' for any finite $n$,
it cannot express ``$\exists$ countably many elements.''
This is possible only with a ``second-order logic'' which allows
existential expressions about sets.

\subsubsection{Models}
The truth of a Boolean expression $\phi$ in propositional logic is computed
by a \emph{truth assignment} $T$ being a mapping $T:X' \to \{0,1\}$ from
a finite subset $X'\subseteq X=\{x_1, x_2, \ldots\}$ of a countably infinite alphabet
$X$ of Boolean variables \cite[§4.1]{Papadimitriou-1994}; 
then $T$ is said to \emph{satisfy} the Boolean expression $\phi$, in symbols 
$T \models \phi,$
if $T(\phi) = 1$ \cite[§1.2]{Hedman-2004}.

However, in first-order logic 
variables, functions, and relations can take much more complex values
than just 0 and 1, or \code{false} and \code{true}.
The analog of a truth assignment for first-order logic is a far more complex mathematical
object called a ``model.''
Let $\Sigma$ be a given vocabulary.
A \emph{$\Sigma$-structure}\index{structure}, or \emph{$\Sigma$-model}\index{model},
$M$ is a pair $(U| \mu)$ where $U$ is a nonempty set $U$, the \emph{universe} of $M$, 
and where $\mu$ is a function called the 
\emph{interpretation}\index{interpretation} of $\Sigma$ 
assigning to each symbol in the vocabulary $\Sigma$ a respective object of the
universe $U$:
\begin{itemize}
	\item
	to each constant $c\in\Sigma$ an element in $U$, i.e.,
	$\mu(c) \in U$;
	
	\item
	to each $n$-ary function $f \in \Sigma$ a function $g:U^n \to U$, i.e.,
	$\mu(f) = g$;
	
	\item
	to each $n$-ary relation $R \in \Sigma$ a subset $S \subset U^n$, i.e.,
	$\mu(R) = S$.
\end{itemize}
Suppose now $M$ to be a $\Sigma$-structure and $\phi$ a sentence over $\Sigma$, 
then we define $M$ to 
\emph{model}
, or \emph{satisfy}, 
$\phi$, in symbols
$ 
	M \models \phi,
$ 
if the interpretation of the sentence $\phi$ as a formula in $M$ is true
\cite[§2.3]{Hedman-2004}.
Moreover, let $\Delta$ be a set of sentences over $\Sigma$. Then we write
\begin{equation}
	\Delta \models \phi
\end{equation}
if any model which
models 
all sentences in $\Delta$, also
models
$\phi$. 
For any $\Sigma$-structure $M$, the \emph{theory of $M$}\index{theory}
is the set $T_M$ of all $\Sigma$-sentences $\phi$ such that $M \models \phi$.
A theory $T$ is \emph{decidable}\index{decidable theory} if there exists
an algorithm which determines in a finite number of steps whether or no a given
$\Sigma$-sentence is in $T$.

\begin{beispiel}
	\label{bsp-number-systems}
	\emph{(Number systems)}
	\cite[§2.4.4]{Hedman-2004}
	Consider the vocabulary of natural numbers, $\Sigma = \{+, \cdot, 1\}$
	having binary functions $+$ and $\cdot$, written as usual in infix-notation,
	e.g., ``$x + y$'' instead of ``+(x,y),'' and constants 0 and 1.
	We let 2 abbreviate $(1+1)$, $x^2$ abbreviate $x\cdot x$, 3 abbreviate $1+(1+1)$,
	and so on. Then we define the $\Sigma$-model $\fett{N}$ as
	\begin{itemize}
		\item
		$\fett{N} = (\mathbb{N}| +, \cdot, 1)$
	\end{itemize}
	Moreover, consider the vocabulary of arithmetic, $\Sigma = \{+, \cdot, 0, 1\}$,
	and the following $\Sigma$-structures:
	\begin{itemize}
		\item
		$\fett{A} = (\mathbb{Z}| +, \cdot, 0, 1)$
		\item
		$\fett{Q} = (\mathbb{Q}| +, \cdot, 0, 1)$
		\item
		$\fett{R} = (\mathbb{R}| +, \cdot, 0, 1)$
		\item
		$\fett{C} = (\mathbb{C}| +, \cdot, 0, 1)$
	\end{itemize}
	Each of these structures interpretes the symbols in $\Sigma$ in the usual way.
	Then any polynomial with natural numbers as coefficients is a $\Sigma$-term.
	Equations such as
	$$
		x^5 + 9x + 3 = 0
	$$
	are $\Sigma$-formulas. 
	(Note that 3 and $x^5$ are not symbols in $\Sigma$,
	they are abbreviations for the $\Sigma$-terms $1+(1+1)$ and 
	$x \cdot (x \cdot (x \cdot (x \cdot x)))$, respectively.)
	
	Then the sentence $\phi$ $=$ $\exists x$ $(x^2 = 2)$ asserts the existence of $\sqrt2$.
	It follows $\fett{Q} \models \neg \phi$ and $\fett{R} \models \phi$, i.e.,
	$\fett{R}$ models $\phi$ but $\fett{Q}$ does not. 
	Likewise, the equation $2x+3 = 0$ has a solution in $\fett{Q}$ but not in $\fett{A}$,
	so 
	$\fett{Q} \models \exists x (2x+3=0)$, whereas 
	$\fett{A} \models \neg\exists x (2x+3=0)$.
	
	To progress from $\mathbb{N}$ to $\mathbb{Z}$ to $\mathbb{Q}$ etc, we add solution
	for more and more polynomials, reaching the end of the line with the complex numbers
	$\mathbb{C}$. The sentence $\exists x(x^2 + 1 = 0)$ distinguishes $\fett{C}$ from
	all the other $\Sigma$-structures in the list.
	The Fundamental Theorem of Algebra states that for any nonconstant polynomial
	$p(x)$ having coefficients in $\mathbb{C}$, the equation $p(x) = 0$ has a solution
	in $\mathbb{C}$. So there is no need to extend to a bigger system,
	by virtue of adding a solution of $x^2+1=0$ to $\mathbb{R}$, we added a solution
	to every polynomial.
	The names of the number systems, starting from the ``natural'' numbers,
	extended by ``negative'' numbers to the integers, further extended by even ``irrational''
	numbers and at last by ``imaginary'' numbers suggests that the systems get
	more and mor complicated as we progress from natural to complex numbers.
	From the point of view of first-order logic, however, this is backwards.
	The structure $\fett{C}$ is the most simple \cite[§5]{Hedman-2004}, 
	whereas the structures $\fett{A}$,
	standing for arithmetic, and $\fett{N}$ are most complex \cite[§8]{Hedman-2004}.
\end{beispiel}

\subsubsection{Proofs and axioms}
Define $\Lambda$ to be the set of logical axioms\index{axiom}\index{logical axioms}
\cite[§5.4]{Papadimitriou-1994}.
Let $\phi$ be a first-order formula and let $\Delta$ be a set of first-order formulas.
Moreover, let $S = \{\phi_1, \phi_2, \ldots, \phi_n\}$ be a finite 
sequence of first-order formulas such that for each formula $\phi_i \in S$,
1 $\leqq$ $i$ $\leqq$ $n$,
we have either (a) $\phi_i \in \Lambda$, or (b) $\phi_i \in \Delta$, 
or (c) there are two expressions 
$\psi$, $(\psi \Rightarrow \phi_i)$ $\in$ $\{\phi_1, \ldots, \phi_{i-1}\}$ $\subseteq$ $S$.
Then, $S$ is called a \emph{(formal) proof}\index{proof} of $\phi = \phi_n$ 
from $\Delta$, and we write
\begin{equation}
	\Delta \vdash \phi.
\end{equation}
Then $\phi$ is called a \emph{$\Delta$-first-order theorem}\index{theorem}.
Hence, a $\Delta$-first-order theorem is an ordinary first-order theorem if
$\Delta = \emptyset$ and we only have to satisfy the logical axioms $\Lambda$.
In other words, a $\Delta$-first-order theorem would be an ordinary first-order theorem
\emph{if} we allowed all formulas in $\Delta$ to be added to our logical axioms $\Lambda$.
This is the reason why the formulas in $\Delta$ are often called
\emph{nonlogical axioms}\index{nonlogical axioms}.

\begin{beispiel}
	\emph{(Group theory)}
	\cite[Ex.\,5.5]{Papadimitriou-1994}
	Let be $\Sigma=\{\circ, 1\}$ and $\Delta$ $=$ \{GT1, GT2, GT3\} where GT1, GT2, GT3
	are the nonlogical axioms
	\begin{equation}
		\begin{array}{l@{\ =\ }ll}
			\mbox{GT1} & 
			\forall x \forall y \forall z \ ((x \circ y) \circ z = x \circ (y \circ z))
			& \mbox{(associativity of $\circ$),}
			\\[.5ex]
			\mbox{GT2} & 
			\forall x \ (x \circ 1 = x)
			& \mbox{(1 is the neutral element),}
			\\[.5ex]
			\mbox{GT3} & 
			\forall x \exists y \ (x \circ y = 1)
			& \mbox{(existence of inverses).}
		\end{array}
	\end{equation}
	These three simple axioms comprise a complete axiomatization of all groups.
	All properties of groups can be deduced from these axioms by formal proofs.
	If we want to axiomatize Abelian groups, we have to add another axiom
	\begin{equation}
		\mbox{GT4 $=$ } \forall x \forall y \ (x \circ y = y \circ x)
	\end{equation}
	On the other hand, if we want to study infinite groups, it suffices to add
	for each $n>1$ the sentence
	\begin{equation}
		\phi_n = \exists x_1 \exists x_2 \cdots \exists x_n 
		\bigwedge_{i \ne j} (x_i \ne x_j).
	\end{equation}
	This infinite sequence of sentences is a complete axiomatization of infinite groups.
	Thus group theory is axiomatizable, but is not \emph{decidable},
	unless we restrict to Abelian groups.
	(Tarski showed in 1946 that any statement in Peano arithmetic can be encoded as a 
	statement in group theory, thus demonstrating that group theory is universal, 
	and that questions about it can be undecidable.)
	In contrast, Gödel's incompleteness theorem is based on the fact that, if number
	theory \emph{were} axiomatizable, then it would be decidable.
\end{beispiel}

\begin{beispiel}
	\label{bsp-axiomatization-number-theory}
	\emph{(Number theory)}
	\cite[§8.1]{Hedman-2004}
	Let the vocabulary $\Sigma_N = \{+, \cdot, 1\}$ be given.
	The theory $T_N$ of the model $\fett{N}$ defined in Example \ref{bsp-number-systems}
	is then first-order number theory.
	If $T_N$ had a decidable axiomatization, then in principle we could use
	formal proof theory to answer every open number theoretic question.
	However, by Gödel's First Incompleteness Theorem, this is not the case.
	Since $T_N$ is undecidable, it does not have a decidable first-order axiomatization
	\cite[Prop.\,5.10]{Hedman-2004}. But there exists a \emph{second-order} theory
	containing $T_N$ which does have a decidable axiomatization.
	For this purpose, we define the set $\Delta$ $=$ \{NT1, NT2, \ldots, NT8\}
	of the following (nonlogical) axioms:
	\[ 
		\begin{array}{l@{:\  }ll}
			\mbox{NT1} & 
			\forall x \ \ \neg (x+1 = 1)
			\\ 
			\mbox{NT2} & 
			\forall x \forall y \ (x + 1 = y + 1 \ \rightarrow \ x = y)
			\\ 
			\mbox{NT3} & 
			\forall x \forall y \ (x + y = y + x)
			& \mbox{(addit.\ commutativity)}
			\\ 
			\mbox{NT4} & 
			\forall x \forall y \ (x + (y+1) = (x + y) + 1)
			& \mbox{(addit.\ associativity)}
			\\ 
			\mbox{NT5} & 
			\forall x \ (x \cdot 1 = x)
			& \mbox{(multipl.\ neutral element),}
			\\ 
			\mbox{NT6} & 
			\forall x \forall y \ (x \cdot y = y \cdot x)
			& \mbox{(multipl.\ commutativity)}
			\\ 
			\mbox{NT7} & 
			\forall x \forall y \ (x \cdot (y+1) = (x \cdot y) + x)
			& \mbox{(distributivity)}
			\\ 
			\mbox{NT8} & 
			\forall^1 S \, (S(1) \wedge \forall x \, (S(x) \rightarrow S(x+1)))
			\, \rightarrow\, \forall x \, S(x)
			& \mbox{(induction)}
		\end{array}
	\] 
	The first seven axioms are all first-order, the Induction Axiom, however, is
	second-order, since it quantifies over a \emph{relation}\index{relation}:
	the expression ``$\forall^1 S$'' is read ``for all unary relations $S$.''
	In effect, the Induction Axiom says that any subset of the universe $\mathbb{N}$
	which contains 1 and is closed under the function $x+1$ necessarily contains
	\emph{every} element in the universe.
\end{beispiel}

\begin{satz}
	\emph{(Soundness and Completeness)}
	For any sentence $\phi$ and any set $\Delta$ of sentences,
	$\Delta \vdash \phi$ if and only if $\Delta \models \phi$.
\end{satz}
\begin{proof}
	\cite[Thm.\,]{Hedman-2004},
	\cite[Thm.\,5.6,\,5.7]{Papadimitriou-1994}
\end{proof}

A proof system with the property 
$\Delta \vdash \phi$ $\Rightarrow$ $\Delta \models \phi$
in the above theorem is called 
\emph{sound}\index{sound}, i.e., the system only proves valid consequences.
A system with the reverse property
$\Delta \models \phi$ $\Rightarrow$ $\Delta \vdash \phi$
is called \emph{complete}\index{complete}
and means that it is able to prove all valid consequences.
The fact that first-order logic is complete was first shown by Gödel
in his Completeness Theorem (1930).

Let be $\Delta_N$ be the set of first-order sentences from the axiomatization
in Example \ref{bsp-axiomatization-number-theory}, with the Induction Axiom NT8
replaced by the set $\mathrm{NT8}_1 = \{\psi_\phi\}$
of $\Sigma_N$-sentences $\psi_\phi$ which, for each $\Sigma_N$-formula in one free variable,
are defined by
\begin{equation}
	\psi_\phi \ = \ 
	( \phi(1) \wedge \forall x\, (\phi(x) \to \phi(x+1))
	\to \forall x \, \phi(x)).
\end{equation}
NT8$_1$ is thus a first-order ``approximazation'' of the Induction Axiom.

\begin{satz}
	\emph{\textbf{(Gödel's First Incompleteness Theorem 1931)}}
	If $T$ is a decidable theory containing $\Delta_N$, then $T$ is incomplete.
\end{satz}
\begin{proof}
	\cite[§8.3]{Hedman-2004}.
\end{proof}

\subsection{Modal logic}
Another extension of propositional logic, besides first-order logic, is
\emph{modal logic}\index{modal logic}
\cite[§27]{Boolos-et-al-2002}.
A modal logic is a logic for handling modalities, i.e., 
concepts like ``possibility'' and ``necessity.''
Formally, a modal logic is a propositional logic established with the
modal operator $\Box$ of \emph{necessity}, satisfying
the following axioms.
\begin{itemize}
	\item[\textbf{N}:]
	\emph{(Necessitation rule)}
	If $A$ is true in propositional logic, then $\Box A$ is also true.
	
	\item[\textbf{K}:]
	\emph{(Distribution Axiom)}
	If $\Box (A \to B)$ then $\Box A \to \Box B$.
\end{itemize}
Usually, a further operator $\Diamond$ of \emph{possibility} defined by the relation
\begin{equation}
	\Diamond A = \neg \Box \neg A,
\end{equation}
meaning that ``it is not necessarily true that not-$A$ is true,''
or shortly ``$A$ is possibly true.''
In fact, the modal operators $\Box$ and $\Diamond$ are dual to one another,
\begin{equation}
	\Box = \neg \Diamond \neg,
	\qquad
	\Diamond = \neg \Box \neg.
	\label{eq-dual-square-lozenge}
\end{equation}
The two axioms \textbf{N} and \textbf{K} yield the weakest modal logic $K$, 
invented by Kripke, there are stronger modal logics supposing more axioms.
\begin{itemize}
	\item[\textbf{T}:]
	\emph{(Reflexivity Axiom)}
	$\Box A \ \to \ A$.
	
	\item[\textbf{4}:]
	\emph{(Non-contingent necessity)}
	$\Box A \ \to \ \Box \Box A$.

	\item[\textbf{B}:]
	\emph{(Symmetry)}
	$A \ \to \ \Box \Diamond A$.

	\item[\textbf{D}:]
	\emph{(Deontity)}
	$\Box A \ \to \ \Diamond A$.
\end{itemize}
Depending on the supposed axioms, these are defined the following modal logics:
\begin{itemize}
	\item
	$K: \ \mathbf{K} + \mathbf{N}$.
	
	\item
	$T: \ \mathbf{K} + \mathbf{N} + \mathbf{T}$.

	\item
	$S4: \mathbf{K} + \mathbf{N} + \mathbf{T} + \mathbf{4}$.

	\item
	$S5: \mathbf{K} + \mathbf{N} + \mathbf{T} + \mathbf{4} + \mathbf{B}$.

	\item
	$D: \ \mathbf{K} + \mathbf{N} + \mathbf{D}$.	
\end{itemize}

The historical roots of modal logic go back to Diodorus' problem concerning the question:
``Will there be a sea battle tomorrow?'' 
According to this question, two propositions are possible answers,
$A$ $=$ ``There will be a sea battle tomorrow'' 
or $\neg A$ $=$ ``There will not be a sea battle tomorrow.''
In propositional logic, $A$ is neither true nor false, in
modal logic we have $\Diamond A = \Diamond \neg A = 1$
(both $A$ and $\neg A$ are possibly true), or by Eq.~(\ref{eq-dual-square-lozenge}),
$\Box \neg A = \Box A = 0$
(both $A$ and $\neg A$ are not necessarily true).

In substructural logics there are also considered propositional structures with
even several modal negations, such as an ``n-type negative pair'' 
of operators or ``p-type negative pair'' operators 
\cite[§8.1]{Restall-2000}.

\section{Discussion}
In this paper the relations of classical and non-classical
concepts of logic are reviewed and unified algebraically.
The systematic subjunction of algebraic restrictions to a lattice, 
viz., fuzzy negation, paraconsistency, non-contradiction,
orthomodularity and distributivity, are shown to yield an algebraic hierarchy tree
of logics as in Figure \ref{fig-logics}.
An essential step to this unifying picture is the
notion of a general fuzzy negation holding the conditions of weak double negation
and antitony, as well as the Boolean boundary condition.
Remarkably, the antitony is equivalent to the \emph{disjunctive} De Morgan law
(Theorem \ref{theo-negation-antitony}),
but does not guarantee the conjunctive De Morgan law.
Thus, in contrast to usual hitherto existing approaches to logic, 
the two De Morgan laws do not play a symmetric role.
Example \ref{bsp-M5-non-conjunctive} shows
the non-contradictory logic $(M_5,\sim)$ in which the conjunctive De Morgan law
in fact does not hold.

With this notion of a fuzzy logic it is possible to prove the
validity of the conjunctive De Morgan law for a strong double negation
satisfying $x=x''$ for all propositions $x \in L$
(Theorem \ref{theo-De-Morgan-fuzzy}).
If the logic is non-contradictory, a strong double negation $x=x''$
even implies the rule \emph{tertium non datur}
(Theorem \ref{theo-tertium-non-datur}), 
a result which has been well-known for implicative logics
(Theorem \ref{theo-implicative-strong-double-negation}).

Thus,
fuzzy logics naturally may be contradictory 
and non-orthomodular (Example \ref{bsp-fuzzy-temperature}),
whereas logics of quantum registers or composite spin-$\frac12$ systems 
are orthomodular, but typically non-distributive.
On the other hand, implicative logics, comprising intuitionistic and Boolean logics,
are distributive, although there exist infinite distributive but non-implicative lattices.
Notably, this concept of logic rules out generalized approaches relaxing
the lattice requirement on a logic to
an orthoalgebra \cite{Dalla-Chiara-1995,Foulis-et-al-1992,Wilce-2005}
requirement, 
yielding ``operational logics'' or ``quantum temporal logics'' 
\cite{Isham-Linden-1994,Omnes-1995,Wilce-2005})
of propositions about ``consistent histories'' 
\cite{Gell-Mann-Hartle-1990}.
In fact, such structures should not be called
``logics''
since, technically, in non-lattice sets least upper or lower bounds 
need not exist such that the De Morgan laws are not expressible,
and philosophically, \emph{logical} relations should be properly distinguished
from \emph{causal} relations \cite[p.\,39]{Joos-et-al-2003}.

What are the benefits of the novel approach presented in this paper,
besides a possibly compelling unified view onto the zoo of logics?
One open question in the foundations of quantum theory is how
quantum logics for composite quantum systems are derived systematically.
Taking into account non-orthomodular fuzzy logics
and perhaps some relaxations which render them
into effect algebras
\cite{%
	Bennett-Foulis-1997,%
	Foulis-2006,
	Foulis-Bennett-1994,
	Schroeck-2005%
}
might point a new way to this outstanding problem.



\end{document}